\documentclass[preprint,3p,11pt,sort&compress]{elsarticle}

\usepackage{hyperref}
\hypersetup{backref, colorlinks=true, linkcolor=blue}
\usepackage{bm}
\usepackage{amssymb,amsmath}
\usepackage{amsthm}
\usepackage{stmaryrd}
\usepackage{stfloats}
\usepackage{subcaption}
\usepackage{placeins}
\usepackage{xcolor}

\usepackage{tikz,float}
\usetikzlibrary{shapes,arrows}
\usetikzlibrary{backgrounds,fit}
\usetikzlibrary{patterns}
\usetikzlibrary{decorations.pathreplacing}
\usetikzlibrary{shapes.multipart}
\usetikzlibrary{positioning}
\usetikzlibrary{shapes,arrows}
\usetikzlibrary{backgrounds,fit}

\usepackage{ulem}

\begin{document}

\begin{center}
	{\bf \Large 	Hybrid Weight Window Method for Global Time-Dependent  \\  \medskip Monte Carlo Particle Transport Calculations}
\end{center}

\author[ncsu,cement]{Caleb A. Shaw}
\ead{cashaw4@ncsu.edu}
\author[ncsu,cement]{Dmitriy Y. Anistratov}
\ead{anistratov@ncsu.edu}
\address[ncsu]{Department of Nuclear Engineering,
	North Carolina State University Raleigh, NC}
\address[cement]{Center for Exascale Monte Carlo Neutron Transport (CEMeNT)}

\begin{frontmatter}
\begin{abstract}
This paper presents a new Monte Carlo (MC) algorithm for time-dependent particle transport problems with global variance reduction based on automatic weight windows (WWs). 
The centers of WWs at a time step are defined by the solution of 
an auxiliary  hybrid MC / deterministic problem formed
 by the low-order second-moment (LOSM) equations.
 The closures for the hybrid LOSM equations are calculated by the MC method.
 The  LOSM equations are  discretized by a scheme of the second-order accuracy in time and space.
 Filtering techniques are applied to reduce noise effects in the LOSM closures.
 The WWs defined with the auxiliary solution give rise to sufficiently uniform 
 MC particle distribution in space on each time step.
 The algorithm is analyzed by means of an analytic transport benchmark.
 We  study performance of the MC algorithm depending on a set parameters of WWs.
 Figure of merit and relative error results are presented, demonstrating the performance of the hybrid MC method and quantifying its computational efficiency.
\end{abstract}

\end{frontmatter}

\section{Introduction}

In the field of particle transport, Monte Carlo methods are considered highly desirable for several reasons, such as enabling continuous energy physics and having no discretization error. They achieve this by simulating computational particles that represent physical particles. Some applications of Monte Carlo transport are
 reactor physics, high-energy-density physics, radiation shielding, atmospheric modeling, astrophysics and medical physics. One limitation is that Monte Carlo solutions converge slowly with increasing number of particles, and they have inherent stochastic error. 
    
    By default the distribution of particles in space is unbalanced, with regions near the source having more particles. Regions where physical particles are less common, such as shielded regions or wave fronts end up undersampled. In general, the amount of stochastic error in a Monte Carlo quantity depends on the number of samples that contributed to it. The non-uniform distribution of particles contributes to the non-uniform spatial distribution of the stochastic error. 
    
    To modify the particle distribution, the computational particles are assigned weights. With weights, computational particles no longer correspond one to one with physical particles. Computational particles in high population regions can be eliminated with a process called rouletting. Conversely, particles that do happen to survive to low population regions can be split into any number of daughter particles. Splitting and rouletting techniques don't bias Monte Carlo solution.

    There are several ways to modify the particle weights and reduce variance in the solution. One option, implicit capture, decreases the weight of particles that would otherwise be absorbed, allowing them to survive through shielding. However, this can lead to many particles with extremely low weight, which will decrease the solution efficiency. Alternatively, splitting and rouletting can be used to keep particle weights constrained within some window about a chosen distribution with the goal of reducing variance \cite{Book_MonteCarlo_1991_Lux_MonteCarloParticleTransportMethodsNeutronandPhotonCalculations,TechReport_2024_LANL_LA-UR-24-24602Rev.1_KuleszaAdamsEtAl}.
    
    The choice of what function to use to define centers of weight windows depends on the type of problem being solved. For problems where a detector response functional is desired, the optimal function is the adjoint flux solution.
     Using the adjoint solution redistributes computational particles toward the detector, decreasing variance at the location of the detector.
     The adjoint solution can be computed using a deterministic solver, leading to a hybrid method \cite{osti_527548,wagner-haghighat}.   
    
    In global problems, the solution quality over the whole domain is of equal importance.
    In this class of problems, the adjoint approach requires two deterministic solves, one forward solve that is used to define the adjoint source and one adjoint solve to obtain the importance function \cite{wagner-peplow-mosher-2014}. Using a normalized forward flux solution to define weight windows has been shown to reduce global variance in steady state neutron transport problems \cite{Cooper_Larsen_2001,Cooper_diss_1999}.
    This work used the hybrid low-order quasi-diffusion (LOQD) (aka variable Eddington factor (VEF)) equations for the angular moments of the angular flux
    with the Eddington factor defining an exact closure \cite{gol'din-cmmp-1964, auer-mihalas-1970}.
    The Eddington factor is computed by Monte Carlo.
    
   In time dependent thermal radiative transfer (TRT) problems, the forward flux LOQD solution has been used to define hybrid weight windows (HWW) \cite{ Wollaber_2008}. For these HWW, 
  the LOQD (VEF) equations were discretized in time with a scheme of the first order of accuracy.
   Multiphysics methods that are second order in time have also been developed for TRT problems \cite{mcclarren2013temperature}.
    The Second Moment (SM) method has been used to develop hybrid implicit Monte Carlo (IMC) schemes for TRT \cite{pozulp-mc2023,pozulp-mc2025}.
    The low-order SM (LOSM) equations are formulated for the scalar flux and current  with exact linear closures \cite{smm-1976}.

    In this paper, we present a new automatic HWW technique  for global variance reduction in  time dependent particle transport problems. At each time step, the centers of weight windows are defined by the
      numerical solution of a hybrid Monte Carlo / deterministic problem
     based on the LOSM equations.
       A hybrid SM method was shown to be effective with low particle counts per cell in steady-state transport problems \cite{vnn-dya-ans-annual-2024, vnn-dya-arxiv-2025}. Additionally, the differential operator of the low-order second moment (LOSM) equations is self-adjoint, which is advantageous when solving multi-dimensional problems.
      The hybrid LOSM equations are discretized in space by second-order finite volume scheme.
      The Crank-Nicolson time-integration method of the second-order accuracy is used for temporal  approximation of the LOSM equations.
      The closures for the LOSM equations are computed by the Monte Carlo solution.
    One version of the developed method applied filtering techniques to the closures of the LOSM equations and their initial condition at each time step to mitigate the effect of
     stochastic noise. The moving average filter and Fourier filtering were both considered.
   We conducted a study of parameters which define the weight windows algorithm  and used the 
   results  to motivate the choices in our Monte Carlo calculations. 
     The performance of the proposed Monte Carlo algorithms with automatic weight windows are analyzed by means of an analytic transport benchmark.
        The developed weight windows  algorithms have been implemented in the   Monte Carlo / Dynamic Code (MC/DC) \cite{Joss}.

     The rest of the paper is organized as follows. 
    Section \ref{sec:ww} describes a general weight windows method
    for global time-dependent transport problems. 
    Section \ref{sec:hybrid-ww-method} presents the weight window method for dynamic Monte Carlo calculations based on auxiliary solution of hybrid LOSM equations.
    The noise filtering techniques are described in Section \ref{sec:filtering}.
    In Section \ref{sec:numerical_results} are numerical results for the parameter study as well as an analysis of the performance of HWW.
     Section \ref{sec:discussion} concludes with a summary of the findings and directions for future work.

\section{Weight Windows for  Global Transport Problems \label{sec:ww}}

In non-analog Monte Carlo  methods, particles are assigned  a weight
that can be changed to reduce variance. 
Weight windows are a way of controlling both the weights, and the spatial distribution of particles. In this scheme for global variance reduction, the weights are altered to obtain a more uniform particle distribution \cite{Cooper_Larsen_2001,Wollaber_2008}. The windows are applied 
when a particle undergoes an event or enters a new mesh cell.
The main parameters of a weight window are its center, ceiling, and floor. The center of the weight window in the $i^{th}$ spatial cell is usually defined  as follows \cite{Cooper_Larsen_2001}:
\begin{equation}
	ww_i^n = \frac{ {\tilde{\phi}_i^n}}{\max\limits_i\tilde{\phi}_i^n},
\end{equation}
where $n$ is the index of the time step of dynamic Monte Carlo, 
$ {\tilde{\phi}_i^n}$  is  the approximate scalar flux in the $i^{th}$ cell at time  $t^n$ 
evaluated by some numerical technique.
The ceiling and floor of the weight window are given by
\begin{subequations} \label{ww-c&f}
\begin{equation} \label{ww-ceiling}
	ww_i^{n,ceiling} = ww_i^n \times \rho \, ,
\end{equation}
\begin{equation} \label{ww-floor}
	ww_i^{n,floor} = \frac{ww_i^n}{ \rho } \, ,
\end{equation}
\end{subequations}
where $\rho \geq 1$ is the window width parameter.

Let us consider a particle having a weight $w_P$.
If the particle is in the $i^{th}$ cell during the $n^{th}$ time interval and $w_P>ww_i^{n,ceiling}$,
the particle will be split into $N_d$ daughter particles, where $N_d = \lceil \frac{w_P}{ww_{i}^{n}}\rceil$. 
The original particles weight is distributed onto its daughter particles.
 If $w_P<ww_i^{n,floor}$, then the particle will be rouletted with survival probability $p =\frac{w_P}{ww_{i}^{n}}$ and if it survives, it will be assigned $w_P = ww_i^n$. Combined, these processes ensure that the weights of all particles in the cell fall within the window of that cell. Note that this process preserves the total weight of all particles and does not introduce any bias into the solution.

Since particles with weights higher than the ceiling are split, it is possible for cells with very low ceilings to experience excessive splitting. \ This also occurs in regions of the domain where the solution has
a large gradient, such as wave fronts. When the window centers, $ww_{i}^{n}$, approach zero, the number of daughter particles, $\lceil\frac{w_P}{ww_{i}^{n}}\rceil$, grows by orders of magnitude. 
Time dependent weight windows can cause excessive splitting \cite{Landman_McClarren_Madsen_Long_2014,Wollaber_2008}. 
 In thermal radiative transfer, where waves are common feature, several modifications have been developed to artificially raise the window center near the wave front \cite{Wollaber_2008,Cooper_diss_1999}. 
   To limit excessive splitting in this work, we define a modified window center given by 
  \begin{equation} \label{ww-center}
	ww_i^n = \left[\frac{ \tilde{{\phi}}_i^n}{\max\limits_i\tilde{\phi}_i^n}\right]\times(1-\varepsilon_{min})+\varepsilon_{min} \, ,
\end{equation}
where $\varepsilon_{min}$ is a parameter that sets the minimum value of the window centers. Using this modification,
  the shape of the original window centers is preserved and scaled to fit above the minimum. 

\section{Hybrid Monte Carlo / Deterministic Method for Time-Dependent Transport Problem  \label{sec:hybrid-ww-method}}

The governing equation for the dynamic behavior of the ensemble of particles is the neutron transport equation. In this work, we focus on problems with one group in energy, in one-dimensional slab geometry. Under these conditions, and assuming isotropic scattering  and source, the time dependent neutron transport equation is given by
\begin{multline}  
    \frac{1}{v}\frac{\partial\psi}{\partial t}(x,\mu,t)+\mu\frac{\partial \psi}{\partial x}(x,\mu,t)+\Sigma_t(x,t)\psi(x,\mu,t)=\\ 
    \frac{1}{2}\Big((\Sigma_s(x,t)+ 
    \nu_f \Sigma_f(x,t) \big) \int_{-1}^1 \psi(x,\mu',t) d \mu' +q(x,t)\Big),
    \label{eqn:HO_transport}
\end{multline}
\[
x\in[0,X],\text{ } \mu\in[-1,1],\text{ } t\geq 0.
\]
This equation is for $\psi$, the particle angular flux. The neutron speed is $v$ and isotropic external source is $q$. The material properties are $\Sigma_t$, $\Sigma_s$, $\Sigma_f$, the total, scattering, and fission cross sections, respectively. $\nu_f$ is the expected number of neutrons due to fission. $x,t$, and $\mu$ are the spatial position, time and directional cosine of particle motion.
The boundary conditions (BCs) are
\begin{equation} \label{t-eq-bc}
    \psi(0,\mu,t) = \psi_L(\mu,t),\text{ for } \mu> 0, \quad
    \psi(X,\mu,t) = \psi_R(\mu,t),\text{ for } \mu< 0, 
\end{equation}
and the initial condition (IC) is
\begin{equation} \label{t-eq-ic}
    \psi(x,\mu,0)  = \psi^{ini}(x,\mu)\, .
\end{equation}

To formulate a hybrid Monte Carlo / deterministic method we apply low-order equations of the second moment method \cite{smm-1976,Cefus-Larsen-ttsp-1989}. 
The low-order second moment (LOSM) equations are derived by  integrating  the transport equation   (Eq. \eqref{eqn:HO_transport})  with weight 1 and  $\mu$
  over $-1\leq\mu \leq 1$.
The LOSM equations for the scalar flux, $\phi(x,t) = \int_{-1}^{1} \psi (x,\mu,t) d \mu$, and current, $J(x,t) = \int_{-1}^{1} \mu \psi (x,\mu,t) d \mu$ are defined by
\begin{subequations} \label{losm}
\begin{equation} 
    \frac{1}{v}\frac{\partial \phi}{\partial t}(x,t) + \frac{\partial J}{\partial x}(x,t) + \big(\Sigma_t(x,t)-\Sigma_s(x,t)-\nu_f\Sigma_f(x,t)\big)\phi(x,t) = q(x,t),
    \label{eqn:balance}
\end{equation}
\begin{equation}
    \frac{1}{v}\frac{\partial J}{\partial t}(x,t) + \frac{1}{3}\frac{\partial \phi}{\partial x}(x,t) +\Sigma_t(x,t)J(x,t) = \frac{\partial F}{\partial x}(x,t),
    \label{eqn:moment}
\end{equation}
\end{subequations}
with the associated BCs
\begin{equation} 
    J(0,t) = -\frac{1}{2}\phi(0,t)+2J_L(t)+P_L(t),\quad J(X,t) = \frac{1}{2}\phi(0,t)+2J_R(t)-P_R(t),
     \label{losm-bcs}
\end{equation}
and ICs
\begin{equation} 
    \phi(x,0) = \phi^{ini}(x), \quad J(x,0) = J^{ini}(x) \, ,
    \label{losm-ics}
\end{equation}
where
\begin{equation}
	J_L = \int_0^1 \mu \psi_L d \mu \, , \quad  J_R = \int_{-1}^0 \mu \psi_R d \mu \, , 
\end{equation}
\begin{equation}
	\phi^{ini}= \int_{-1}^1  \psi^{ini} d \mu \, , \quad  	J^{ini}= \int_{-1}^1 \mu  \psi^{ini} d \mu \, .
\end{equation}
The closures for the system of high-order transport and LOSM equations  are defined by
\begin{equation} \label{F}
   F(x,t) = \int_{-1}^{1}\left(\frac{1}{3}-\mu^2\right)\psi(x,\mu,t)d\mu \, ,
\end{equation}
and 
\begin{equation*}
   P_L(t) = \int_{-1}^{1}\left(\frac{1}{2}-|\mu|\right)\psi(0,\mu,t)d\mu,\quad
   P_R(t) = \int_{-1}^{1}\left(\frac{1}{2}-|\mu|\right)\psi(X,\mu,t)d\mu \, .
\end{equation*}

To discretize the LOSM equations  (Eq. \eqref{losm}) in time, we apply the $\theta$-weighted time-integration method to obtain
\begin{subequations} \label{theta-scheme}
\begin{multline}
	\frac{1}{v\Delta t^n}\left(\phi^n-\phi^{n-1}\right) + \theta \bigg[ \frac{dJ^n}{dx} + \big ( \Sigma_t^n-\Sigma_s^n-\nu_f^n\Sigma_f^n \big)\phi^n \bigg] =\\
	\theta q^n + (\theta-1) \bigg[ \frac{dJ^{n-1}}{dx}+(\Sigma_t^{n-1}-\Sigma_s^{n-1}-\nu_f^{n-1}\Sigma_f^{n-1} )\phi^{n-1} - q^{n-1} \bigg] \, ,
	\label{eqn:LOSM-theta-1}
\end{multline}
\begin{multline}
	\frac{1}{v\Delta t^n} \left(J^n-J^{n-1}\right)+\theta \bigg[ \frac{1}{3}\frac{d\phi^n}{dx}+\Sigma_t^nJ^n \bigg] =\\ 
	\theta \frac{dF^n}{dx} + (\theta -1) \bigg[ \frac{1}{3}\frac{d\phi^{n-1}}{dx}+\Sigma_t^{n-1}J^{n-1} - \frac{dF^{n-1}}{dx} \bigg]\, ,
	\label{eqn:LOSM-theta-2}
\end{multline}
\end{subequations}
where $\phi^n = \phi(x,t^n)$, $J^n=J(x,t^n)$, $n$ is the index of the time layer, and $\Delta t^n =t^n - t^{n-1}$ is the $n^{th}$ time step.
The time-integration  scheme \eqref{theta-scheme} with $\theta =1$ is the Backward Euler (BE)  method of the first-order accuracy. For $\theta=\frac{1}{2}$, it reduces to the second-order 
Crank-Nicolson (CN) method.

	The  LOSM equations  are discretized in space by a second-order finite volume scheme.
	A spatial mesh $\{x_{i}\}_{i=0}^{I}$ with $I$ intervals is defined on the domain. 
    We assume that cross sections and source are piece-wise constant on the set of mesh cells.
	The balance equation 	(Eq. \eqref{eqn:LOSM-theta-1}) is integrated over the $i^{th}$ spatial cell $[x_{i-1}, x_{i}]$.
 	The first-moment equation (Eq. \eqref{eqn:LOSM-theta-2})  is integrated over half cells.
 	The discretized LOSM equations are given by
 	\begin{subequations} \label{losm-disc}
 		 	\begin{multline}  \label{losm-0}
 			\frac{1}{v\Delta t^n}\big( \phi_i^n -\phi_i^{n-1} \big) \Delta x_i	 + 
 			\theta \bigg[ J_{i}^n -J_{i-1}^n + ( \Sigma_{t,i}^n-\Sigma_{s,i}^n- \nu_{f,i}^n\Sigma_{f,i}^n )\Delta x_i \phi_i^n \bigg] = \\
 			\theta q_i^n + (\theta-1)\bigg[ J_{i}^{n-1} - J_{i-1}^{n-1}  + ( \Sigma_{t,i}^{n-1} -\Sigma_{s,i}^{n-1} - \nu_{f,i}^{n-1} \Sigma_{f,i}^{n-1}  )\Delta x_i\phi_i^{n-1}  - q_i^{n-1} \bigg] \, ,   \\
            i=1,\ldots,I \, ,
 		\end{multline}
 		\begin{multline} \label{losm-1}
 			\frac{1}{v\Delta t^n}\big( J_{i}^n  -J_{i}^{n-1} \big)   \Delta \hat  x_{i} +
 			\theta \bigg[  \frac{1}{3} \Big( \phi_{i+1}^n - \phi_{i}^n \Big)   + 
 			\hat \Sigma_{t,i}^n  \Delta \hat  x_{i} J_{i}^n \bigg] = \\
 			\theta \big( F_{i+1}^n-F_{i}^n \big) + 
 			(\theta-1) \bigg[  \frac{1}{3} \Big( \phi_{i}^{n-1}- \phi_{i+1}^{n-1} \Big)   + 
 			\hat \Sigma_{t,i}^{n-1} \Delta \hat  x_{i} J_{i}^{n-1}  -   F_{i+1}^{n-1} + F_{i}^{n-1}  \bigg] \, , \\
            i=0,\ldots,I
 		\end{multline} 
 		\begin{equation}  \label{losm-disc-bc}
 			J_{0}^n = -\frac{1}{2}\phi_0^n+2J_L^n+P_L^n ,\quad J_{I}^n = \frac{1}{2}\phi_{I+1}^n + 2J_R^n-P_R^n,
 		\end{equation}
 		\begin{equation*}
 			\hat \Sigma_{t,i}^n =\frac{  \Sigma_{t,i}^n \Delta x_{i} +  \Sigma_{t,i-1}^n \Delta x_{i-1} }{ \Delta x_{i} + \Delta x_{i-1} } \, ,  \quad 
              			\Delta \hat x_{i} = \frac{1}{2}( \Delta x_{i+1} +  \Delta x_{i}) \, , 
 		\end{equation*}
 		\begin{equation*}
 			\Delta x_{i} =    x_{i} - x_{i-1} \, , \quad
        i=1,\ldots, I \, , \quad
       \Delta x_{0} = \Delta x_{I+1} = 0 \, ,
 		\end{equation*}
	\end{subequations}    
 where $J_{i}^n = J^n(x_i)$ is the cell-edge current,
\begin{equation}
	\phi_i^n =\frac{1}{\Delta x_i} \int_{x_{i-1}}^{x_{i}}\phi^n(x)dx  \, , 
\end{equation}
\begin{equation} \label{F_i-n}
	F_i^n= \frac{1}{\Delta x_i} \int_{x_{i-1}}^{x_{i}}F^n(x)dx \,  
\end{equation}
are cell-average values of the scalar flux and closure function ($i=1,\ldots, I$), respectively, and 
their values at domain boundaries are $\phi_0^n= \phi^n(x_0)$, $F_0^n= F^n(x_0)$, 
$\phi_{I+1}^n= \phi^n(x_{I+1})$, and $F_{I+1}^n= F^n(x_{I+1})$.

	The closure terms of the discretized LOSM equations defined by $\{F_i^n\}_{i=0}^{I}$, $P_L^n$, and $P_R^n$
	are computed by means of the Monte Carlo solution during the $n^{th}$ time step. 
    Each time step is treated individually. The Monte Carlo solution is also used to calculate (i) the initial conditions for  the $n^{th}$ time step defined by 
	$\{\phi_i^{n-1}\}_{i=0}^{I+1}$  and 	$\{J_{i}^{n-1}\}_{i=0}^{I}$
	and (ii) $\{F_i^{n-1}\}_{i=0}^{I+1}$   in the case of the CN scheme.
    	This defines the hybrid Monte Carlo/ deterministic method based on the LOSM equations for global time-dependent transport calculations.
	We apply the numerical solution of the hybrid LOSM (HLOSM) equations  (Eqs. \eqref{losm-disc})   on each time step as an auxiliary  solution $\{ \tilde \phi_i^n \}_{i=1}^I$ to define centers of the weight windows according to Eq. \eqref{ww-center}.

Computing the initial conditions and closures for the HLOSM problem involves evaluating certain grids functions of    
$\phi$, $J$, and $F$.
The discretized LOSM equations specify that these quantities are evaluated at the discrete time layers. For the initial conditions, a mixture of cell-average and cell-edge quantities is needed at the time layer. 

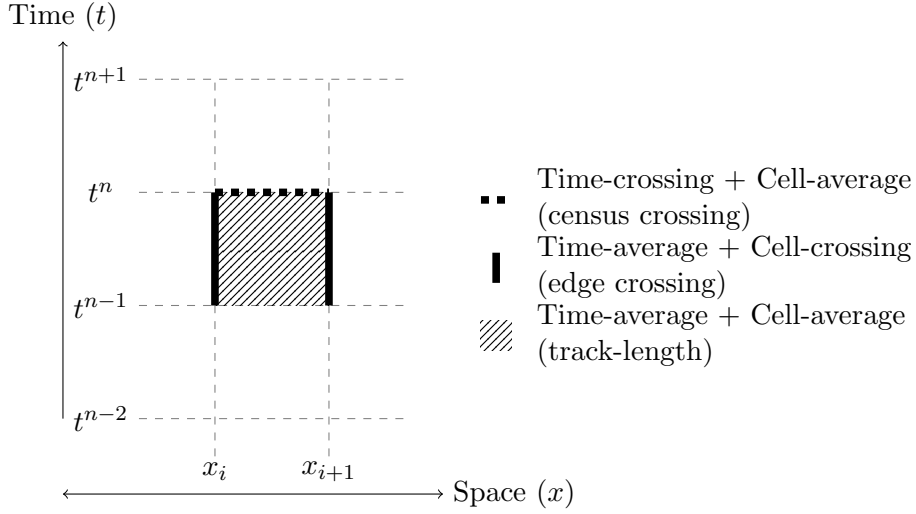
\begin{figure}[h]
    \centering
    \begin{tikzpicture}

        \draw[step=1.5cm, gray, thin,dashed] (-1, -0.5) grid (2.5, 4.5);

        \draw[<->] (-2, -1) -- (3, -1) node[right] {Space ($x$)};
        \draw[->] (-2, 0) -- (-2, 5) node[above] {Time ($t$)};

        \node at (0, -0.7) {$x_{i}$};
        \node at (1.5, -0.7) {$x_{i+1}$};

        \node at (-1.5, 0) {$t^{n-2}$};
        \node at (-1.5, 1.5) {$t^{n-1}$};
        \node at (-1.5, 3) {$t^n$};
        \node at (-1.5, 4.5) {$t^{n+1}$};

        \fill[pattern=north east lines, pattern color=black] (0, 1.5) rectangle (1.5, 3);

        \draw[black,line width=1mm,dashed] (0, 3) -- (1.5, 3) node[above] {};
        
        \draw[black,line width=1mm] (0, 1.5) -- (0.0, 3) node[above] {};
        \draw[black,line width=1mm] (1.5, 1.5) -- (1.5, 3) node[above] {};

        \fill[pattern=north east lines, pattern color=black] (3.5, 0.9) rectangle (3.9, 1.3);
        \node[anchor=west, align=left] at (4.1, 1.1) {Time-average + Cell-average \\ (track-length)};

        \draw[black,line width=1mm,dashed] (3.5, 2.9) -- (3.9, 2.9);
        \node[anchor=west, align=left] at (4.1, 2.9) {Time-crossing + Cell-average \\(census crossing)};
        
        \draw[black,line width=1mm] (3.7, 1.8) -- (3.7, 2.2);

        \node[anchor=west, align=left] at (4.1, 2.0) {Time-average + Cell-crossing\\ (edge crossing)};
    \end{tikzpicture}
    \caption{ \label{fig:tally_locations} Tally options in Monte Carlo space-time cell}

\end{figure}
Figure \ref{fig:tally_locations} shows the various locations in space-time where tallies can be performed. Track-length tallies are used to compute the final Monte Carlo solution 
 
\begin{equation} \label{phi-mc}
\big< \phi \big>_i^n = \frac{1}{\Delta t^n \Delta x_i} \int_{t^{n-1}}^{t^n} \int_{x_i}^{x_{i+1}} \phi dx dt \, .
\end{equation}
as an average over the $i^{th}$ spatial cell and $n^{th}$ time interval. In this type of tally, the paths left by particles are accumulated to estimate integral quantities. The quantities needed for the hybrid problem are not integrated in time, so track-length tallies will not directly provide the needed information. Crossing tallies are used to evaluate quantities that are discrete in either space or time, but not both. This type of tally requires a particle to pass through the discrete point, either in time or space, to contribute. Surface crossings tallies generally have greater statistical error than track-length tallies for the same number of source particles due to crossing events being less frequent. This represents a trade off, since this tally can directly estimate some of the quantities needed for the HLOSM problem. To obtain coefficients of the discretized HLOSM equations,
  time-crossing tallies,
  \begin{equation} \label{ct-mc}
      \llbracket \phi \rrbracket^n_i = \frac{1}{\Delta x_i} \int_{x_i}^{x_{i+1}} \phi(x,t^n) dx \, ,
      \end{equation}
are applied to compute average values in spatial cells at the time layer $t^n$.
Since the likelihood of a particle crossing both a space and time surface at the same time is exceedingly rare, we cannot evaluate the current at the cell edge and time layer directly. 
To compute the current at $x_i$, and $t^{n-1}$,
 we interpolate linearly in space, from  the cell-average currents in the cells adjacent to $x_i$ at the time layer.

To compute weight windows at the beginning of the $n^{th}$ time step when the Monte Carlo solution at $t^n$ is not available yet, 
the HLOSM equations are defined using the closures obtained from the Monte Carlo solution over the previous time step, namely,
$\{F_i^{n-1}\}_{i=0}^{I+1}$,
$P_L^{n-1}$, and $P_R^{n-1}$. At this stage, the weight windows are determined by the solution of semi-implicit HLOSM equations
given by Eq. \eqref{losm-0} and
 	\begin{subequations} \label{losm-disc-semi}
 		\begin{multline} \label{losm-1-semi}
 			\frac{1}{v\Delta t^n}\big( J_{i}^n  -J_{i}^{n-1} \big)   \Delta \hat  x_{i} +
 			\theta \bigg[  \frac{1}{3} \Big( \phi_{i+1}^n - \phi_{i}^n \Big)   + 
 			\hat \Sigma_{t,i}^n  \Delta \hat  x_{i} J_{i}^n \bigg] = \\
 			\theta \big( F_{i+1}^{n-1}-F_{i}^{n-1} \big) + 
 			(\theta-1) \bigg[  \frac{1}{3} \Big( \phi_{i}^{n-1}- \phi_{i+1}^{n-1} \Big)   + 
 			\hat \Sigma_{t,i}^{n-1} \Delta \hat  x_{i} J_{i}^{n-1}  -   F_{i+1}^{n-1} + F_{i}^{n-1}  \bigg] \, , \\
            i=0,\ldots,I
 		\end{multline} 
 		\begin{equation}  \label{losm-disc-bc-semi}
 			J_{0}^n = -\frac{1}{2}\phi_0^n+2J_L^n+P_L^{n-1} ,\quad J_{I}^n = \frac{1}{2}\phi_{I+1}^n + 2J_R^n-P_R^{n-1},
 		\end{equation}
        \end{subequations}

During the $n^{th}$ time step ($t^{n-1}< t \leq t^{n}$),   the weight windows are updated 
 $u_{ww}^n$ 
times by means of the solution to the fully implicit HLOSM equations (Eqs. \eqref{losm-disc})  using closures computed from Monte Carlo particle histories simulated before the update.
	 Let $H^n$ be the total number of source particles on the $n^{th}$ time step. The $p^{th}$ update occurs  
 after the number of histories reached the value given by
\begin{equation} \label{Hpn}
    H_p^n =f_pH^n \, , \quad p=1,\ldots,u_{ww}^n \, ,
\end{equation}
where $f_p$ is the   fraction of the total number of particle histories specified for the $p^{th}$ update.

\section{Statistical Noise Filtering of Closures and Initial Conditions \label{sec:filtering}}

The error in solution of the HLOSM equations  has two components: (1) discretization error due to approximation of the differential equations
and (2) statistical error due to stochastic noise  in  closures and   initial conditions  computed by a Monte Carlo algorithm.
The  effects of noise can be pronounced when these quantities are computed from  small number of particle histories. 
To reduce noise effects in the hybrid solution, we formulate a modified version of the HLOSM equations defined 
 with closures and initial conditions processed by noise filtering methods. 
 We note that   the hybrid solution is used only to define weight windows and hence the solution of the HLOSM method with filtering     does not bias the Monte Carlo solution.
 Two types of filtering are used: (1)  a moving average  (MA) filter which operates in the spatial domain of the problem and
  (2) Fourier filtering operating in the frequency domain. 
  Filtering is applied to  grid functions of closure terms 
  $\{ F_i^{n} \}_{i=1}^{I}$, $\{F_i^{n-1}\}_{i=1}^{I}$,
  and  initial conditions $\{\phi_i^{n-1} \}_{i=1}^{I}$, and $\{J_i^{n-1}\}_{i=0}^{I+1}$,

 The MA filter  is a technique that modifies a  grid function by averaging its values over a  spatial interval
	to reduce high frequency noise \cite{Smith}. 
The  filter is applied to a vector $\{ \chi_m \}_{m=1}^{M}$ representing a grid function over a uniform mesh.
It  is defined by
\begin{equation}
	\label{eqn:moving_average}
	\tilde \chi_m  = \frac{1}{2k+1}\sum_{\ell =-k}^k \chi_{m-\ell} \, .
\end{equation} 
where $\tilde \chi_j $ is the filtered $j^{th}$ element of the grid function and $k \in \mathbb{N}$ is the parameter of the filter base which specifies the spatial interval for averaging.
A filter with a large interval results in more smoothing at the expense of preserving details in the grid function. The resolution can be adjusted by optimizing the filter base. 
Near  boundaries of the spatial domain, the filter base  extends outside of the domain. 
The parameter $k$  is adaptively decreased  to avoid this issue.

Fourier filtering is a method that  applies discrete Fourier transform  (DFT) to  a grid function,  uses a filter the frequency domain, and
transforms the grid function to its original domain. The DFT of $\{ \chi_m \}_{m=1}^{M}$ is 
given by
\begin{equation}
	\label{eqn:fft}
	\mathcal{X}_j = \sum_{m=1}^{M} \chi_m  e^{- i 2 \pi \frac{j}{M} m } \, ,
\end{equation}
The Fourier transform is truncated to remove the high frequency modes using a cutoff value, $\varkappa$,
\begin{equation}
	\tilde{\mathcal{X}}_j = 
	\begin{cases}
		\mathcal{X}_j, & j \leq \varkappa \, , \\
		0, & j > \varkappa \, .
	\end{cases}
\end{equation}
This is also known as a low-pass filter when in the context of signal processing. The truncated frequency representation is then 
transformed
 back into original domain, yielding the filtered 
grid function $\{ \tilde \chi_m \}_{=1}^{I}$, where
\begin{equation}
	\label{eqn:ifft}
        \tilde \chi_m =  \frac{1}{M}\sum_{j=1}^{M} \tilde{\mathcal{X}}_j   e^{i 2 \pi \frac{j}{M} m} \, .
\end{equation}
The  large-scale structure of the solution  is defined by  the lower frequency modes.
 The noise effects are typically represented by the high frequency components that 
are set to zero. As a result, the filtered solution 
  preserves the original structure without high frequency oscillations.

\section{Numerical Results}
\label{sec:numerical_results}

In this section, we present numerical results for various aspects of the hybrid problem. First, a description of the test is given. Hybrid solutions computed with CN and BE time discretizations are shown, illustrating the advantage of the higher order CN discretization. The results of the parameter study in $u_{ww}^n$, $\epsilon_{min}$, and $\rho$ are then presented, motivating the choice of parameters. Following that is an analysis of the quality of the hybrid LOSM solution, and comparison to the hybrid Monte Carlo solution. The effects of filtering techniques on the closure are also presented here. Finally, a comparison in terms of particle track density, relative variance, and figure of merit between analog Monte Carlo, and hybrid weight windows is presented.
\subsection{Test Problem}
To analyze the proposed hybrid weight window method, we use a dynamic neutron transport benchmark problem with a semi-analytic solution \cite{azurv1}.
This is a 1D slab geometry problem for
  infinite homogeneous medium.
The  material cross sections are
$\Sigma_s = \frac{1}{3} \,   cm^{-1}$, $\Sigma_f = \frac{1}{3} \,   cm^{-1}$, $\Sigma_t~=~1~cm^{-1}$, and $\nu_f = 2.3$.
The speed of particles is $v=1  \frac{cm}{s}$.
A point  source is located in the domain center at $x=0$ and  produces a particle impulse at  $t=0$ 
with strength 1 $\frac{particle}{cm^3  \, s}$. 
There are no particles in the domain at $t=0$.
The time interval of the problem is $0\leq t\leq 20s$. 
The spatial mesh consists of 201 uniform cells covering the region $-20.5cm\leq x \leq 20.5cm$. 
The test is calculated with the time step $\Delta t=1$ s.

\begin{figure}[H]
    \centering
    \includegraphics[width=0.49\linewidth]{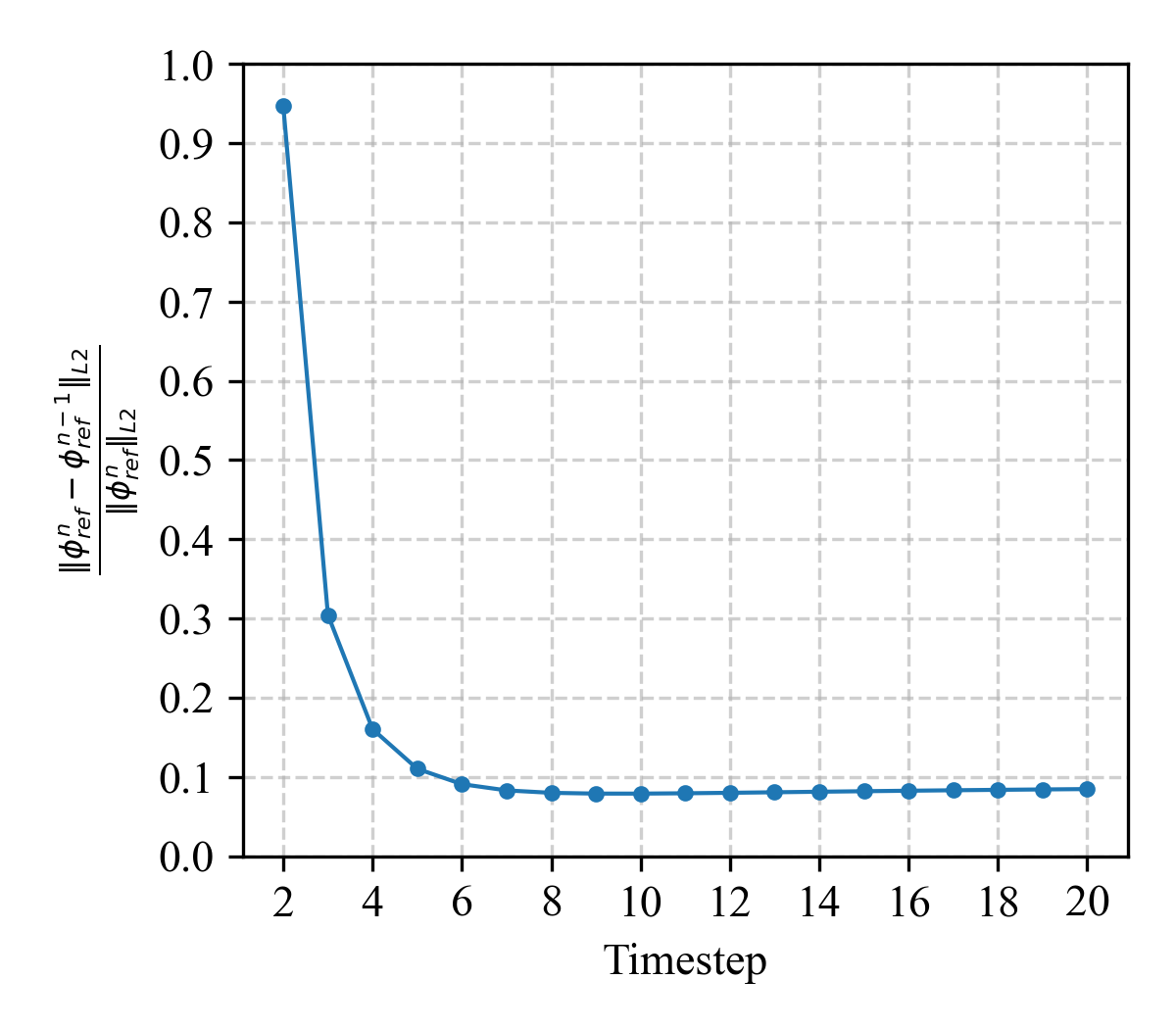}
    \caption{\label{fig:rel_change} Relative solution change}
\end{figure}

This problem consists of a multiplying material and is supercritical. 
The point impulse source  forms a particle wave traveling in the domain.
 The growing population of particles necessitates  a population control technique   to avoid overfilling particle banks. We use the uniform combing technique \cite{osti_1889957}.  The obtained numerical solutions are compared to the benchmark semi-analytic solution.
 The solution in this problem changes significantly in both the shape and magnitude while the initial burst of neutrons spreads out. 
Figure \ref{fig:rel_change} shows 
a relative change rate in the scalar flux versus time, $\alpha^n=\frac{\|\phi^n-\phi^{n-1}\|_{L2}}{\|\phi^n\|_{L2}}$,  computed using the benchmark solution. 
There are two stages in evolution of the solution.
During the initial stage ($n \le 5$),  $\alpha^n$    reduces from 0.9 to 0.1. On the second stage, $\alpha^n$ stabilizes near 0.08.

\subsection{Second Order Time Discretization}
\begin{figure}[H]
    \centering
    \begin{subfigure}[b]{0.45\textwidth}
        \centering
        \includegraphics[width=\textwidth]{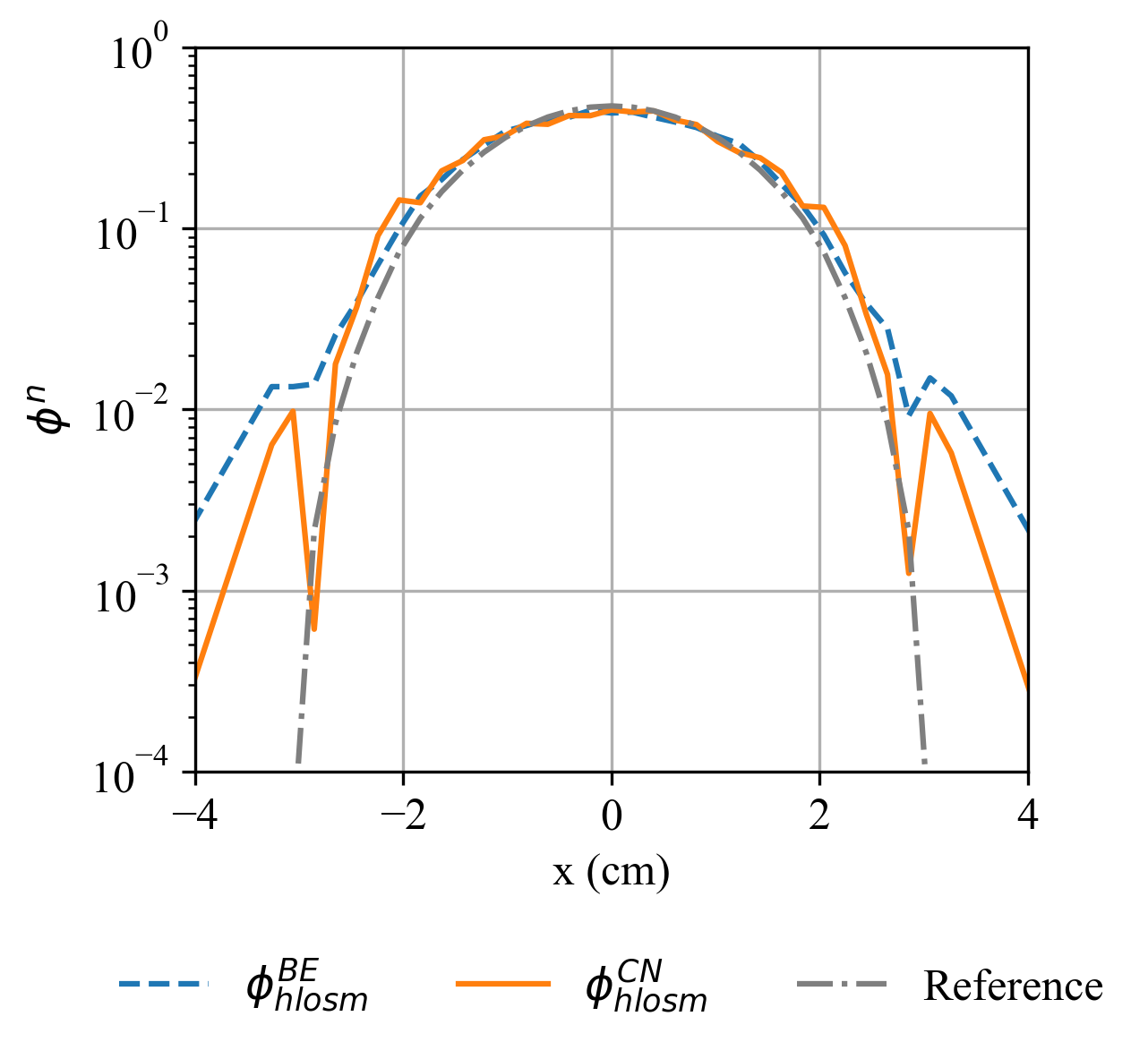}
        \caption{$n=3$}
        \label{fig:aux_be_cn_3}
    \end{subfigure}
    \hfill
    \begin{subfigure}[b]{0.45\textwidth}
        \centering
        \includegraphics[width=\textwidth]{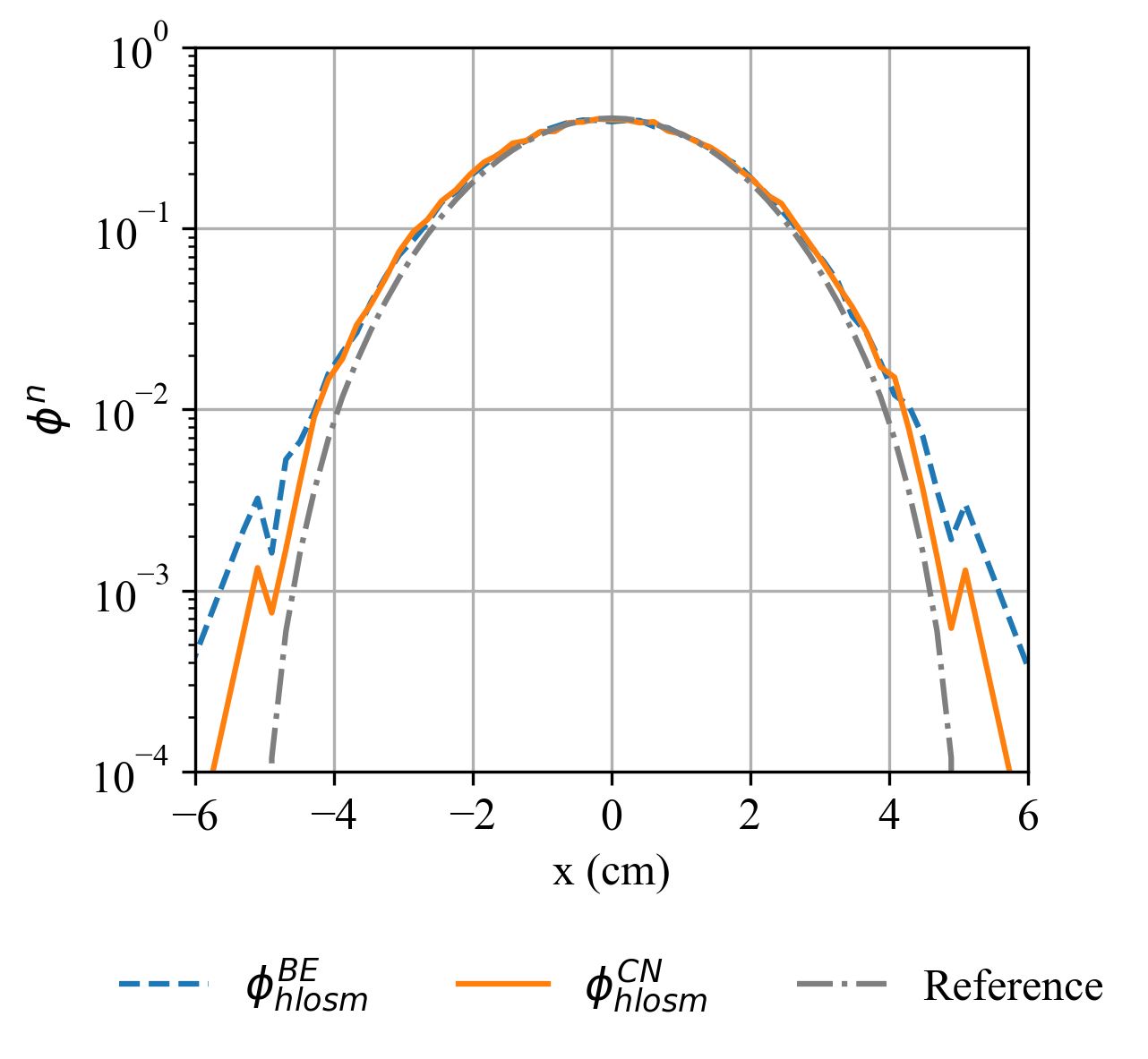}
        \caption{$n=5$}
        \label{fig:aux_be_cn_5}
    \end{subfigure}
    \hfill
    \begin{subfigure}[b]{0.45\textwidth}
        \centering
        \includegraphics[width=\textwidth]{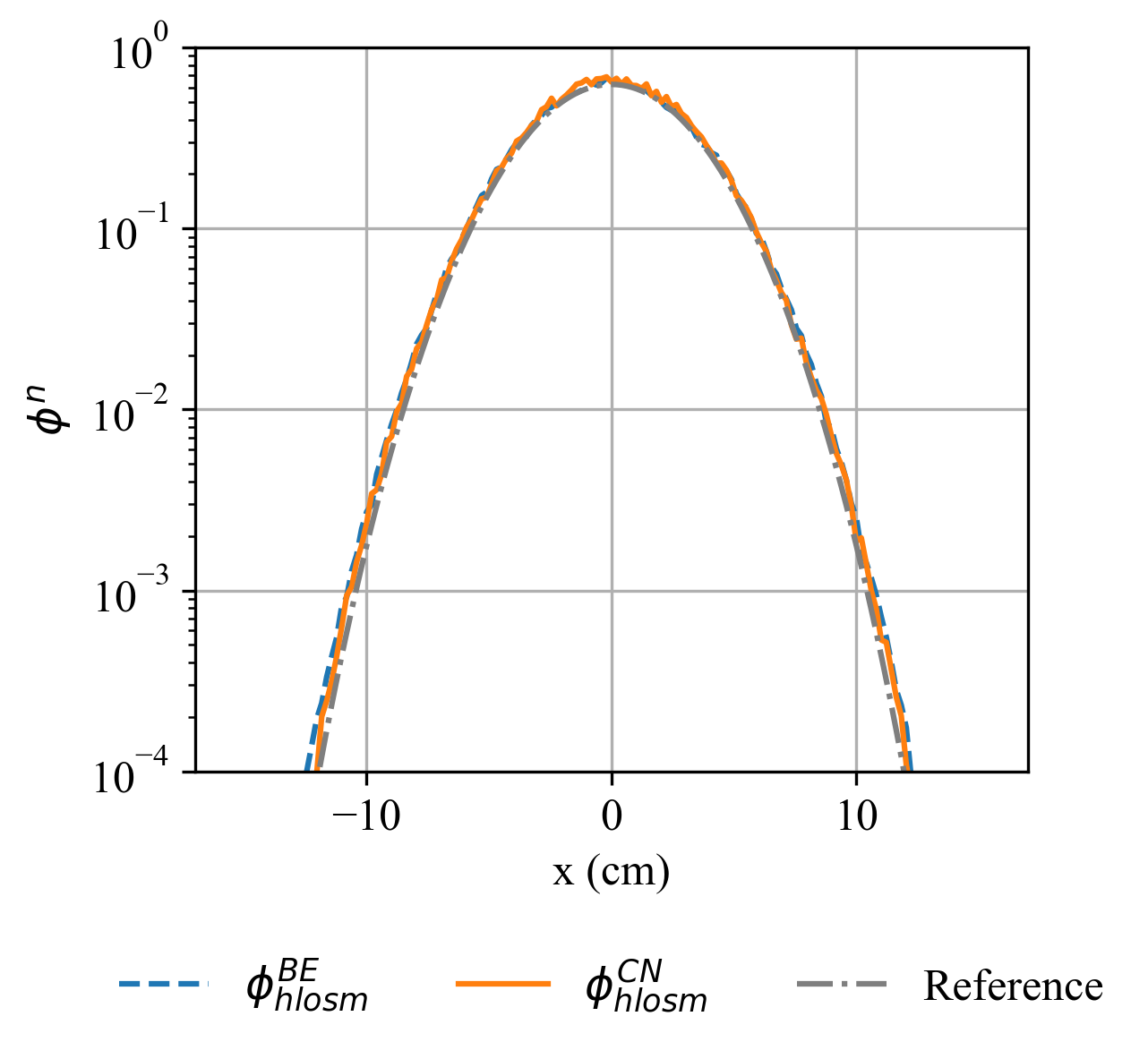}
        \caption{$n=16$}
        \label{fig:aux_be_cn_10}
    \end{subfigure}
    \caption{ Hybrid solutions $ {\phi}_{hlosm}$ computed with the BE and CN time integration schemes and the benchmark solution   at  time steps $n=3,5,10$}
    \label{fig:Prob_1_aux_be_cn}
\end{figure}

In this study, the Monte Carlo algorithm applies weight windows defined by the numerical  solution
of the HLOSM equations, $ {\phi}_{hlosm}$, discretized in time by the CN scheme ($\theta=\frac{1}{2}$). The CN scheme is a high-order time-integration method and is more accurate   than the BE scheme.
Figure \ref{fig:Prob_1_aux_be_cn} 
presents  $ {\phi}_{hlosm}^{BE}$ and $ {\phi}_{hlosm}^{CN}$, the hybrid solutions at time steps $n=3,5,10$ computed by solving the HLOSM equations approximated by BE and CN methods respectively. 
The closures are calculated by Monte Carlo with $H^n=10^4$ particle histories.
  The results demonstrate   that the CN scheme resolves the wave front better than BE scheme.

\pagebreak
\subsection{Analysis of Weight Window Parameters}

\begin{figure}[H]
    \centering
    \begin{subfigure}[b]{0.485\textwidth}
        \centering
        \includegraphics[width=\textwidth]{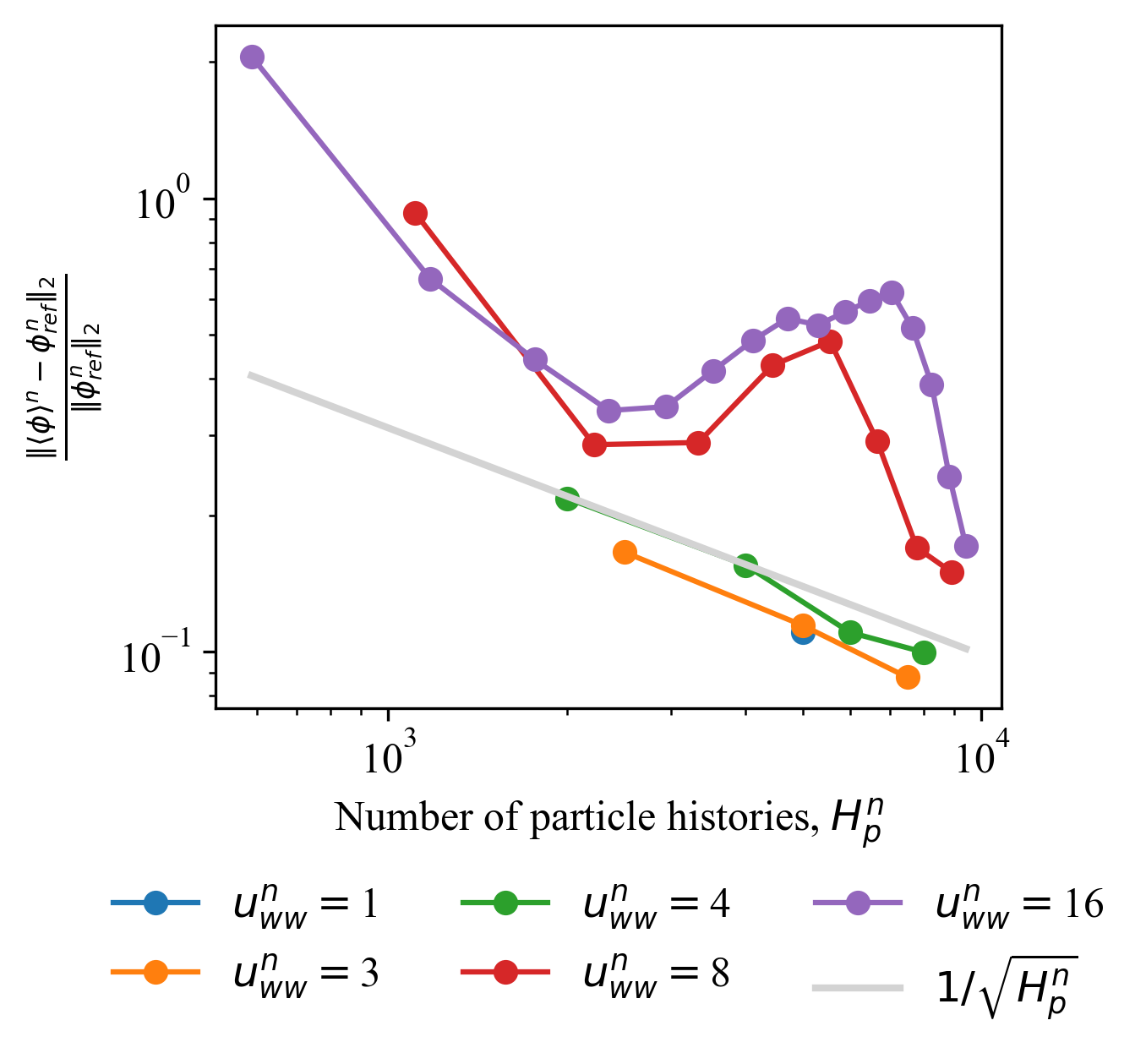} 
            \caption{$n=4$} 
    \end{subfigure}
    \hfill
    \begin{subfigure}[b]{0.485\textwidth}
        \centering
        \includegraphics[width=\textwidth]{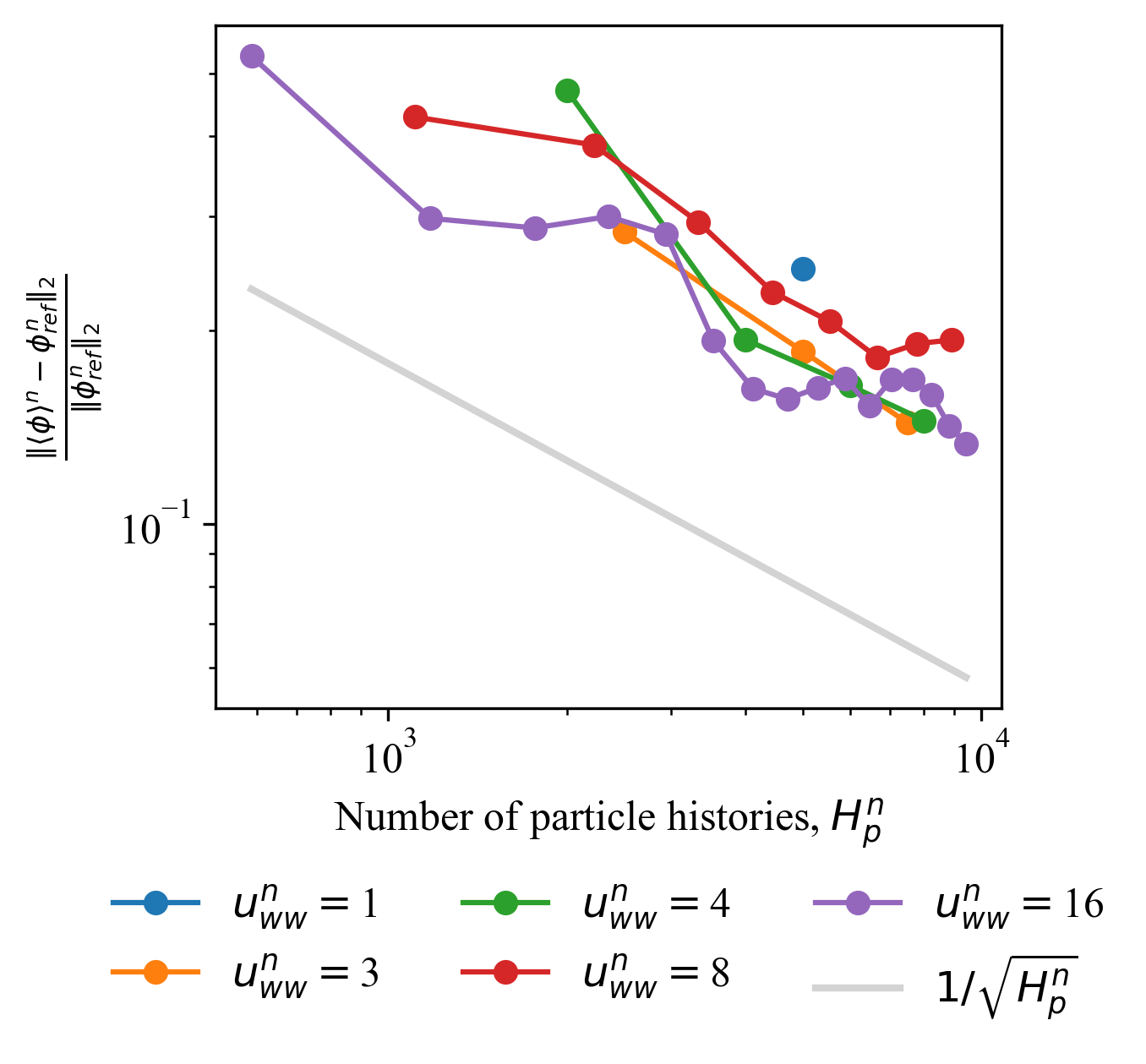} 
          \caption{$n=19$} 
    \end{subfigure}
    \caption{
 Relative $L_2$-norm of the error in the Monte Carlo scalar flux ($\langle\phi\rangle^n$) versus number of histories simulated before each weight-window updates   over  time steps $n=4$ and $n=19$ for different numbers of updates $u_{ww}^n$. }
    \label{fig:MC_relative_error}
\end{figure}

\begin{figure}[H]
    \centering
    \begin{subfigure}[b]{0.485\textwidth}
        \centering
        \includegraphics[width=\textwidth]{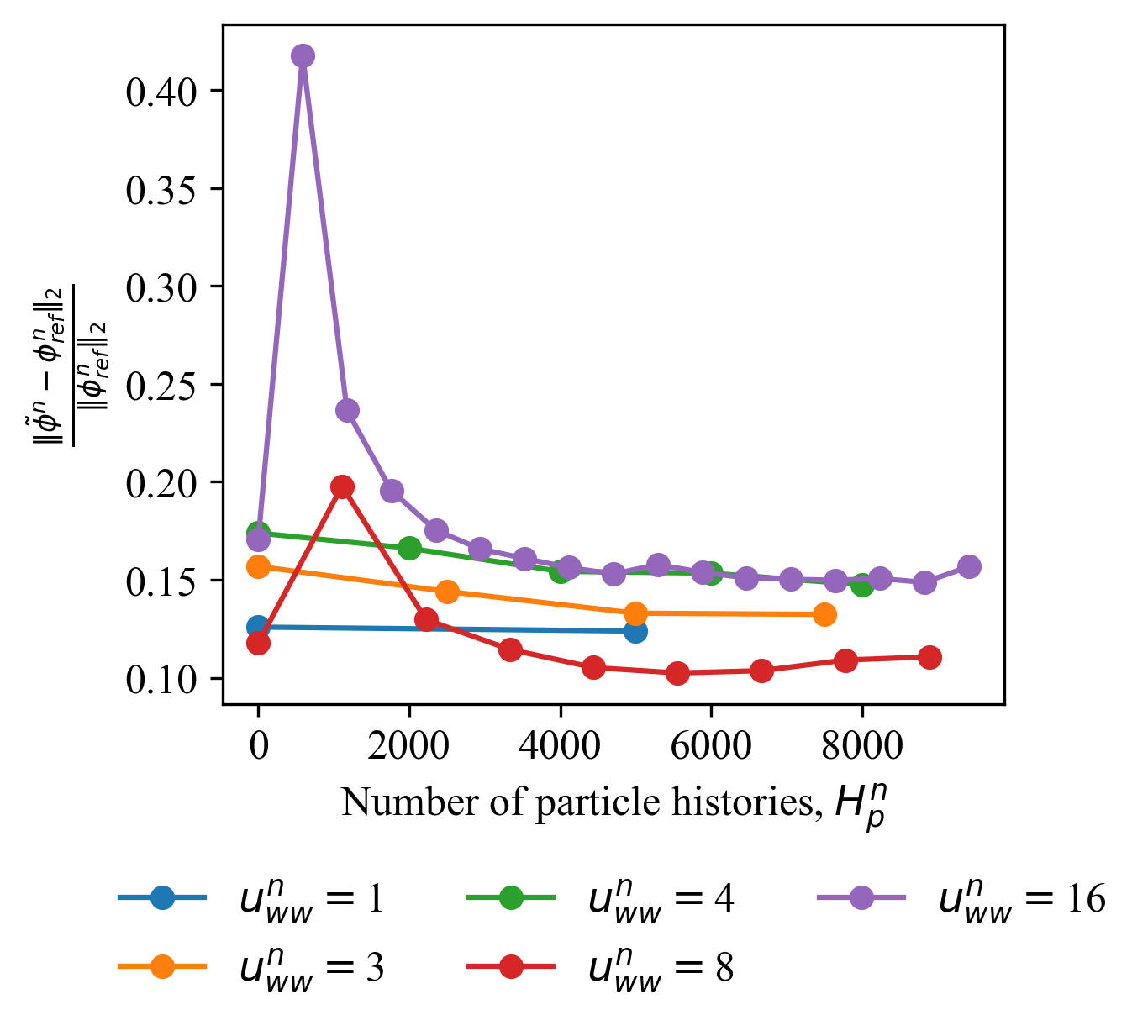}
      \caption{$n=4$}
    \end{subfigure}
    \hfill
    \begin{subfigure}[b]{0.485\textwidth}
        \centering
        \includegraphics[width=\textwidth]{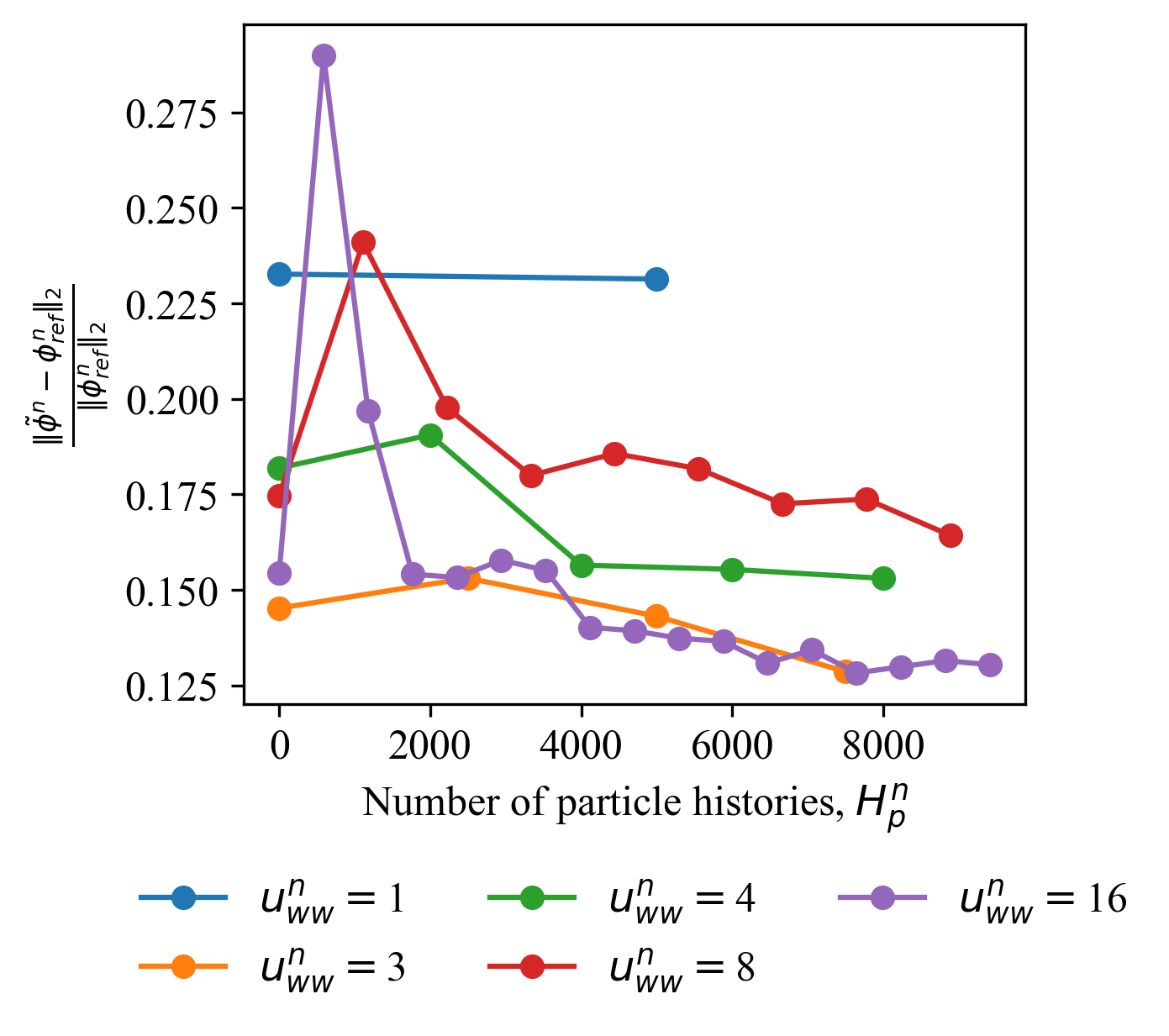} 
     \caption{$n=19$}
    \end{subfigure}
    \caption{
Relative $L_2$-norm of the error in the hybrid scalar flux ($\phi_{hlosm}$) versus number of histories simulated before each weight-window updates   over  time steps $n=4$ and $n=19$ for different numbers of updates $u_{ww}^n$.
 }    
 \label{fig:hybrid_relative_error}
\end{figure}

The hybrid weight window method has several parameters.
First, we consider the parameters $u_{ww}^n$   which define  the number    of updates of weight window centers using
 recalculated solution of the HLOSM equations.
Determining when and how many times to do these updates is based on two competing interests. Updating too soon means that the updated closures will 
have larger statistical error. 
On the other hand, if the weight windows are updated too late, the update could have little impact on the final tallies. 
To demonstrate effects of weight windows updates on the Monte Carlo solution,
we present the results for the test problem computed with 
different numbers of updates $u_{ww}^n$.
The test is solved  using $H^n = 10^4$ particle histories on each time step and
with  $f_p=\frac{p}{u_{ww}^n +1}$.
Note that the set $\{ f_p \}_{p=1}^{u_{ww}^n}$ defines frequency of updates in terms of number of  particle  histories.
The parameters of weight windows are  $\rho=1.25$ and 
$\epsilon_{min}=10^{-3}$.
The center of weight windows are computed by means of    numerical solution of the HLOSM equations discretized with the CN temporal scheme  and without filtering of closures and initial conditions. 
Figure \ref{fig:MC_relative_error}  presents the relative $L_2$-norm of the error in the scalar flux computed by  Monte Carlo algorithm, $\langle\phi\rangle^n$,  over two different time steps ($n=4,19$).
The relative error is  calculated using the benchmark solution, $\phi_{ref}$.
The results in  Figure \ref{fig:MC_relative_error} show that  the error decreases as $\frac{1}{\sqrt{H_p^n}}$ in the case of a small number of updates executed with relatively large fractions of  particle histories.
We also present the results of calculations with a very large number of updated $u_{ww}=8,16$. 
Such a number of updates is not efficient, but these cases were used
to analyze   feedback effects of frequent updates of weight windows on  Monte Carlo solution. The results show that
when the number of updates is large, the Monte Carlo solution converges, but there are stages over which the error increases. 
This indicates  potential instabilities of computations with large number and frequent updates of weight windows.
Figure \ref{fig:hybrid_relative_error} shows the relative error in the solution of HLOSM 
equations, $\phi_{hlosm}$,
which is used as an auxiliary solution $\tilde \phi$
to define weight windows (Eq. \eqref{ww-center})
in  Monte Carlo calculations presented in 
Figure \ref{fig:MC_relative_error}.
As more histories run, the error in closures decreases. This leads to improvement in  accuracy of numerical solution of the HLOSM equations.  
At some number of histories, the stochastic component of error in HLOSM solution becomes smaller then the error due to discretization in space and time.
As a result, the error limits to some value.
The point where the stochastic error of the hybrid solution becomes smaller than the effects of  discretization error depends on the mesh and number of histories \cite{vnn-dya-ans-annual-2024}.

\begin{figure}[h]
    \centering
    \begin{subfigure}[b]{0.45\textwidth}
        \centering
        \includegraphics[width=\textwidth]{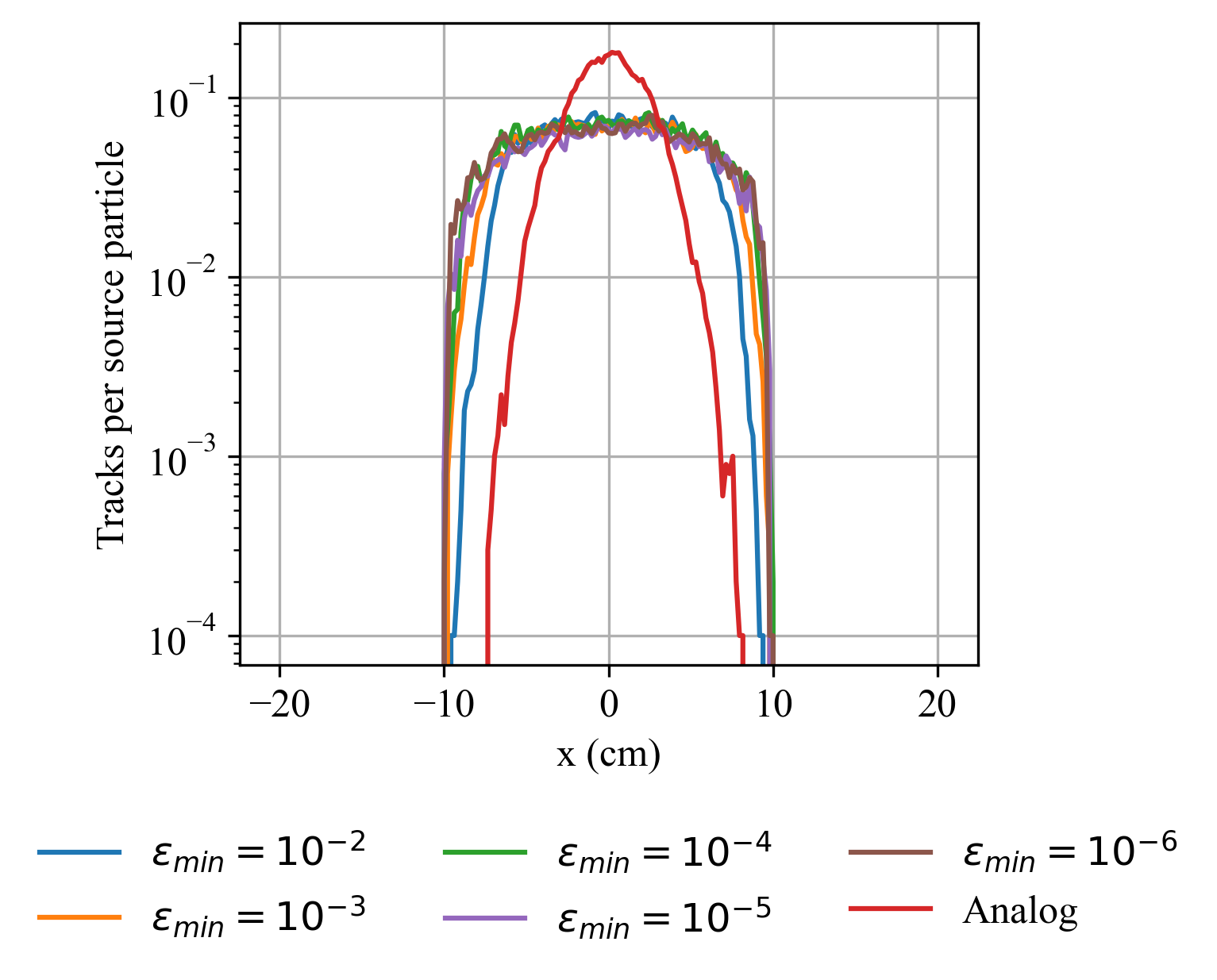} 
       \caption{$n=10$}
    \end{subfigure}
    \hfill
    \begin{subfigure}[b]{0.45\textwidth}
        \centering
        \includegraphics[width=\textwidth]{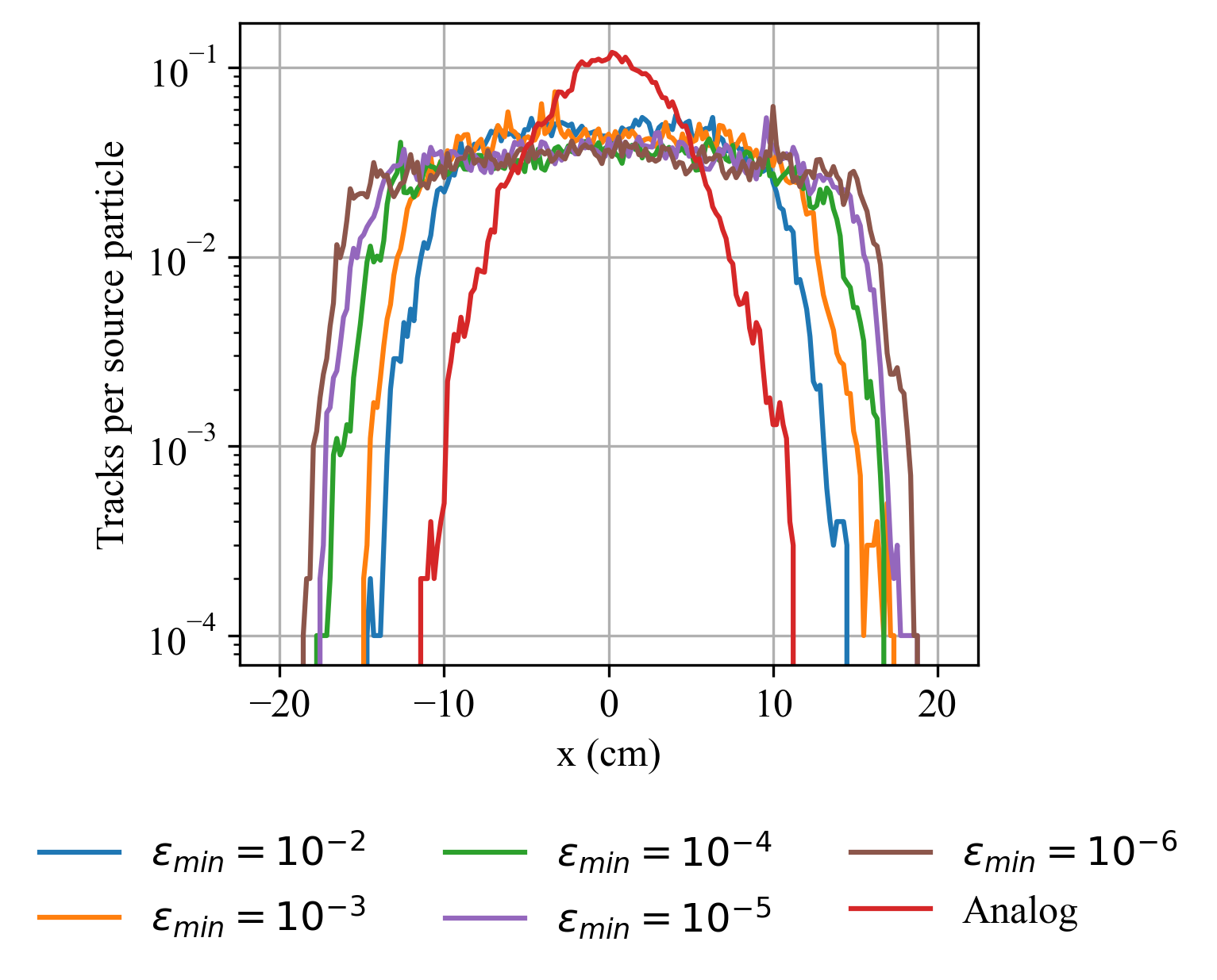} 
         \caption{$n=20$}
    \end{subfigure}
    \caption{
     Number of Monte Carlo particle tracks per source particle versus spatial position 
    at the time steps $n=10, 20$ in Monte Carlo simulations performed  with   weight windows using 
    $\epsilon_{min} \in \{10^{-2},10^{-3},10^{-4},10^{-5},10^{-6}\}$
    }
    \label{fig:epsilon_tracks1}
\end{figure}

The next parameter of interest is $\varepsilon_{min}$ 
which defines
the minimum value of the weight window center
to limit the amount of splitting.
Decreasing the parameter $\epsilon_{min}$ 
leads to better resolution of the wave front
at the expense of additional particles and therefore
increases computational cost.
Figure \ref{fig:epsilon_tracks1} shows the spatial distribution of particle tracks per source particle  at  time steps $n=10, 20$ solved with $H^n=10^4$ source particles and  window windows defined with  $\rho=1.25$ and various values of $\epsilon_{min}$. The  analog particle tracks
in Monte Carlo simulations without weight windows technique
  are peaked in the center, near the initial source of particles. 
This leads to  underpredicting of wavefront position. 
Applying weight windows  improves the resolution of the wave front with Monte Carlo particles.
At the  time step n=20, the best resolution is achieved with $\epsilon_{min}=10^{-6}$. In this case, the total particle population at the end of the time step grows by orders of magnitude.
  One source particle with initial weight  $w=1$ would need to become approximately $10^6$ particles each with weight $w\approx10^{-6}$.

\begin{figure}[h]
    \centering
    \begin{subfigure}[b]{0.45\textwidth}
        \centering
        \includegraphics[width=\textwidth]{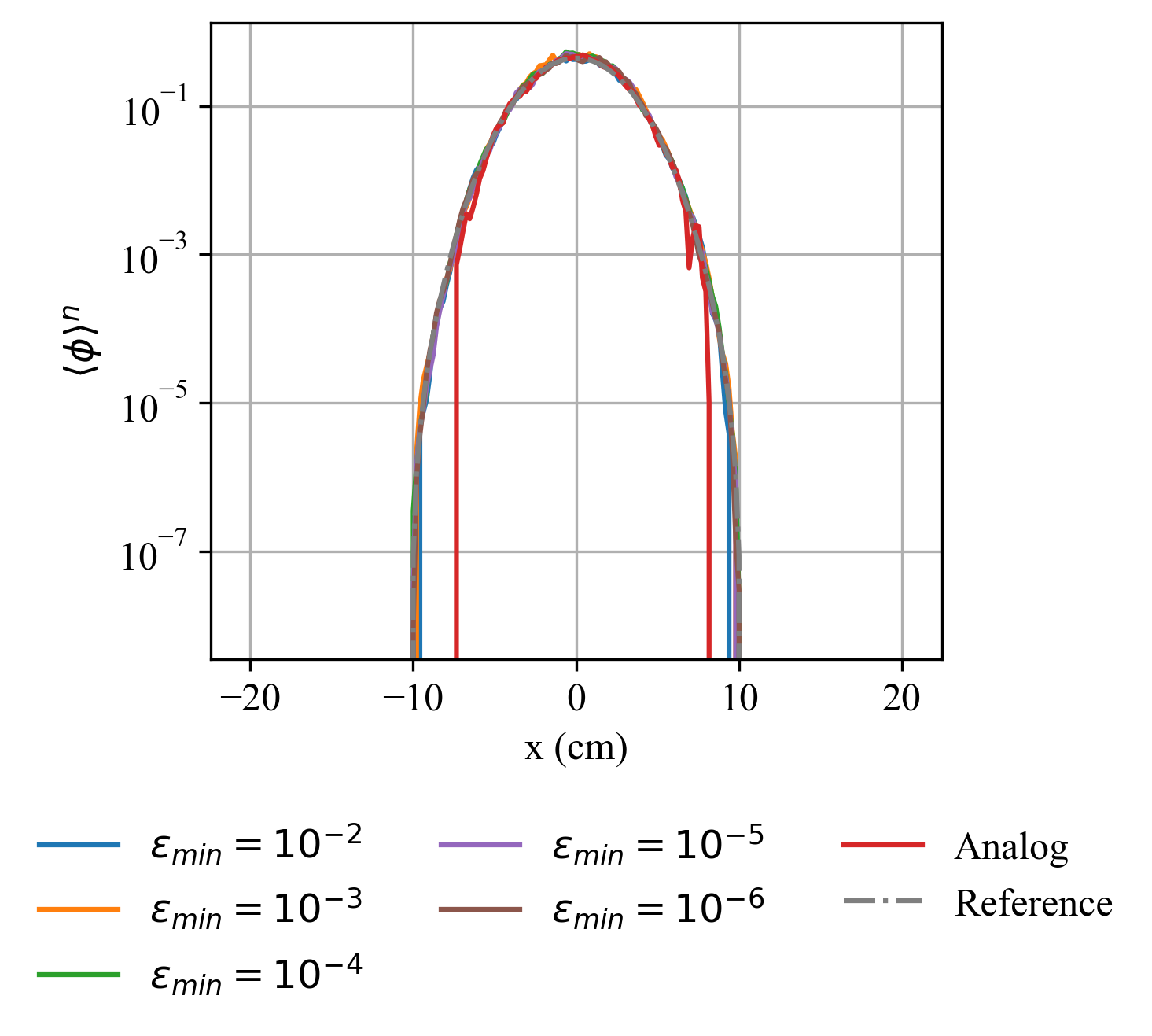} 
       \caption{$n=10$}
    \end{subfigure}
    \hfill
    \begin{subfigure}[b]{0.45\textwidth}
        \centering
        \includegraphics[width=\textwidth]{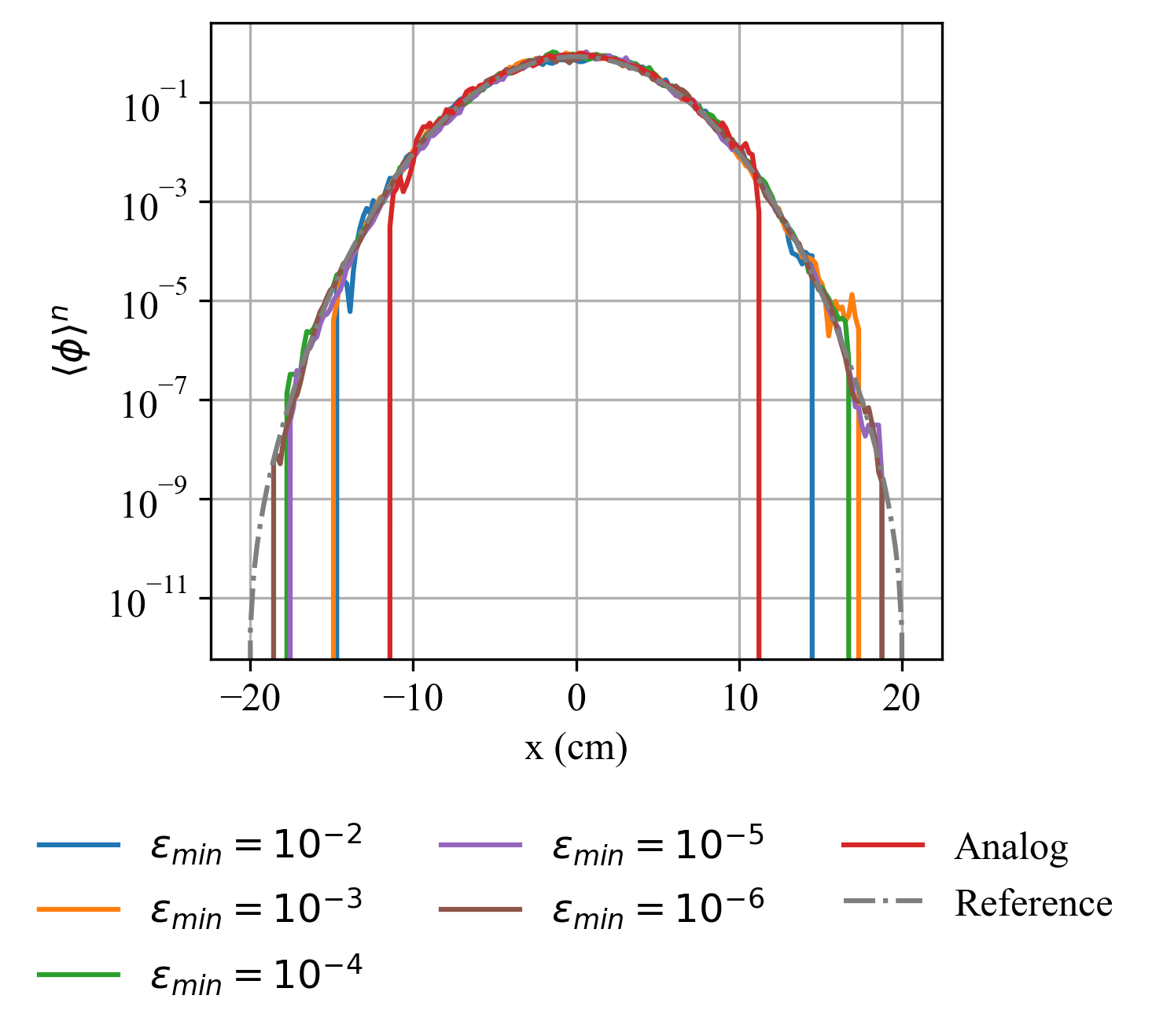} 
       \caption{$n=20$}
    \end{subfigure}
    \caption{
     Monte Carlo scalar flux, $\langle\phi\rangle^n$, versus spatial position 
    at the time steps $n=10, 20$ computed using   weight windows with 
    $\epsilon_{min} \in \{10^{-2},10^{-3},10^{-4},10^{-5},10^{-6}\}$}
    \label{fig:epsilon_solution}
\end{figure}

    The population control technique  applied in these calculations removes relatively high weight particles $w\approx1$ in the center near the source region. 
  This can lead to a greater total number of particle tracks with a lower maximum track density when using weight windows with population control. 
Figure \ref{fig:epsilon_solution} presents the Monte Carlo  scalar flux, $\langle\phi\rangle^n$ 
(Eq. \eqref{phi-mc}), demonstrating the effects of weight windows with different $\epsilon_{min}$ on position of wave front.

\FloatBarrier

\begin{figure}[H]
    \centering
    \begin{subfigure}[b]{0.45\textwidth}
        \centering
        \includegraphics[width=\textwidth]{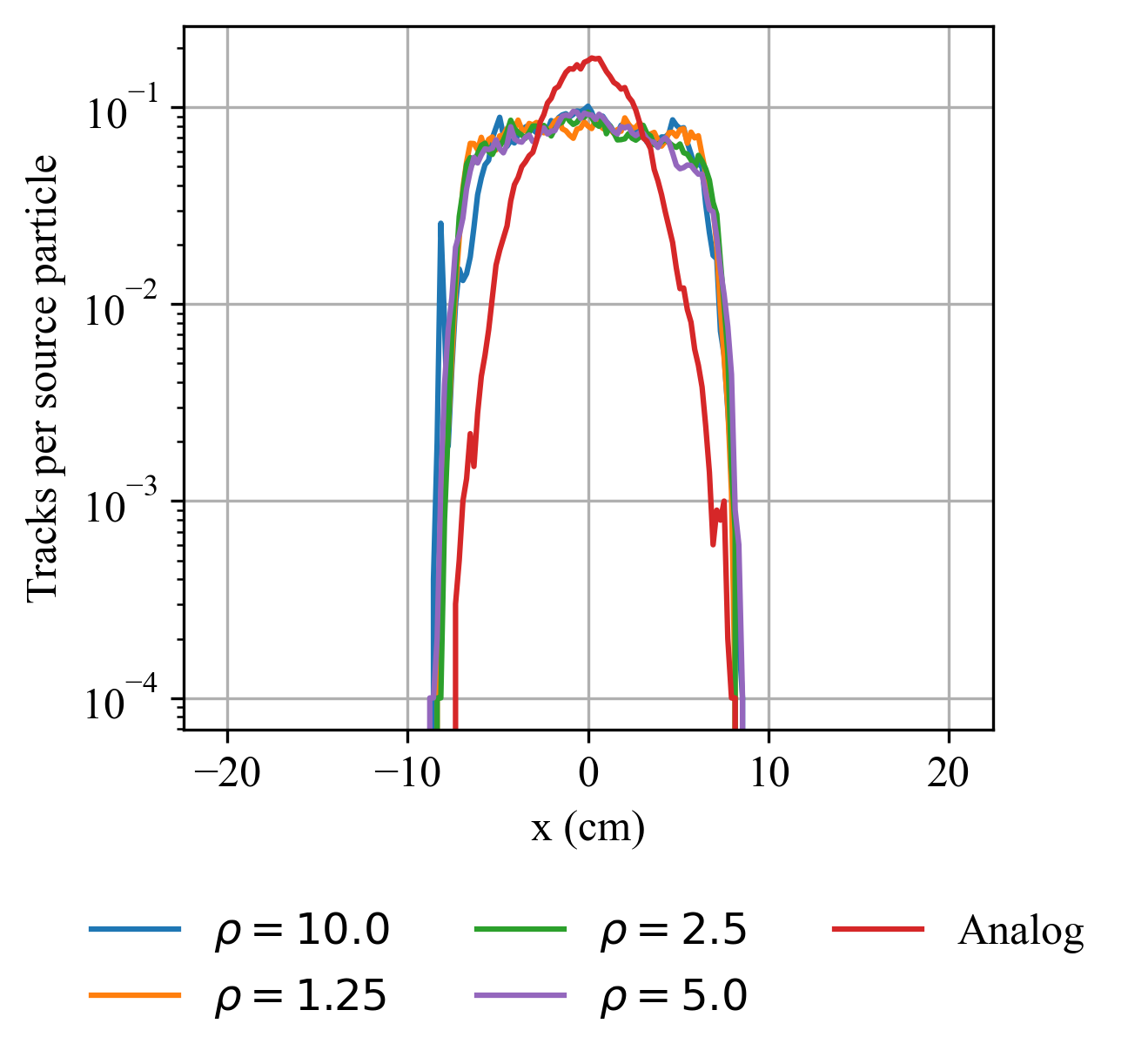} 
        \caption{$n=10$}
    \end{subfigure}
    \hfill
    \begin{subfigure}[b]{0.45\textwidth}
        \centering
        \includegraphics[width=\textwidth]{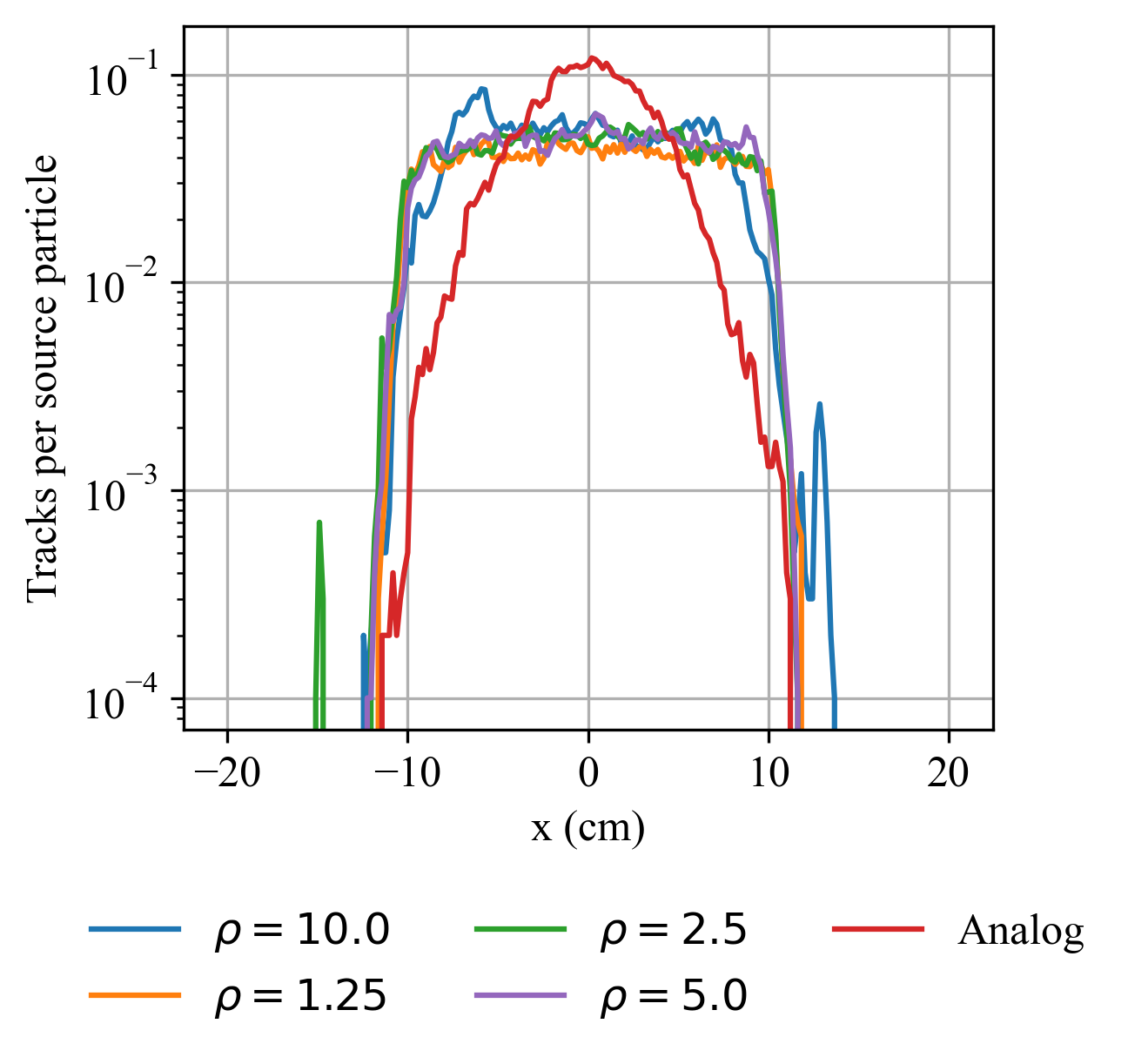} 
       \caption{$n=20$}
    \end{subfigure}
       \caption{
    Number of Monte Carlo particle tracks per source particle versus spatial position 
    at the time steps $n=10, 20$ in Monte Carlo simulations performed  with   weight windows using 
    $\rho \in \{1.25,2.5,5,10\}$}
    \label{fig:rho_tracks1}
\vspace{0.4cm}
    \centering
    \begin{subfigure}[b]{0.45\textwidth}
        \centering
        \includegraphics[width=\textwidth]{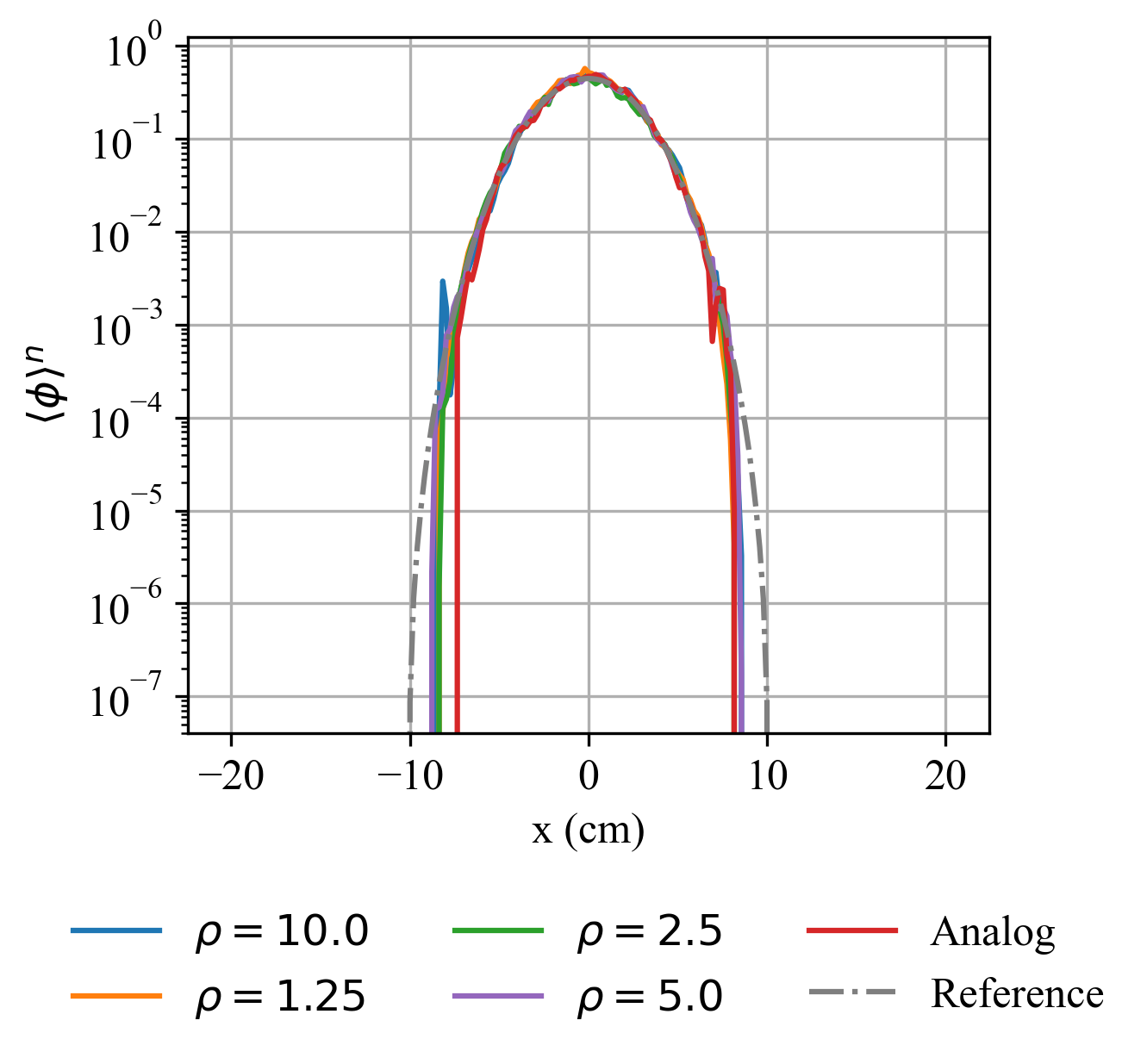} 
       \caption{$n=10$}
    \end{subfigure}
    \hfill
    \begin{subfigure}[b]{0.45\textwidth}
        \centering
        \includegraphics[width=\textwidth]{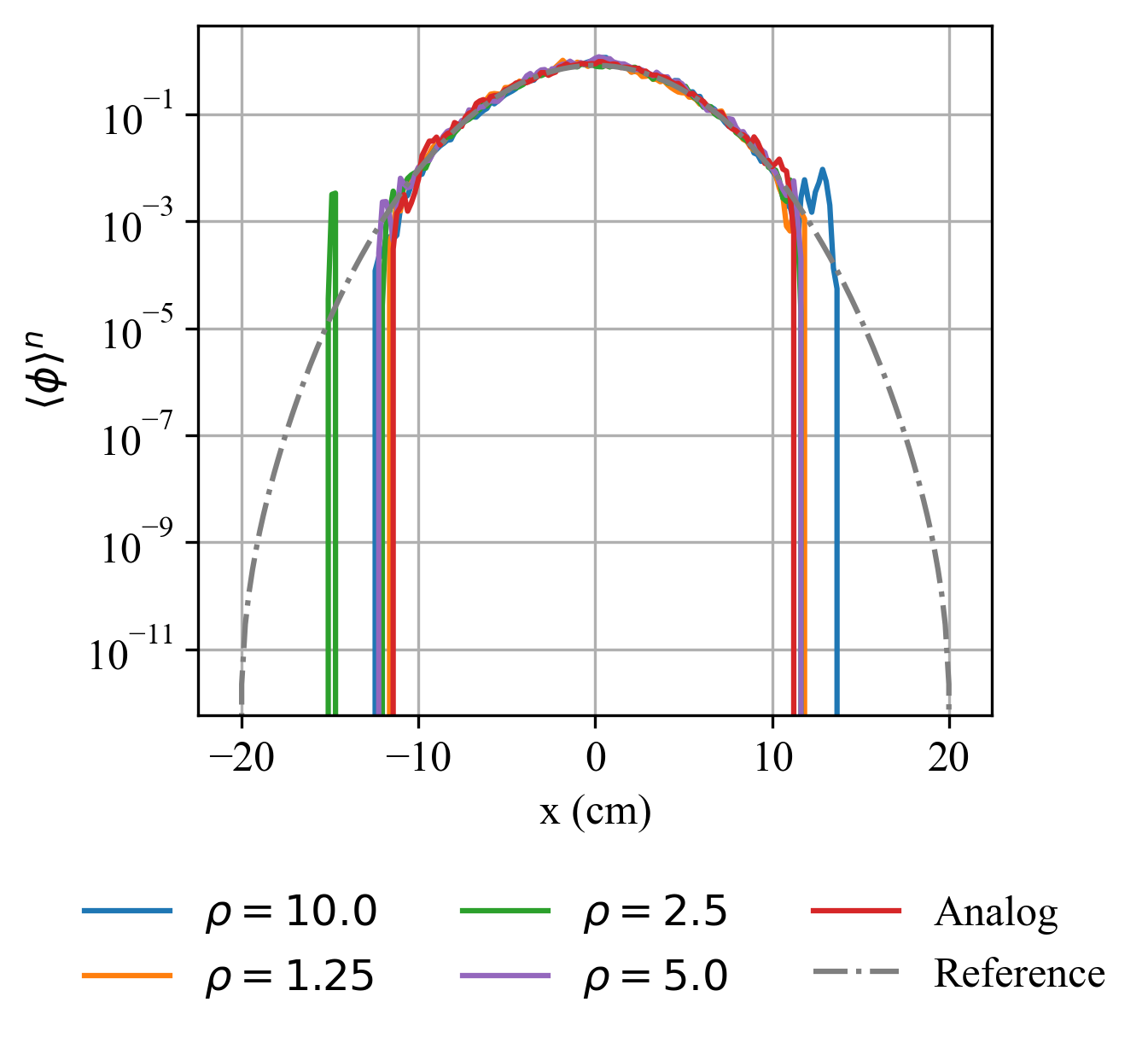}
         \caption{$n=20$}
    \end{subfigure}
    \caption{
    Monte Carlo scalar flux, $\langle\phi\rangle^n$, versus spatial position 
    at the time steps $n=10, 20$ computed using   weight windows with 
    $\rho \in \{1.25,2.5,5,10\}$}
    \label{fig:rho_solution}
\end{figure}

We now vary the parameter $\rho$ which 
defines the window width and controls over how tightly to constrain the particle weights to the desired distribution. 
Calculations are performed
with $\rho \in \{1.25,2.5,5,10 \}$.
The tests are solved with $H^n=10^4$, $\epsilon_{min}=10^{-3}$,  two updates of weight windows ($u_{ww}^n=2$) and  $f_p=\frac{p}{u_{ww}^n +1}$.
In Figure \ref{fig:rho_tracks1} the number of particle tracks per source particle is shown. 
The spatial distribution of particle tracks obtained with   $\rho=10$ is the least uniform. Calculations with $\rho=1.25,2.5,5$ yields particle distributions that are similar to each other.
Figure \ref{fig:rho_solution} shows the Monte Carlo solution using the considered set of $\rho$. Unlike, $\epsilon_{min}$, the ability for particles to reach the wave front is not significantly impacted by this parameter.

\FloatBarrier

\subsection{Numerical Solution of HLOSM Equations} 
\label{sec:numerical_solution}

\begin{figure}[H]
    \centering
    
    \begin{subfigure}[b]{0.45\textwidth}
        \centering
        \includegraphics[width=\textwidth]{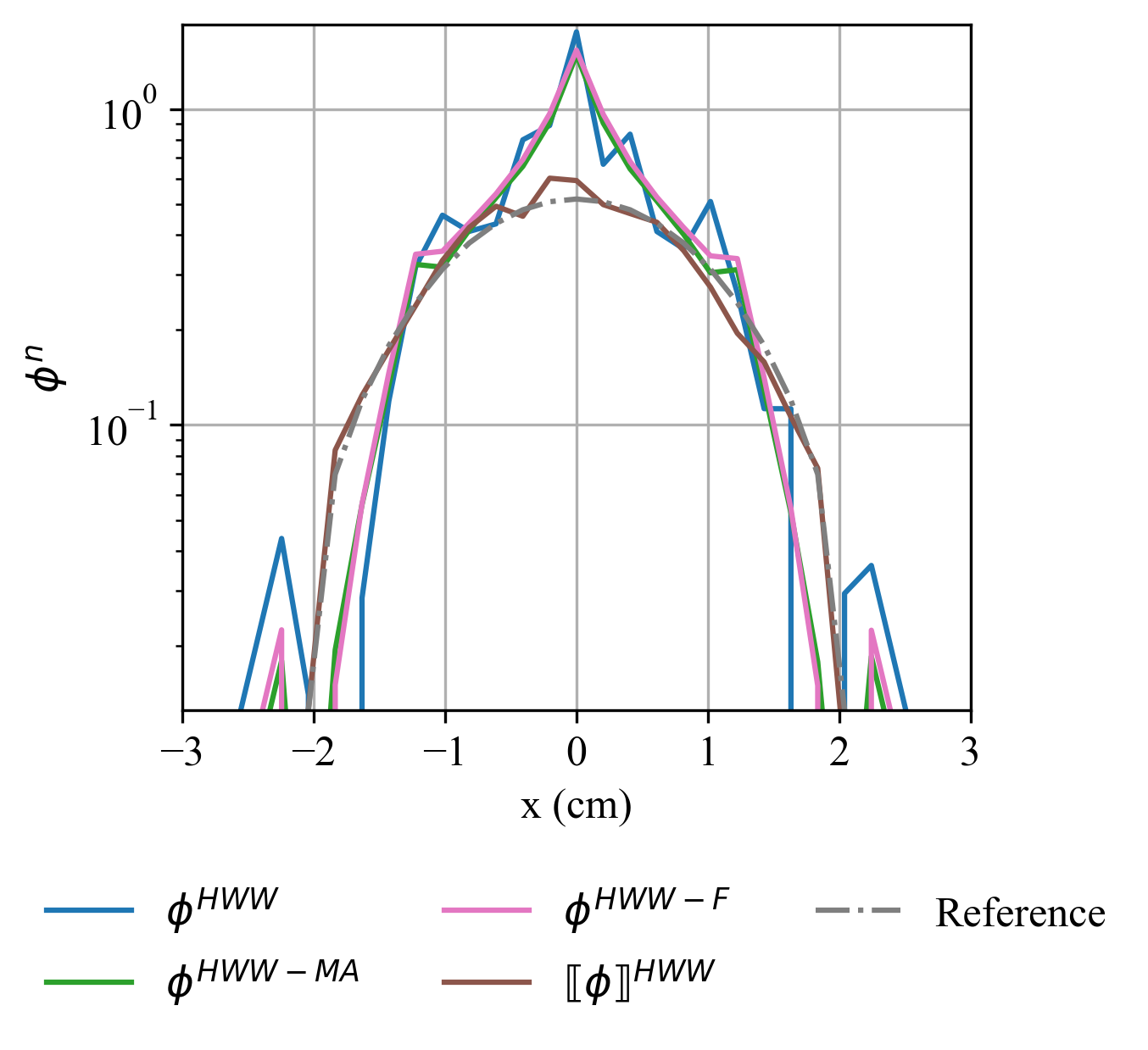}
        \vspace{-0.5cm}
          \caption{\label{fig:problem1_soln_a_filtering}
       $ {\phi}_{hlosm}^n$ and $\llbracket\phi\rrbracket^n$ at $n=2$ }
        \vspace{0.5cm}
    \end{subfigure}
    \hfill
    \begin{subfigure}[b]{0.45\textwidth}
        \centering
        \includegraphics[width=\textwidth]{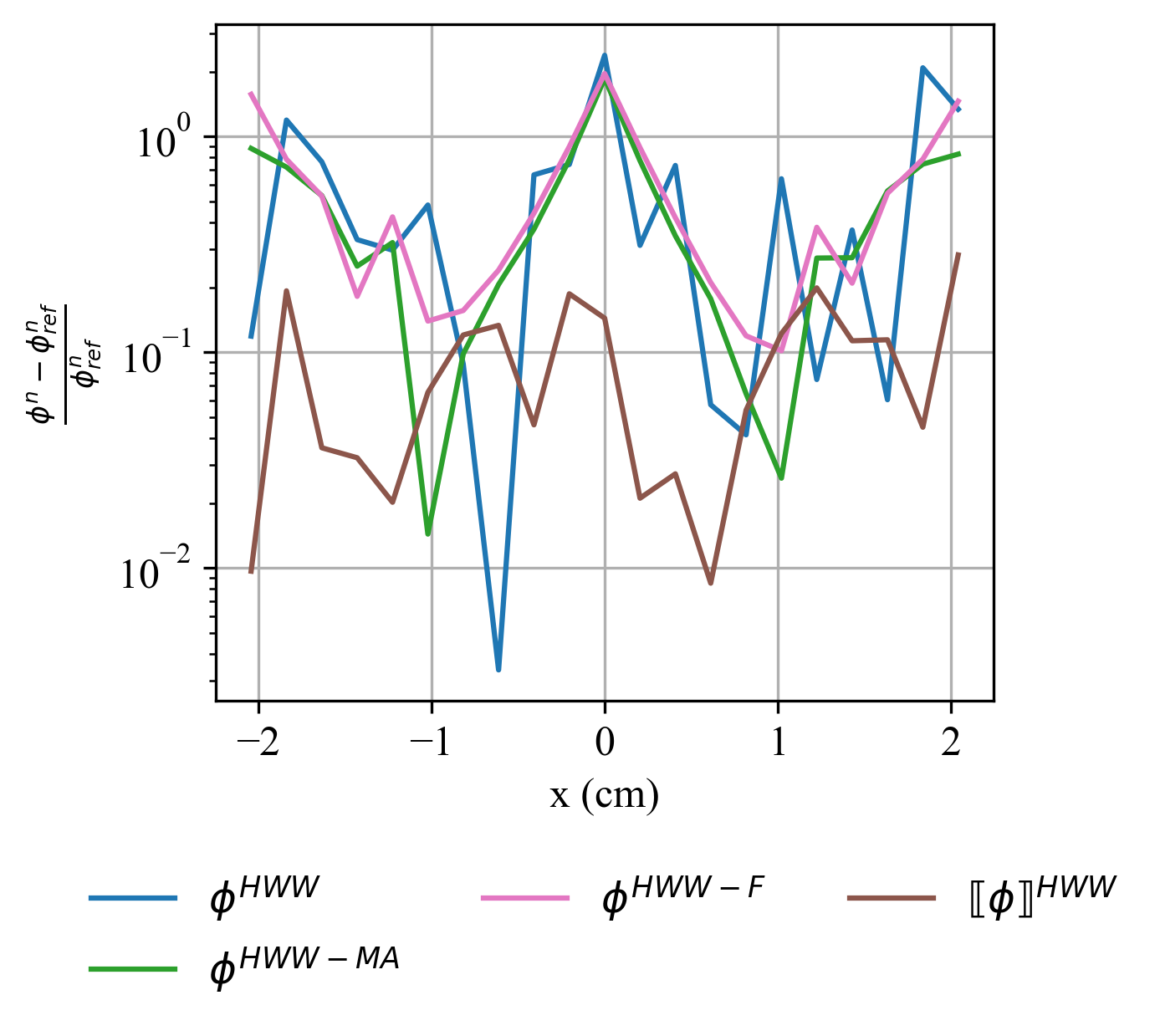}
       \vspace{-0.5cm}
        \caption{\label{fig:problem1_soln_b_filtering}
        Relative  errors at $n=2$}
    \vspace{0.5cm}
    \end{subfigure}
    \caption{
Numerical solutions $\phi_{hlosm}^n$, $\phi_{mc}^n$ and  their relative errors at $n=2$}
    \label{fig:Prob_1_filtering_n=2}
\end{figure}

\begin{figure}
    \begin{subfigure}[b]{0.45\textwidth}
        \centering
        \includegraphics[width=\textwidth]{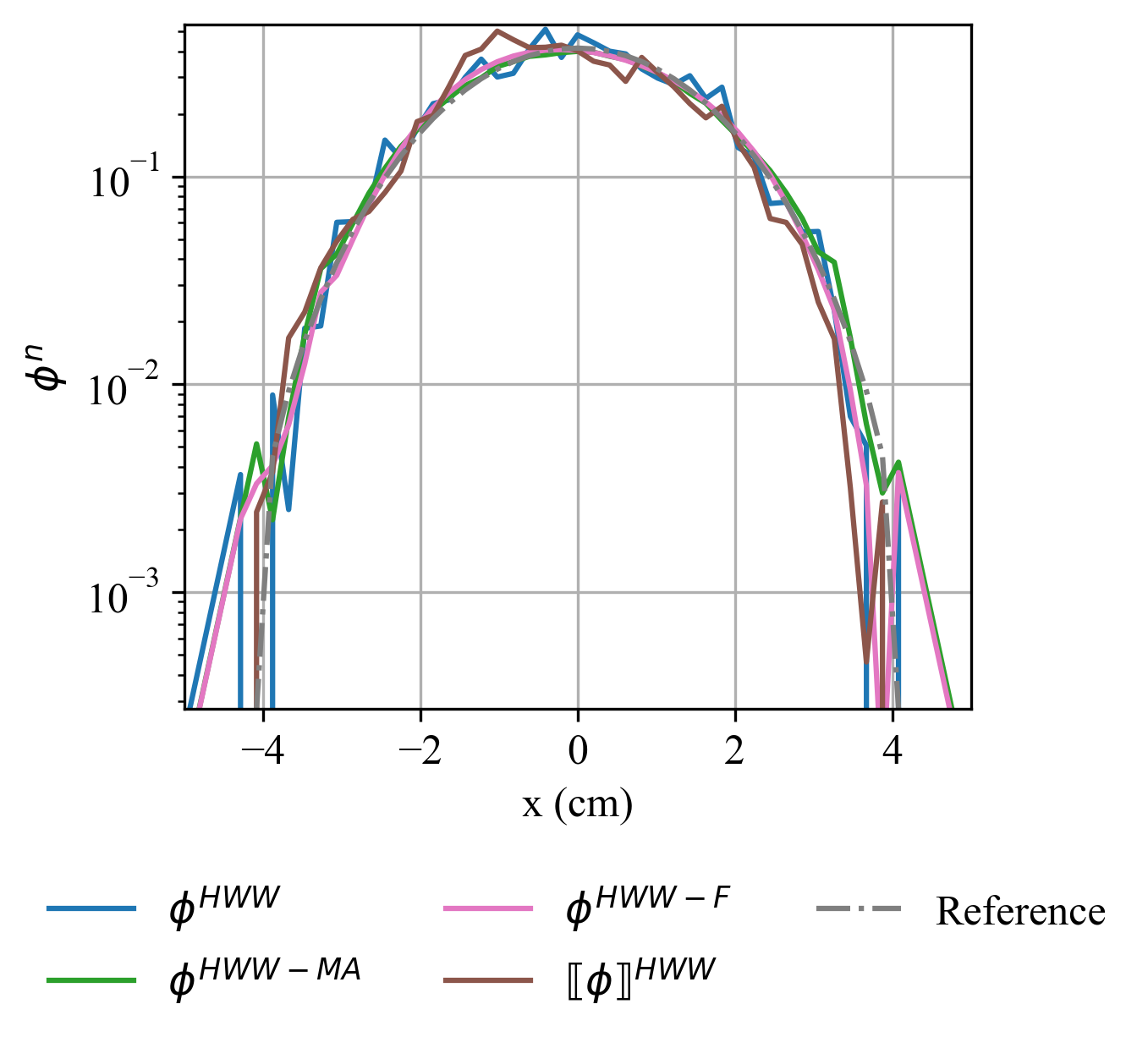}
         \vspace{-0.5cm}
        \caption{\label{fig:problem1_soln_c_filtering}
      $ {\phi}_{hlosm}^n$ and $\llbracket\phi\rrbracket^n$ at $n=4$ }
        \vspace{0.5cm}
    \end{subfigure}
    \hfill
    \begin{subfigure}[b]{0.45\textwidth}
        \centering
        \includegraphics[width=\textwidth]{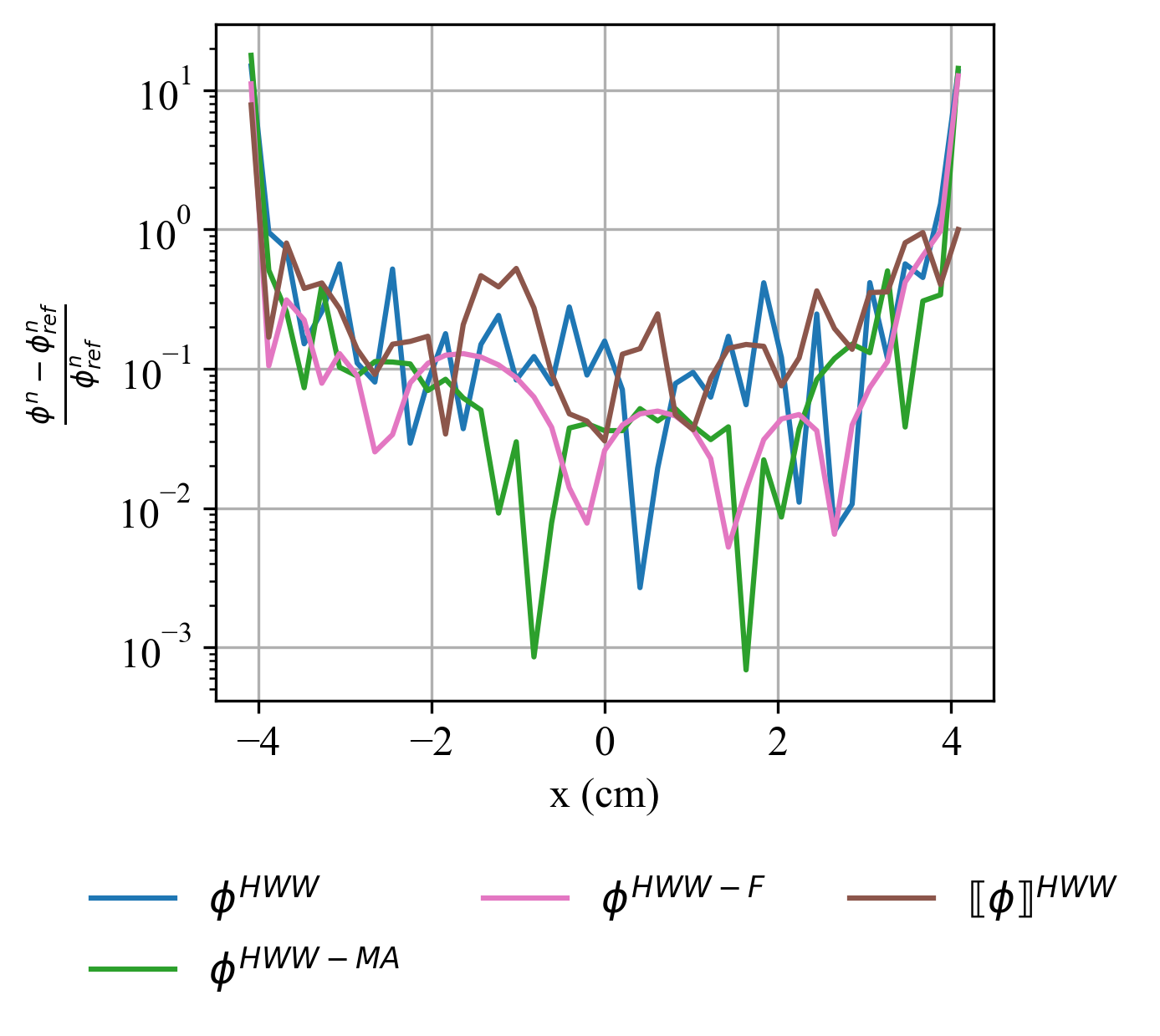}
         \vspace{-0.5cm}
        \caption{\label{fig:problem1_soln_d_filtering}
        Relative errors at  $n=4$ }
        \vspace{0.5cm}
    \end{subfigure}

        \caption{
Numerical solutions $\phi_{hlosm}^n$, $\phi_{mc}^n$ and  their relative errors at  $n=4$}
    \label{fig:Prob_1_filtering_n=4}
\end{figure}

\begin{figure}
    \begin{subfigure}[b]{0.45\textwidth}
        \centering
        \includegraphics[width=\textwidth]{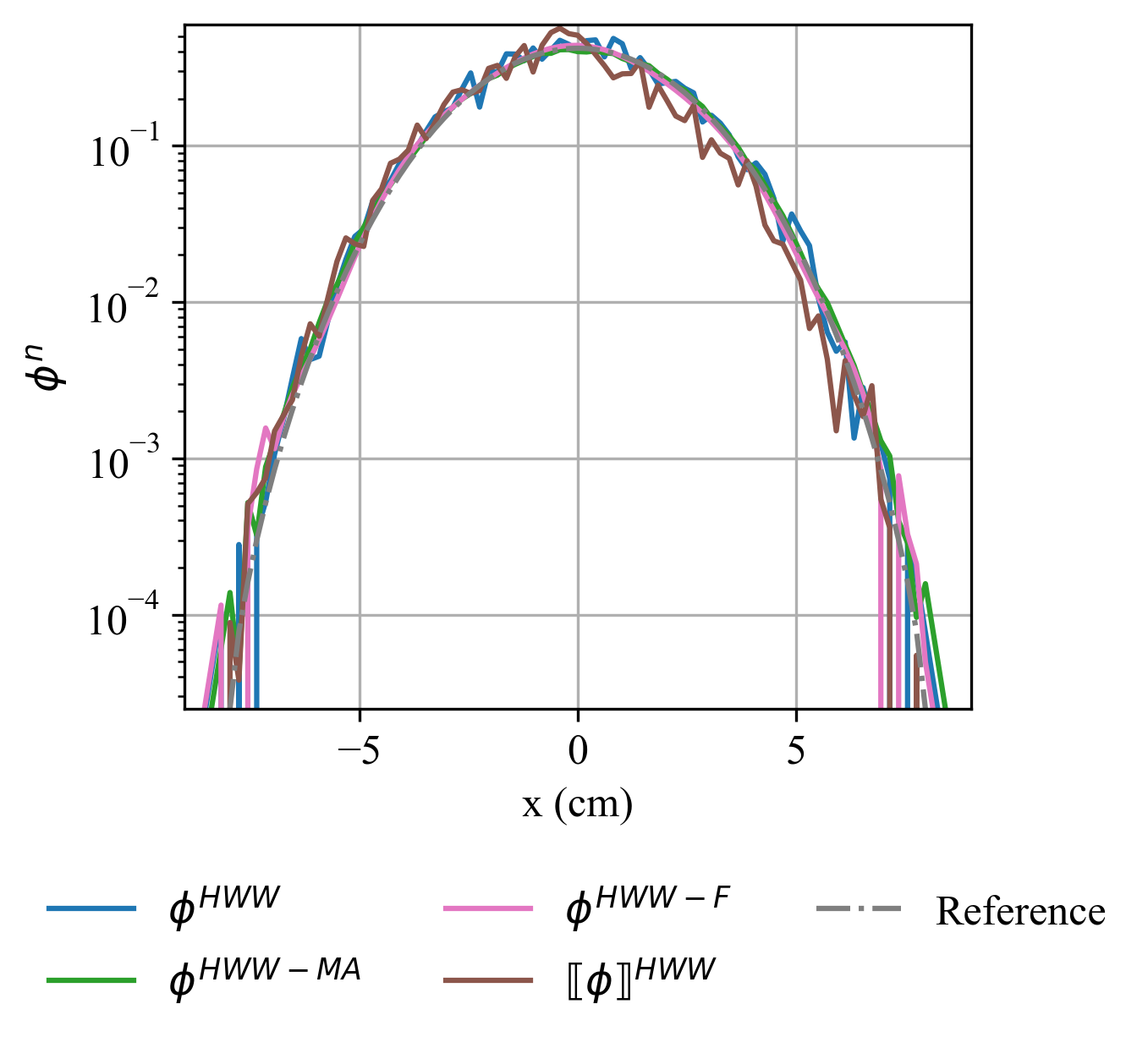}
         \vspace{-0.5cm}
        \caption{\label{fig:problem1_soln_e_filtering}
      $ {\phi}_{hlosm}^n$ and $\llbracket\phi\rrbracket^n$ at $n=8$ }
        \vspace{0.5cm}        
    \end{subfigure}
    \hfill
    \begin{subfigure}[b]{0.45\textwidth}
        \centering
        \includegraphics[width=\textwidth]{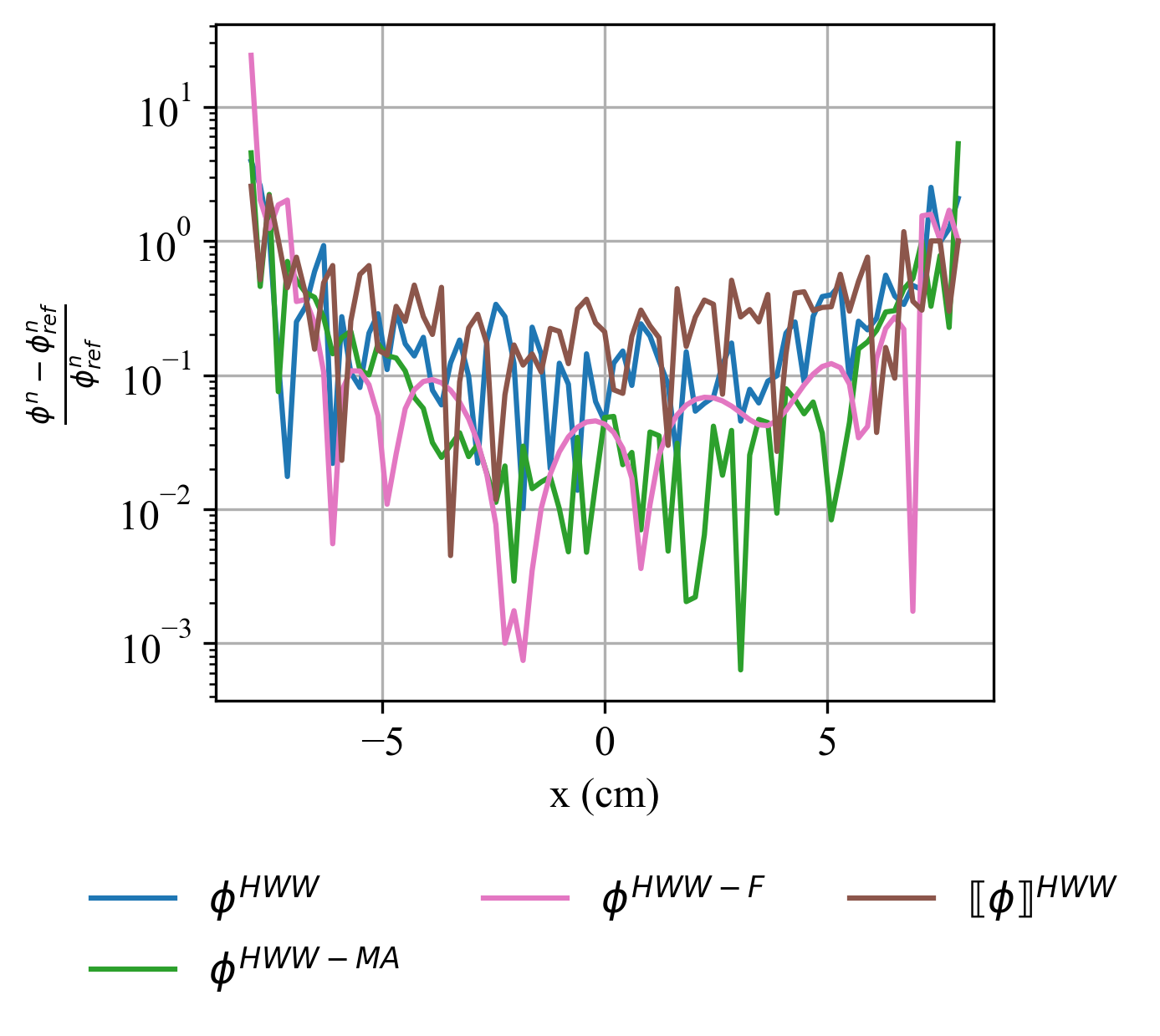}
         \vspace{-0.5cm}
        \caption{\label{fig:problem1_soln_f_filtering}
        Relative errors at  $n=8$ }
        \vspace{0.5cm}       
    \end{subfigure}
    
    \caption{
Numerical solutions $\phi_{hlosm}^n$, $\phi_{mc}^n$ and  their relative errors  at  $n=8$}
    \label{fig:Prob_1_filtering_n=8}
\end{figure}

\begin{figure}[H]
    \centering
    \includegraphics[width=0.45\linewidth]{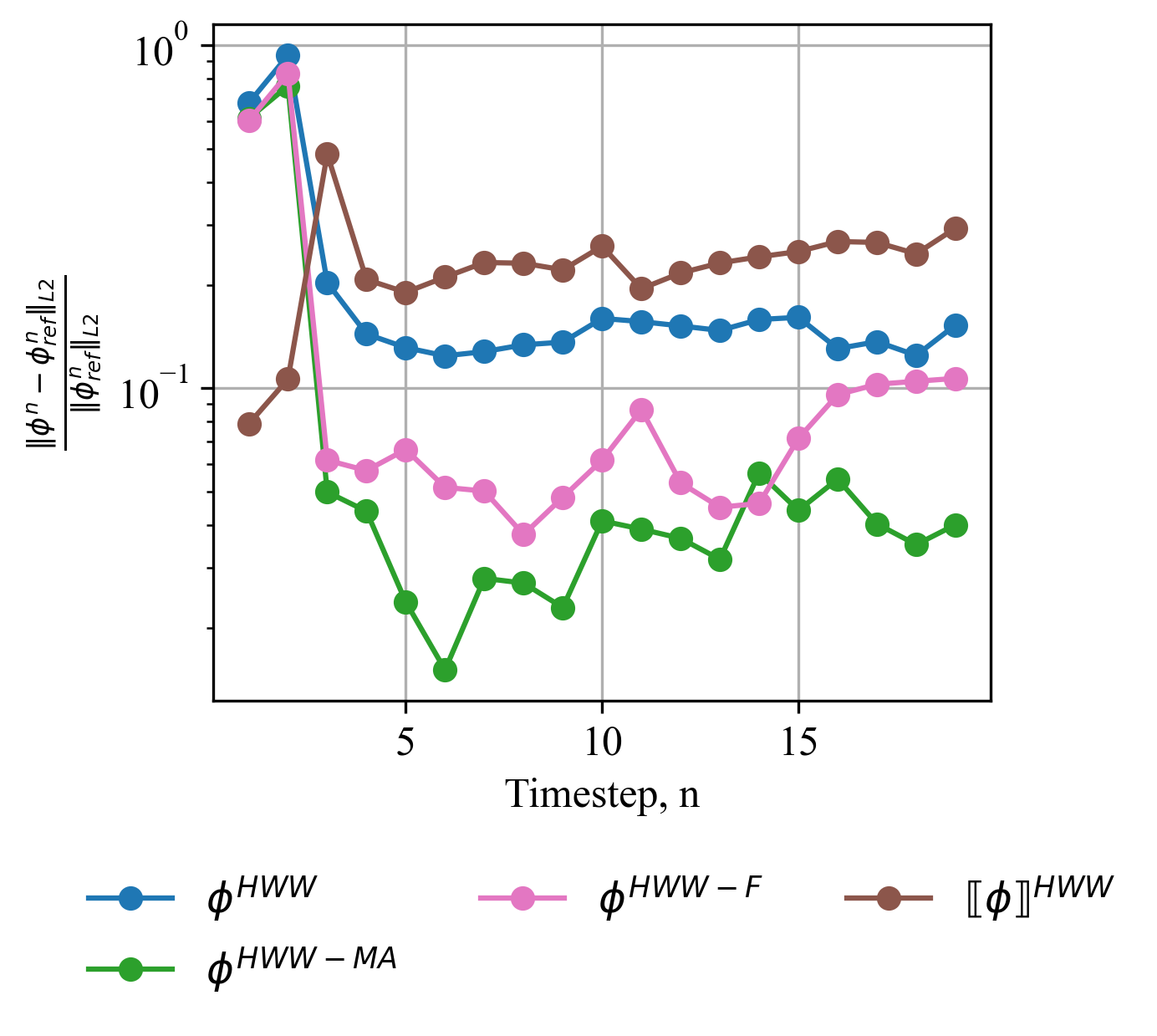}
    \caption{\label{fig:filtering_total}
     Relative $L_2$-norm of errors in  $\llbracket\phi\rrbracket^n$, 
 $ {\phi}_{hlosm}^n$ with and without filtering}
\end{figure}

In this section, we analyze the hybrid solution of the HLOSM equations used as an auxiliary numerical solution to define weight window centers.
The weight windows are
updated $u_{ww}^n=3$ times with update fractions $f_p=\frac{p}{u_{ww}^n +1}$   and hence
$ f_1=0.25,f_2=0.5,f_3=0.75$.
The values of $\rho$ and $\epsilon_{min}$ are $1.25$ and $10^{-3}$ respectively, and the number of source particles used is  $H^n=10^4$. 
Figures \ref{fig:Prob_1_filtering_n=2}, \ref{fig:Prob_1_filtering_n=4}, and \ref{fig:Prob_1_filtering_n=8}  
show
the hybrid solution computed using the
HWW algorithm, $ {\phi}^{HWW}$, the Monte Carlo
time-crossing tally solution, $\llbracket\phi\rrbracket^{HWW}$ (Eq. \eqref{ct-mc}), 
and  their relative errors     for the time steps $n=2,4,8$. These figures    also include  two hybrid solutions computed using filtered initial conditions and closures. One hybrid solution, $ {\phi}^{HWW-MA}$,
uses the MA filter with $k=3$
and the other, $ {\phi}^{HWW-F}$, the  Fourier filter with $\varkappa=30$. At time step $n=2$, the effects of discretizing the delta function source are still present, and they are not improved by filtering. At later times $n=4,8$, the unfiltered hybrid solution has similar or less error than the 
Monte Carlo solution. The filtered hybrid solutions at $n=4,8$ achieve lower relative error than either unfiltered hybrid or 
Monte Carlo solutions.
Figure \ref{fig:filtering_total} shows 
the relative $L_2$-norm of errors in $ {\phi}_{hlosm}^n$  and $\llbracket\phi\rrbracket^n$ as functions of time step.
By this measure, after $n=2$ the hybrid solution is 
more accurate than the Monte Carlo solution, and the filtered hybrid solutions are even better.  The computationally cheap MA
filter performs best overall.

\begin{figure}[h!]
	\centering
	\begin{subfigure}[b]{0.45\textwidth}
		\centering
    \includegraphics[width=\linewidth]{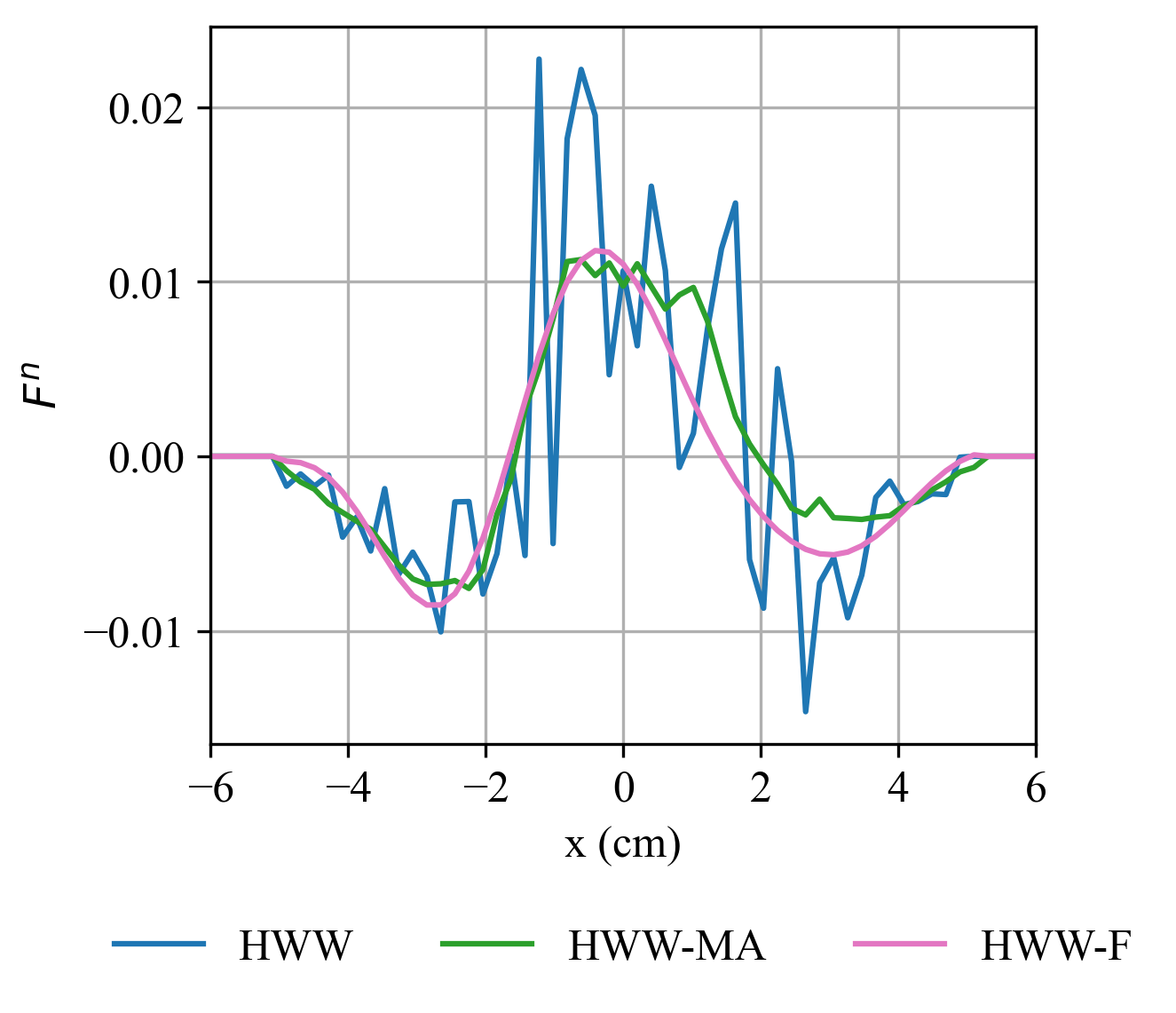}
    \vspace{-0.5cm}
		\caption{$n=5$ \label{fig:closure_filtering5}}
	\end{subfigure}
	\hfill
	\begin{subfigure}[b]{0.45\textwidth}
		\centering
      \includegraphics[width=\linewidth]{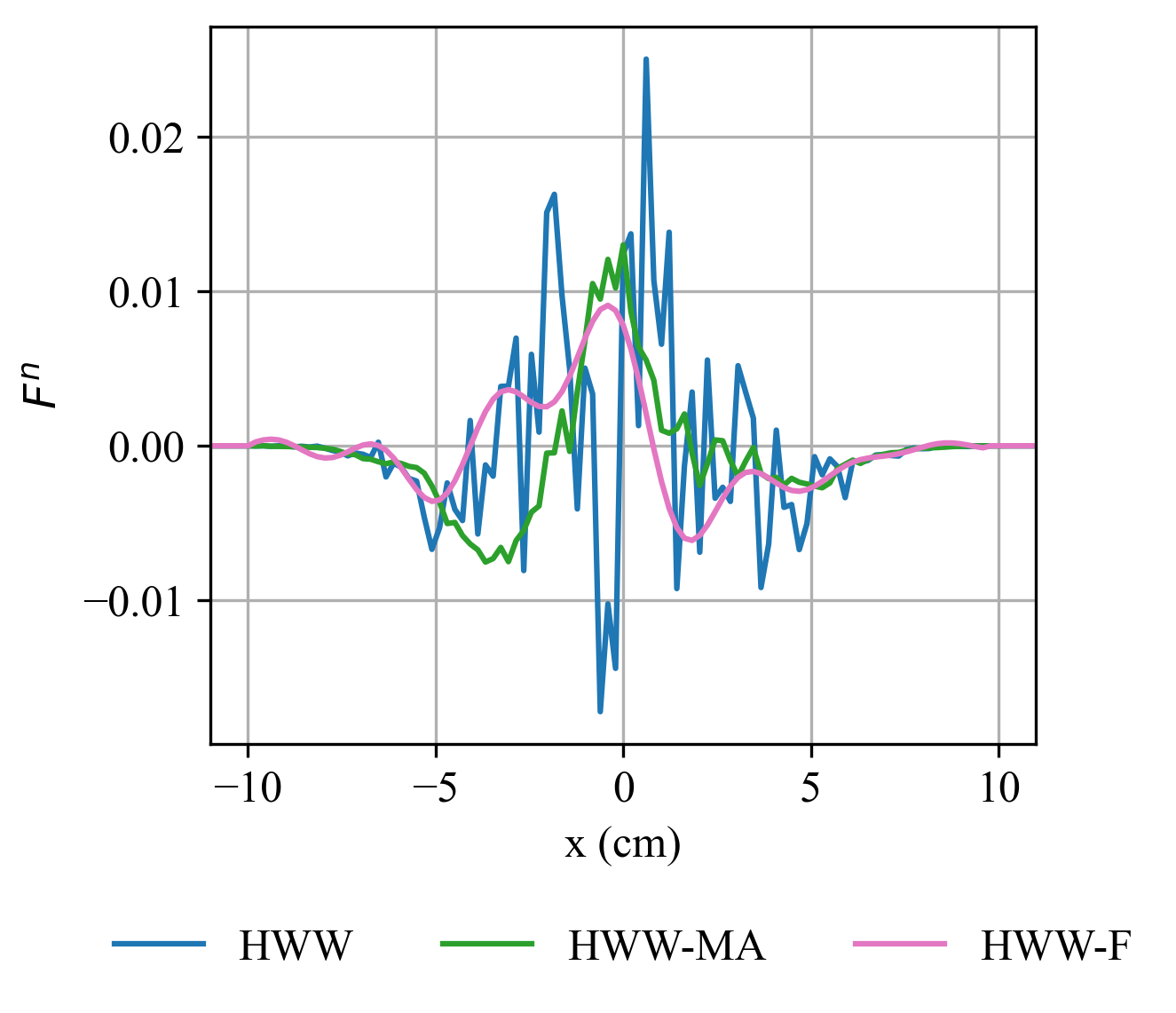}
    \vspace{-0.5cm}
		\caption{$n=10$ \label{fig:closure_filtering10}}
	\end{subfigure}
    \vspace{-0.25cm}
    \caption{Closure function $F^n$ and effects of filtering  \label{fig:closure_filtering}}
\end{figure}

Figure \ref{fig:closure_filtering}  shows
the function $F^n$ which defines the closure in the HLOSM equations. This figure  demonstrates the cell-average values of  $F^n$ at time steps $n=5, 10$  computed by Monte Carlo as well as this closure function after applying  filter  techniques.  HWW denotes the unfiltered closure, HWW-MA the MA filter, and HWW-F the Fourier filter. The parameters of the MA and Fourier filters are $k=3$ and  $\varkappa=30$, respectively.
At the time step $n=5$, both filters lead to  functions with a similar large-scale shape. 
The Fourier filter performs better in removing statistical noise at this stage. 
At  later times as the particles spread out, the stochastic noise increases.
At $n=10$ the filters with the given parameters 
 generate functions with different structures while
 reducing significantly  the noise in the closure.  


\subsection{Numerical Solution of Monte Carlo Algorithm with Hybrid Weight Windows}

We now present numerical   solutions of Monte Carlo algorithms  with (a) the hybrid weight windows (HWW) using $ {\phi}_{hlosm}^n$,  (b)
the lagged weight windows  (LWW) defined with the Monte Carlo solution $\langle\phi\rangle^{n-1}$ from the previous time step,  and (c)  the analog Monte Carlo without weight windows.  We also include the versions of the HWW algorithms with MA and Fourier filtering, denoted HWW-MA and HWW-F respectively. The parameter of the MA filter is $k=3$. The Fourier filter is applied with   $\varkappa=30$.
The test problem is solved with 
$H_p^n=10^4$ 
source particles and the same weight window parameters as in Section \ref{sec:numerical_solution}; $\rho=1.25$, $\epsilon_{min}=10^{-3}$, and $u_{ww}^n=3$ with $f_p=\frac{p}{u_{ww}^n +1}$.

\begin{figure}[t]
    \centering
    \begin{subfigure}[b]{0.45\textwidth}
        \centering
        \includegraphics[width=\textwidth]{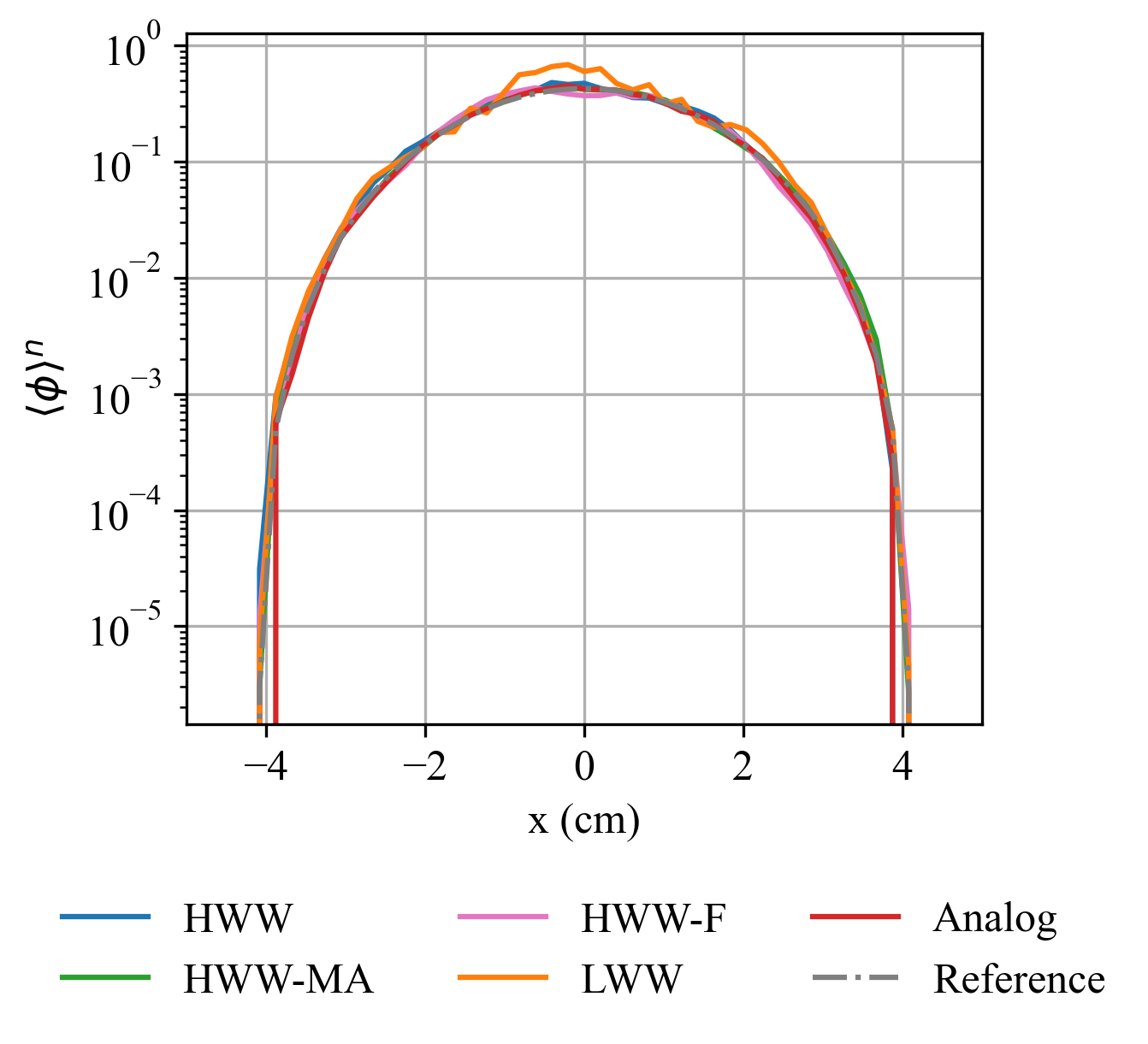}
        \caption{$n=4$}
        \label{fig:problem1_soln_4}
    \end{subfigure}
    \hfill
    \begin{subfigure}[b]{0.45\textwidth}
        \centering
        \includegraphics[width=\textwidth]{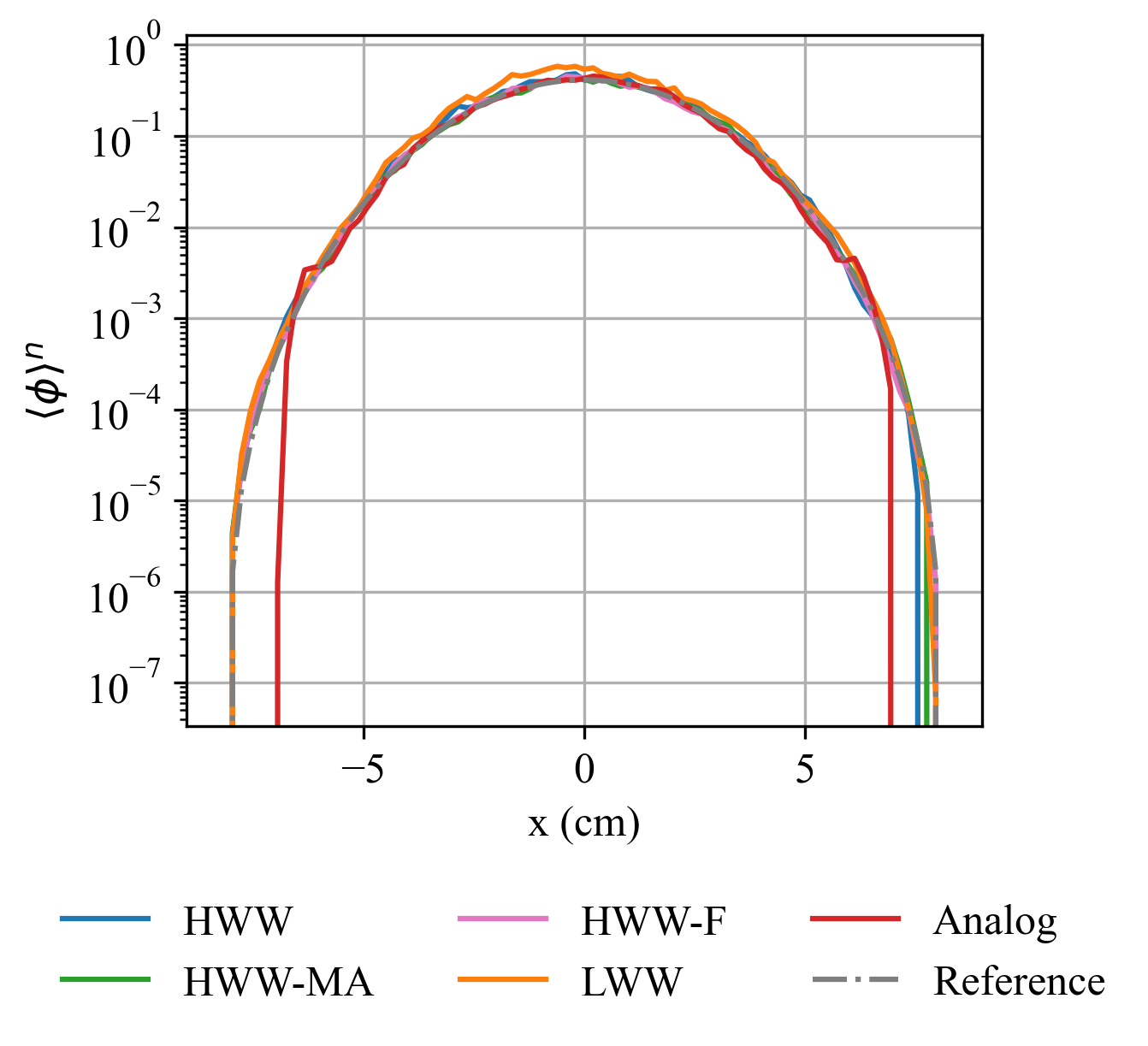}
        \caption{$n=8$}
        \label{fig:problem1_soln_8}
    \end{subfigure}
    \hfill
    \begin{subfigure}[b]{0.45\textwidth}
        \centering
        \includegraphics[width=\textwidth]{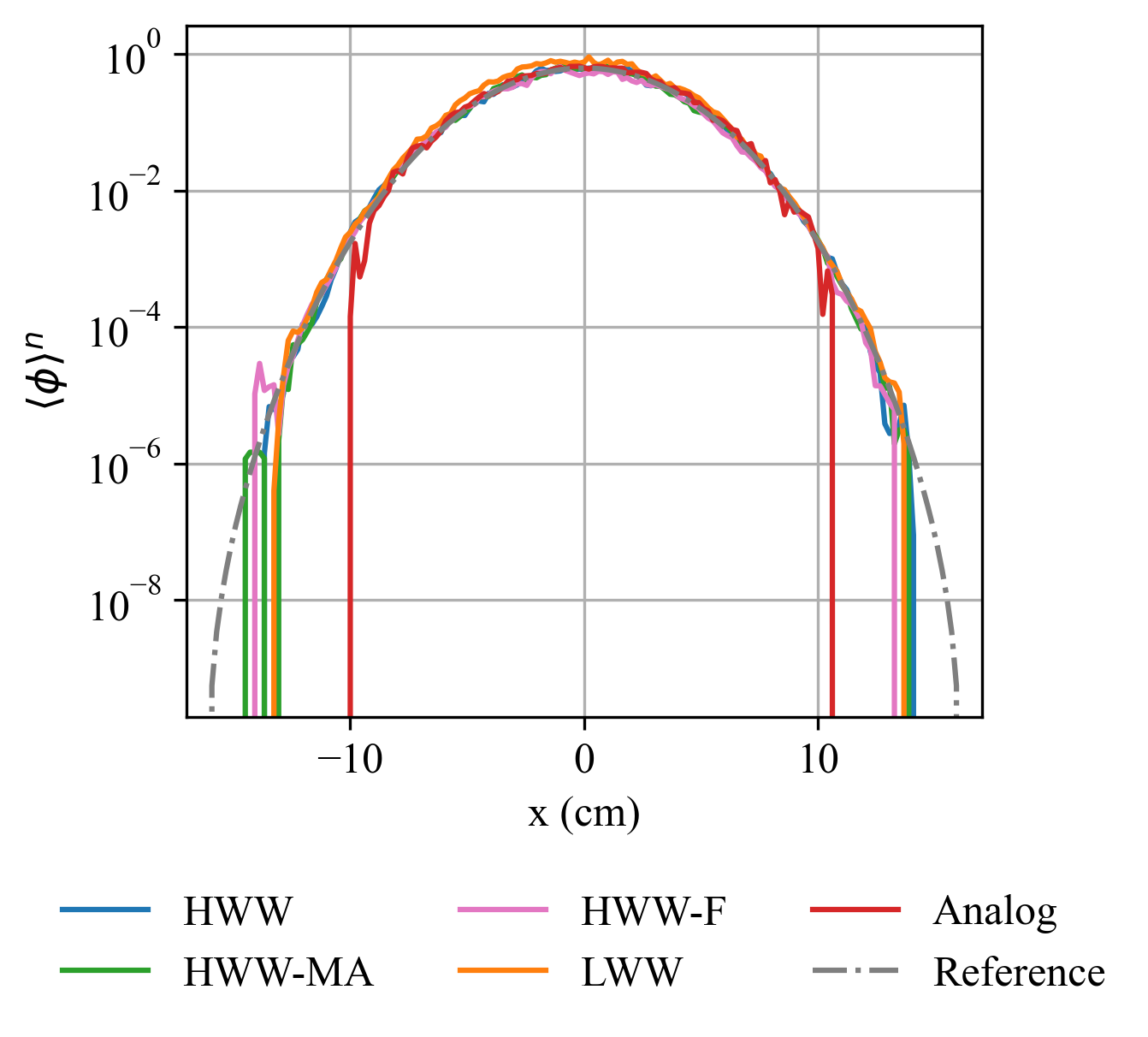}
        \caption{$n=16$}
        \label{fig:problem1_soln_16}
    \end{subfigure}
    \caption{\ $\langle\phi\rangle^n$, solution of Monte Carlo algorithms  at  time steps $n=4,8,16$}
    \label{fig:Prob_1_phi}
\end{figure}
\begin{figure}[t!]
    \centering
    \begin{subfigure}[b]{0.45\textwidth}
        \centering
        \includegraphics[width=\textwidth]{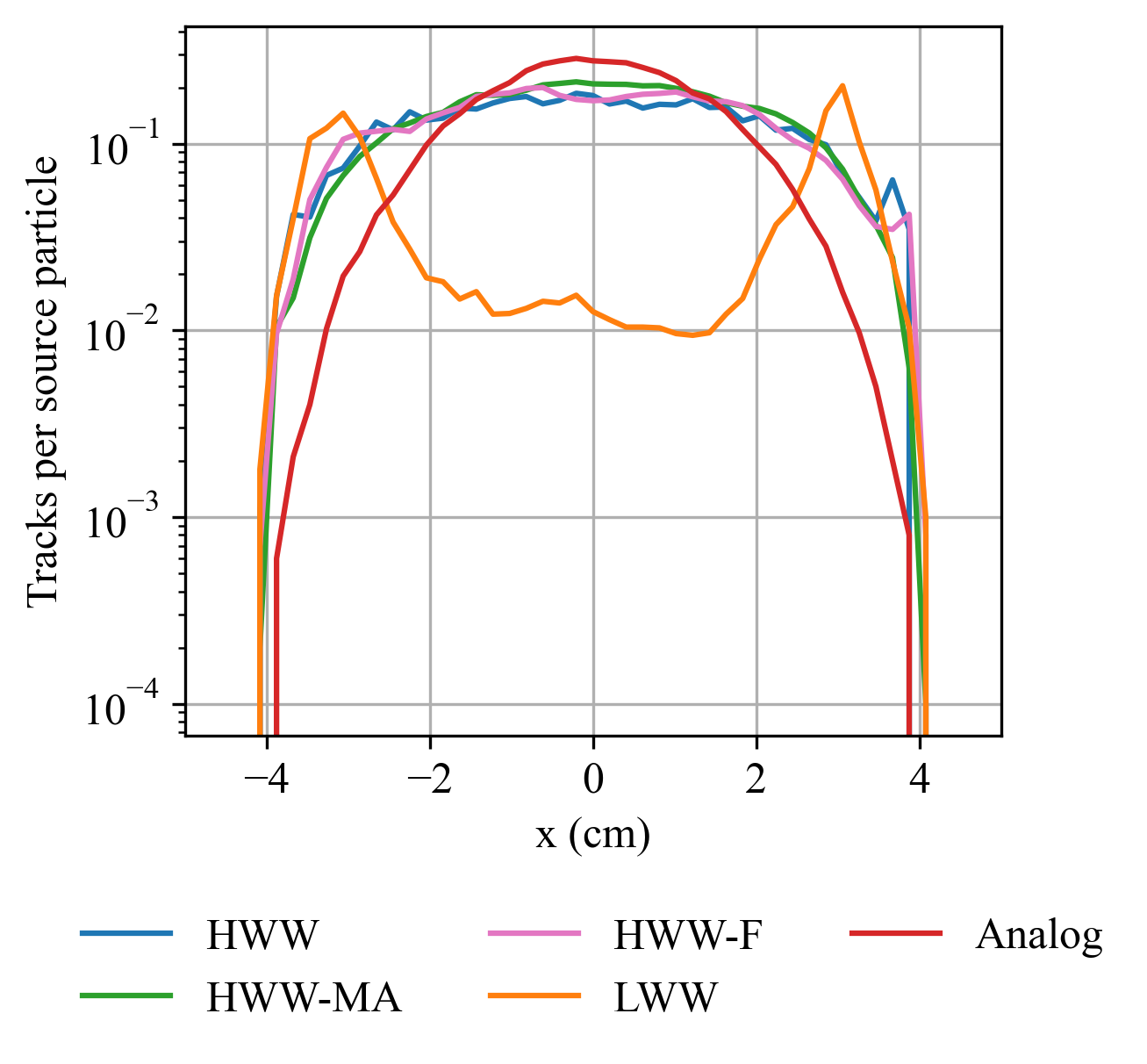}
        \caption{$n=4$}
        \label{fig:problem1_track_4}
    \end{subfigure}
    \hfill
    \begin{subfigure}[b]{0.45\textwidth}
        \centering
        \includegraphics[width=\textwidth]{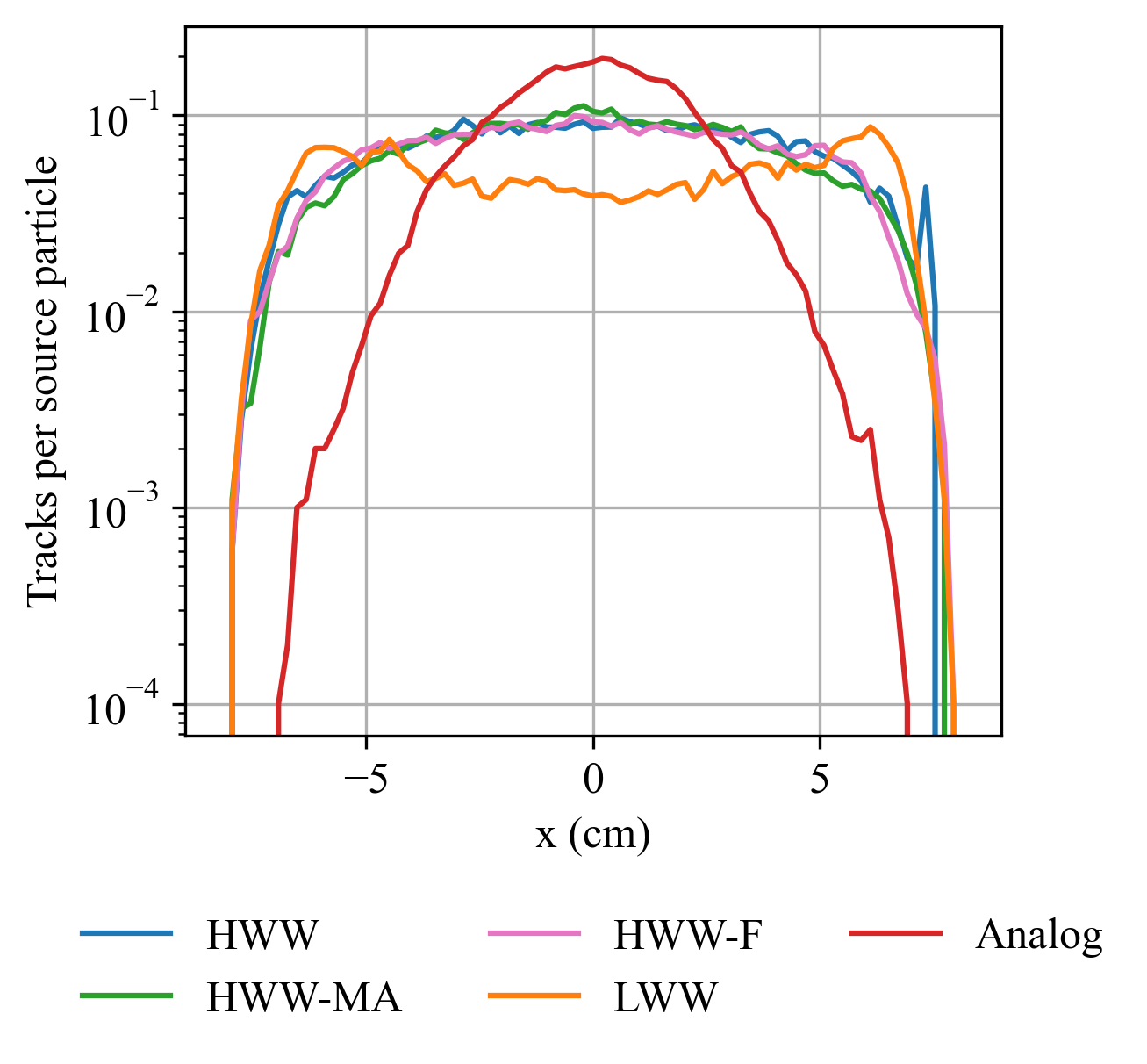}
        \caption{$n=8$}
        \label{fig:problem1_track_8}
    \end{subfigure}
      \hfill  
    \begin{subfigure}[b]{0.45\textwidth}
        \centering
        \includegraphics[width=\textwidth]{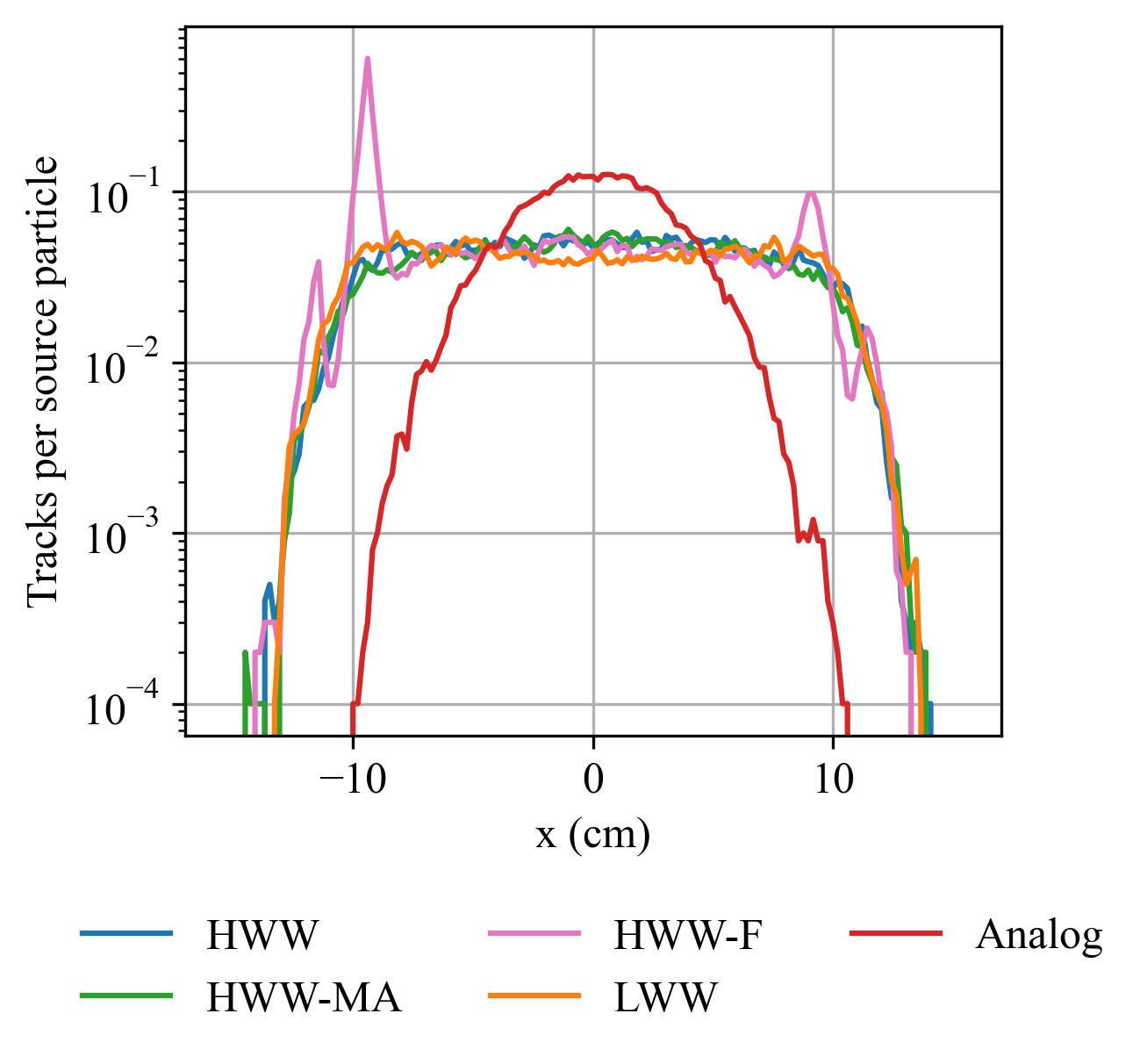}
        \caption{$n=16$}
        \label{fig:problem1_track_16}
    \end{subfigure}
    \caption{Number of particle tracks per source particle at  time steps $n=4,8,16$
        }
    \label{fig:Prob_1_tracks}
\end{figure}
\begin{figure}[H]
    \centering
    \begin{subfigure}[b]{0.45\textwidth}
        \centering
        \includegraphics[width=\textwidth]{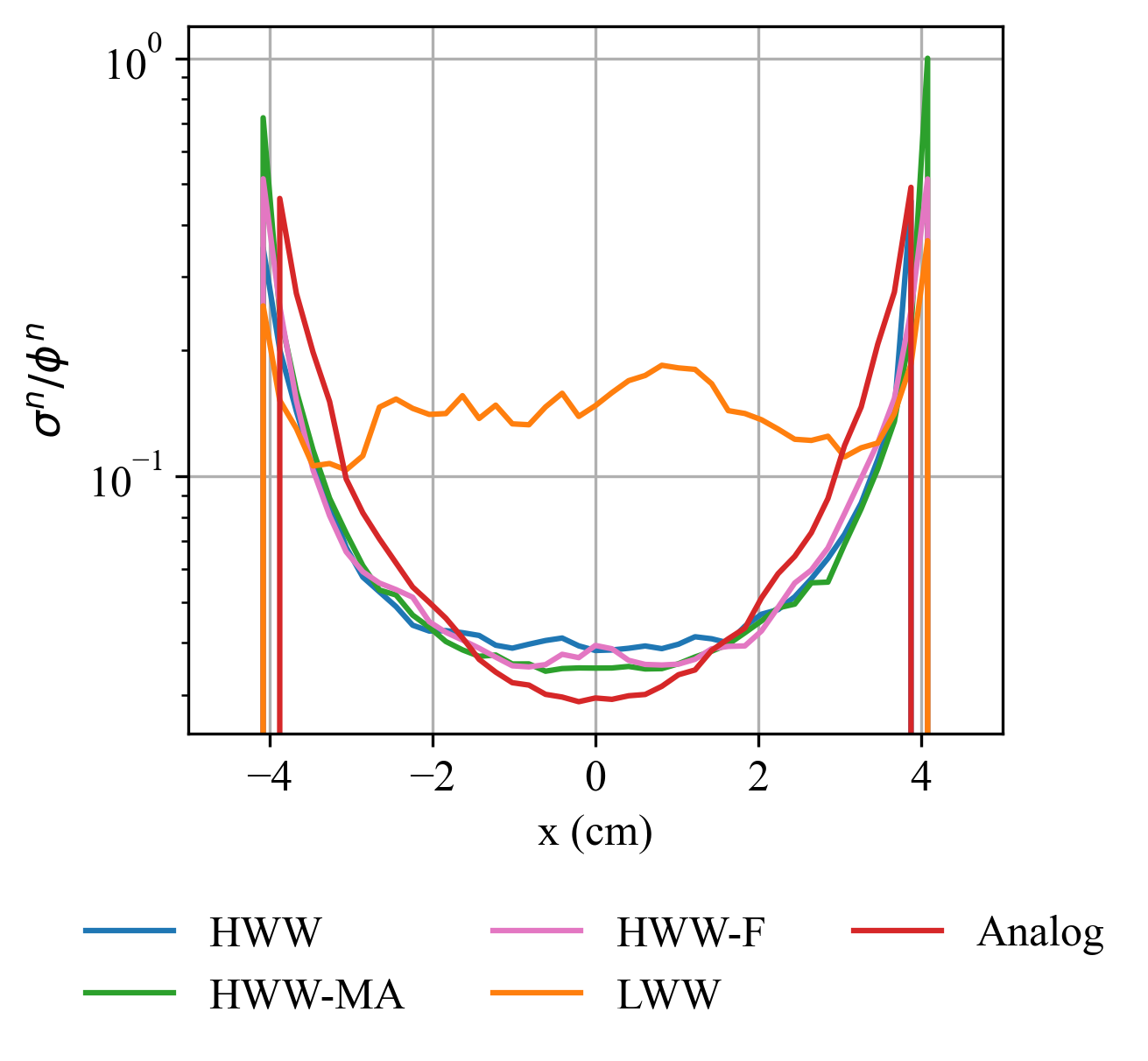}
             \vspace{-0.8cm}
        \caption{$n=4$}
        \label{fig:problem1_sdev_4}
    \end{subfigure}
    \hfill
    \begin{subfigure}[b]{0.45\textwidth}
        \centering
        \includegraphics[width=\textwidth]{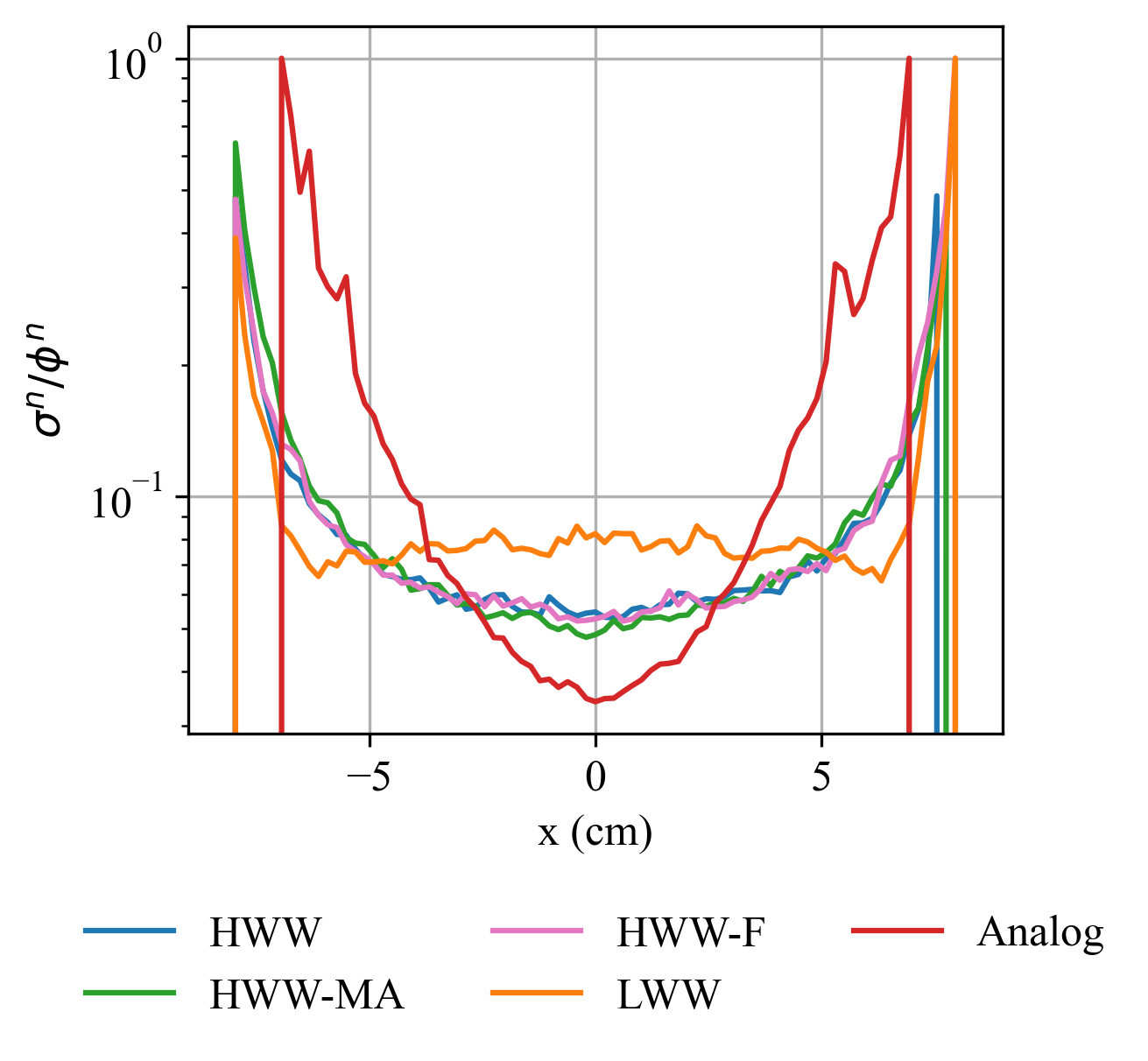}
             \vspace{-0.8cm}
        \caption{$n=8$}
        \label{fig:problem1_sdev_8}
    \end{subfigure}
    \hfill
    \begin{subfigure}[b]{0.45\textwidth}
        \centering
        \includegraphics[width=\textwidth]{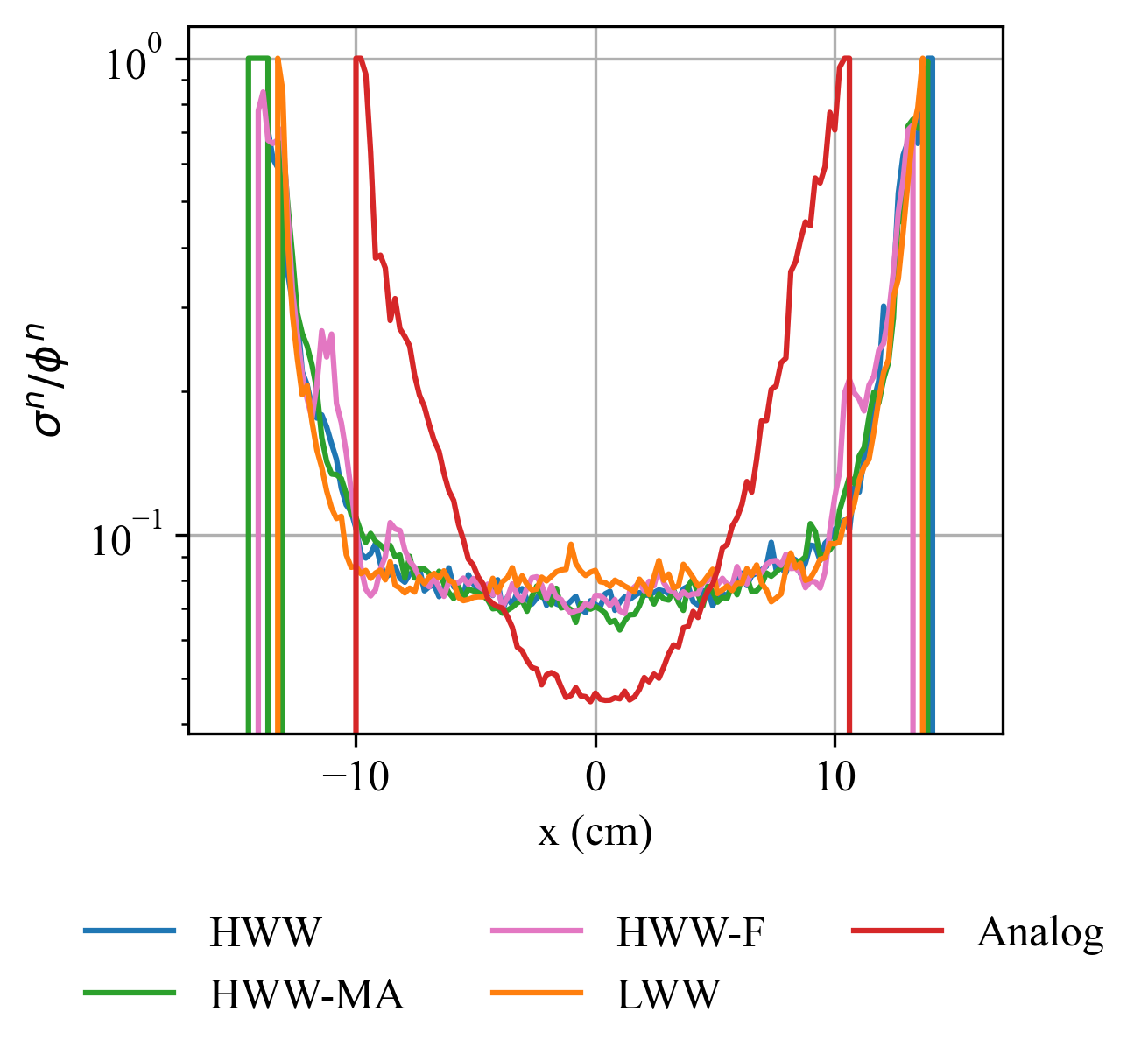}
             \vspace{-0.8cm}
        \caption{$n=16$}
        \label{fig:problem1_sdev_16}
    \end{subfigure}
         \vspace{-0.2cm}
    \caption{The relative standard deviation  $\frac{\sigma_i^n}{\phi_i^n}$ versus spatial position at  time steps $n=4,8,16$
        }
    \label{fig:Prob_1_rel_sdev}
\vspace{0.2cm}
    \centering
    \includegraphics[width=0.45\textwidth]{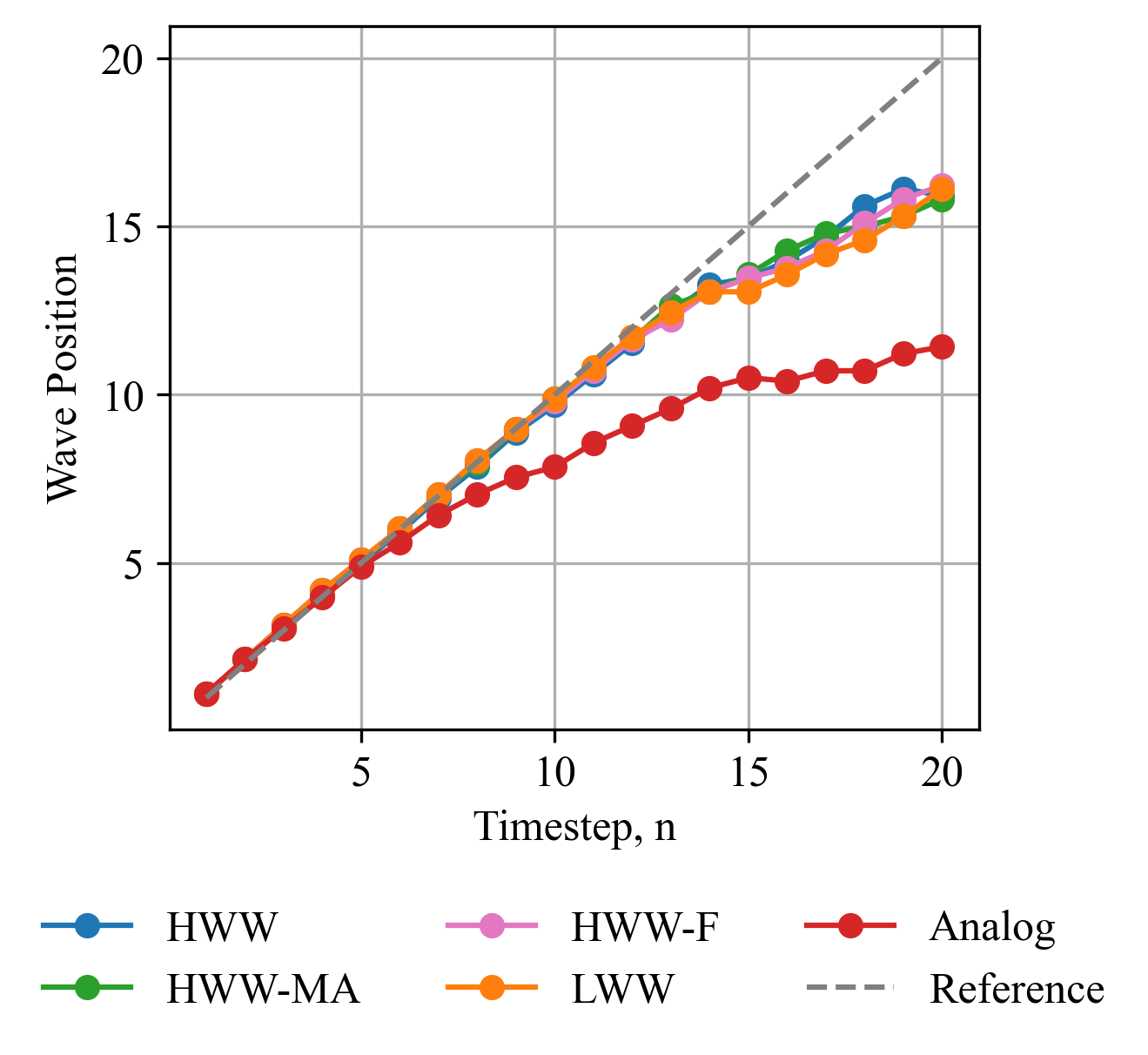}
             \vspace{-0.6cm}
   \caption{ \label{fig:wave_position}
   Furthest cell with particles versus time for the HWW algorithm with and without filtering techniques, the LWW and analog algorithms
 } 
\end{figure}

Figure \ref{fig:Prob_1_phi} shows  three different Monte Carlo solutions $\langle\phi\rangle^n$  versus spatial position at  time steps $n=4,8,16$. 
No filtering techniques were used by the HWW algorithm.
The effects of weight windows on Monte Carlo particle distribution are illustrated in Figure \ref{fig:Prob_1_tracks}, which 
 presents the number of particle tracks per source particle used by each method on these time steps.
The relative standard deviation  $\frac{\sigma_i^n}{\phi_i^n}$ of the Monte Carlo solutions versus spatial position  is presented in Figure \ref{fig:Prob_1_rel_sdev}.
 At $n=4$, the particle tracks  of the LWW algorithm  are very non-uniform, leading to a higher stochastic error in the center.  The lagged solution under-predicts the wave front position, causing more splitting and increased particle track density in this region.   By $n=8$, when the solution is changing less, the effect persists but to a lesser degree.  The HWW algorithm does not generate  as many particles at the wave front as the lagged weight windows, but particles still get much further than 
 those of the analog Monte Carlo. The solution of the HWW algorithm has
 the most uniform relative standard deviation.
 The variance of the analog Monte Carlo solution
 is lowest where the most physical particles are, and the opposite appears true for the solutions
 of the LWW algorithm.
 Here  the HWW algorithm achieves
 a balance, that leads to a high quality solution everywhere. 
 The solutions presented in Figure \ref{fig:Prob_1_phi} have similar wave front positions for all methods at $n=4,8$. At $n=16$, the solutions are similar for $|x|\leq10$. Near the wave front, the HWW and LWW solutions are more accurate than the analog solution. 
Figure \ref{fig:wave_position} shows the wavefront position for each of the methods studied as a function of time. This is defined as the outer edge of furthest cell with particles in it. This figure shows that analog Monte Carlo is unable to accurately track the wavefront beyond $5s$, and WW are essential to achieve an accurate wavefront position at later times. The point where the WW algorithms depart from the reference is determined by the $\varepsilon_{min}$ parameter.

\begin{figure}[H]
    \centering
    \includegraphics[width=0.6\textwidth]{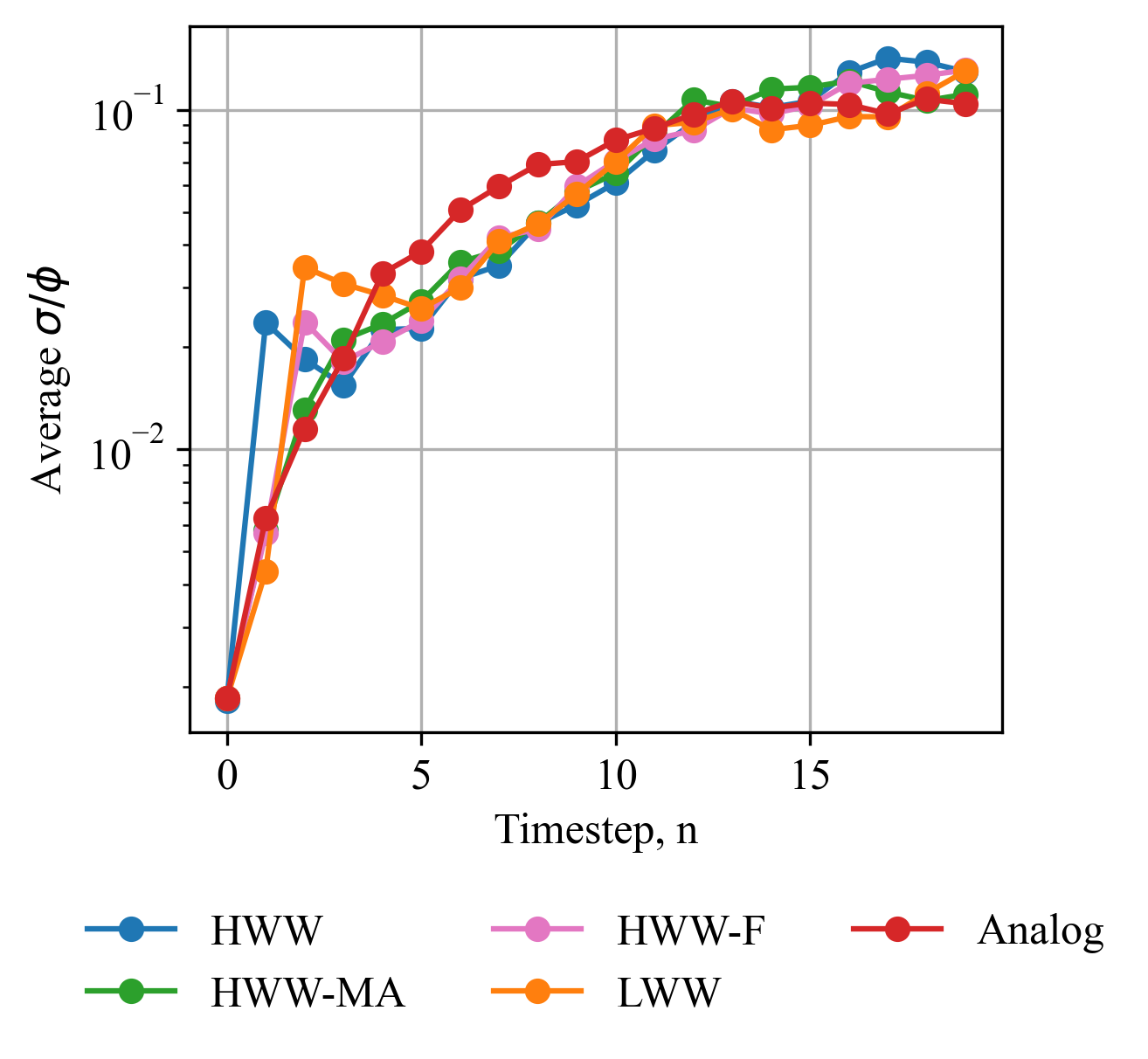} 
   \caption{ \label{fig:relative_sdev}
   Average relative standard deviation for the HWW algorithms with and without filtering techniques, the LWW and analog algorithms
 } 
\end{figure}

\begin{figure}[H]
    \centering
    \includegraphics[width=0.6\textwidth]{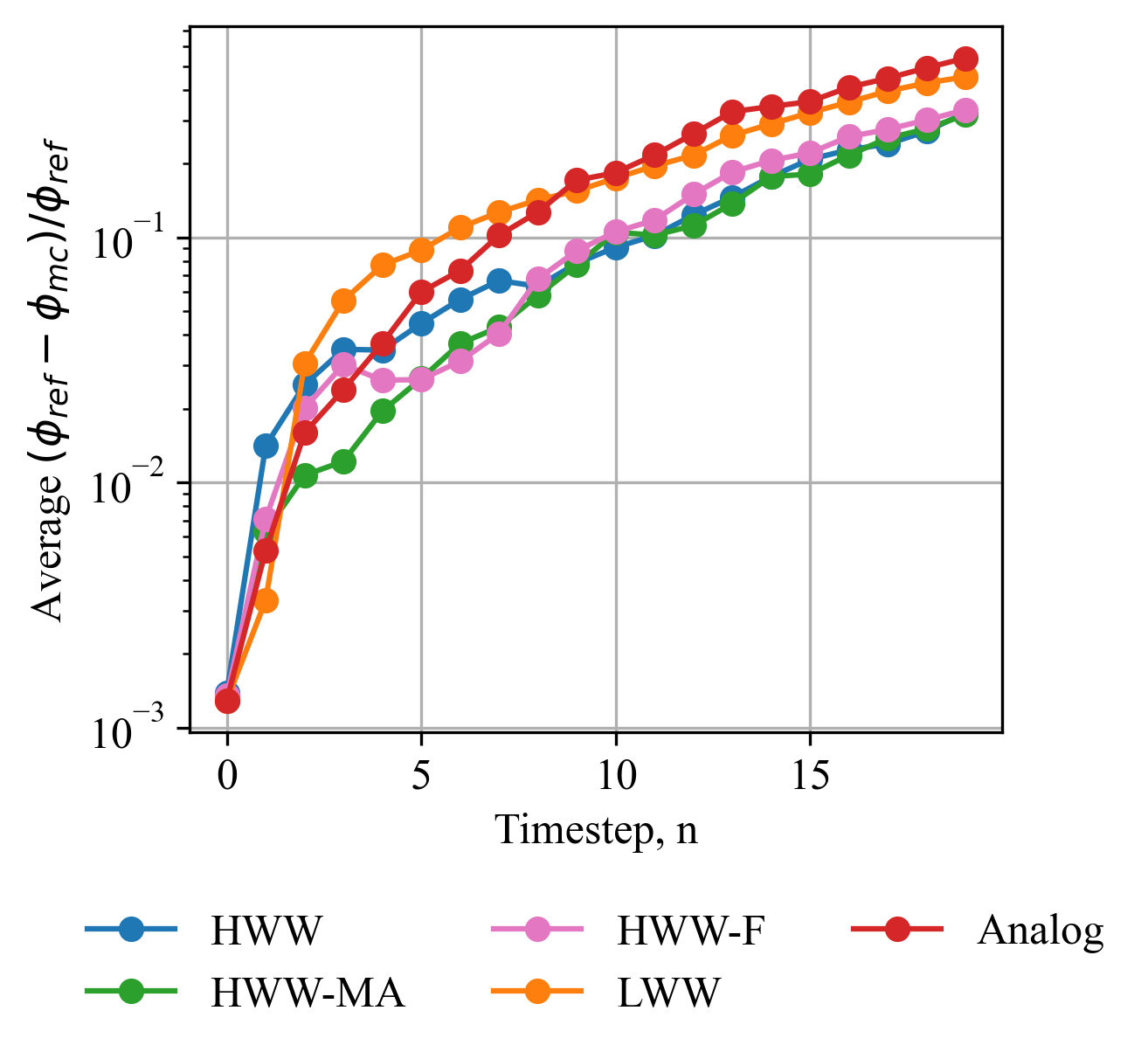} 
   \caption{ \label{fig:relative_error}
   Average relative error for the HWW algorithms with and without filtering techniques, the LWW and analog algorithms
 } 
\end{figure}

To quantify the accuracy of  Monte Carlo solutions, we use the relative standard deviation and the relative error.
The average of cell-wise relative standard deviation is defined as
\begin{equation}
    \overline{RSD}^n =\frac{1}{I}\sum_{i=1}^{I}\frac{\sigma_i^n}{\phi_{i,mc}^n}.
\end{equation}
This measure does not account for regions where the reference solution is nonzero and Monte Carlo particles have not reached. The average cell-wise relative error, defined as,
\begin{equation}
    \overline{RE}^n=\frac{1}{I}\sum_{i=1}^{I}\frac{|\phi_{i,mc}^n-\phi_{i,ref}^n|}{\phi_{i,ref}^n}
\end{equation}
is an explicit measure of the solution quality. A spatial average of both measures is computed on each timestep and plotted in Figures \ref{fig:relative_sdev} and \ref{fig:relative_error}. In Figure \ref{fig:relative_sdev} the positive effects of weight windows are seen between $n=4$ and $n=11$, where the more uniform relative standard deviation improves the average. Figure \ref{fig:relative_error} shows that the error caused by particles not making it to the edges of the problem is significant and dominates the analog solution. Also in Figure \ref{fig:relative_error} we can see the LWW solution has increased error caused by the uneven redistribution of particles during early timesteps. The HWW achieve the lowest relative error after $n=4$, because of the more uniform particle distribution and more particles reaching the wavefront.

To measure the efficiency of each algorithm, we define two figures of merit (FOM). One FOM is computed using relative standard deviation and one using relative error defined, respectively, as follows:
\begin{equation}
    FOM_{\sigma}^n = \frac{1}{\tau^n\left[\frac{1}{I}\sum_{i=1}^{I}\left(\frac{\sigma_i^n}{\phi_{i,mc}^n}\right)^2\right]} \, ,
    \label{eqn:FOM_sigma}
\end{equation}
\begin{equation}
    FOM_{Error}^n = \frac{1}{\tau^n\left[\frac{1}{I}\sum_{i=1}^{I}\left(\frac{|\phi_{i,mc}^n-\phi_{i,ref}^n|}{\phi_{i,ref}^n}\right)^2\right]} \, ,
    \label{eqn:FOM_error}
\end{equation}
where $\tau^n$ is the runtime of the $n^{th}$ census.

\begin{figure}[H]
    \centering
    \includegraphics[width=0.6\textwidth]{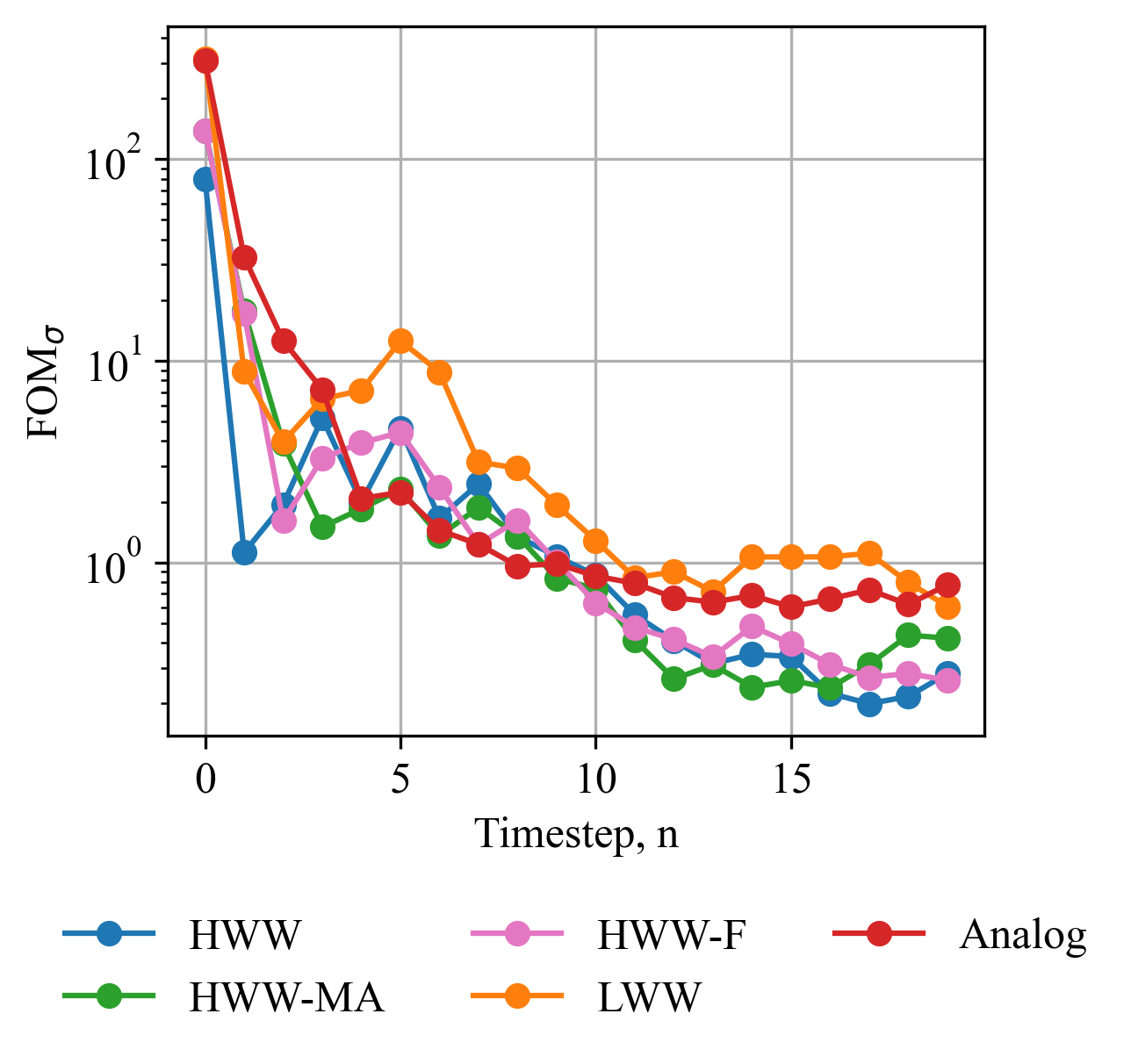} 
   \caption{ \label{fig:FOM_sdev}
   $FOM_{\sigma}$ for the HWW algorithm with and without filtering techniques, the LWW and analog algorithms
 } 
\vspace{0.75cm}
    \centering
    \includegraphics[width=0.6\textwidth]{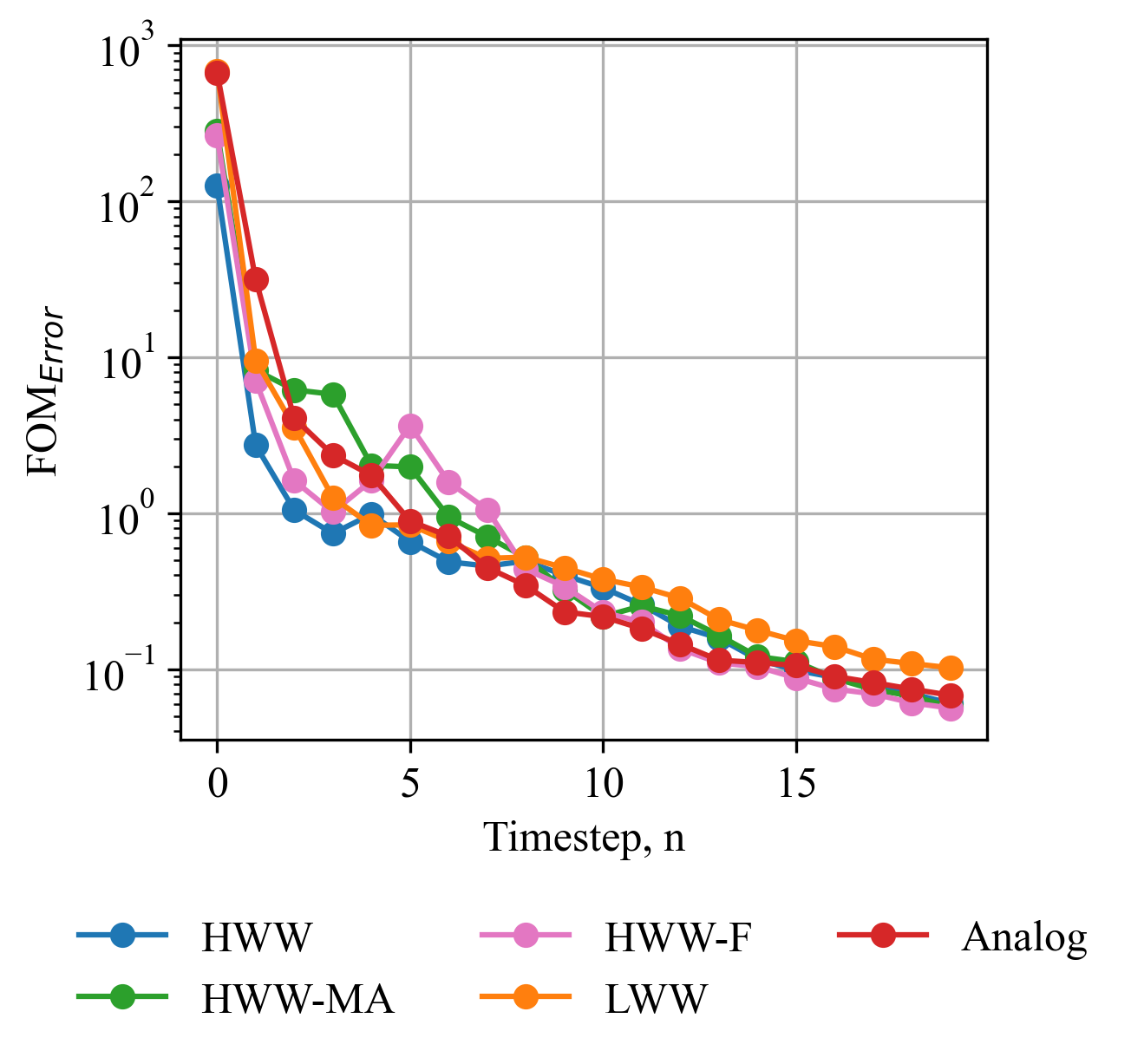} 
   \caption{ \label{fig:FOM_error}
   $FOM_{Error}$ for the HWW algorithm with and without filtering techniques, the LWW and analog algorithms
 } 
\end{figure}

\begin{figure}[H]
    \centering
    \includegraphics[width=0.6\textwidth]{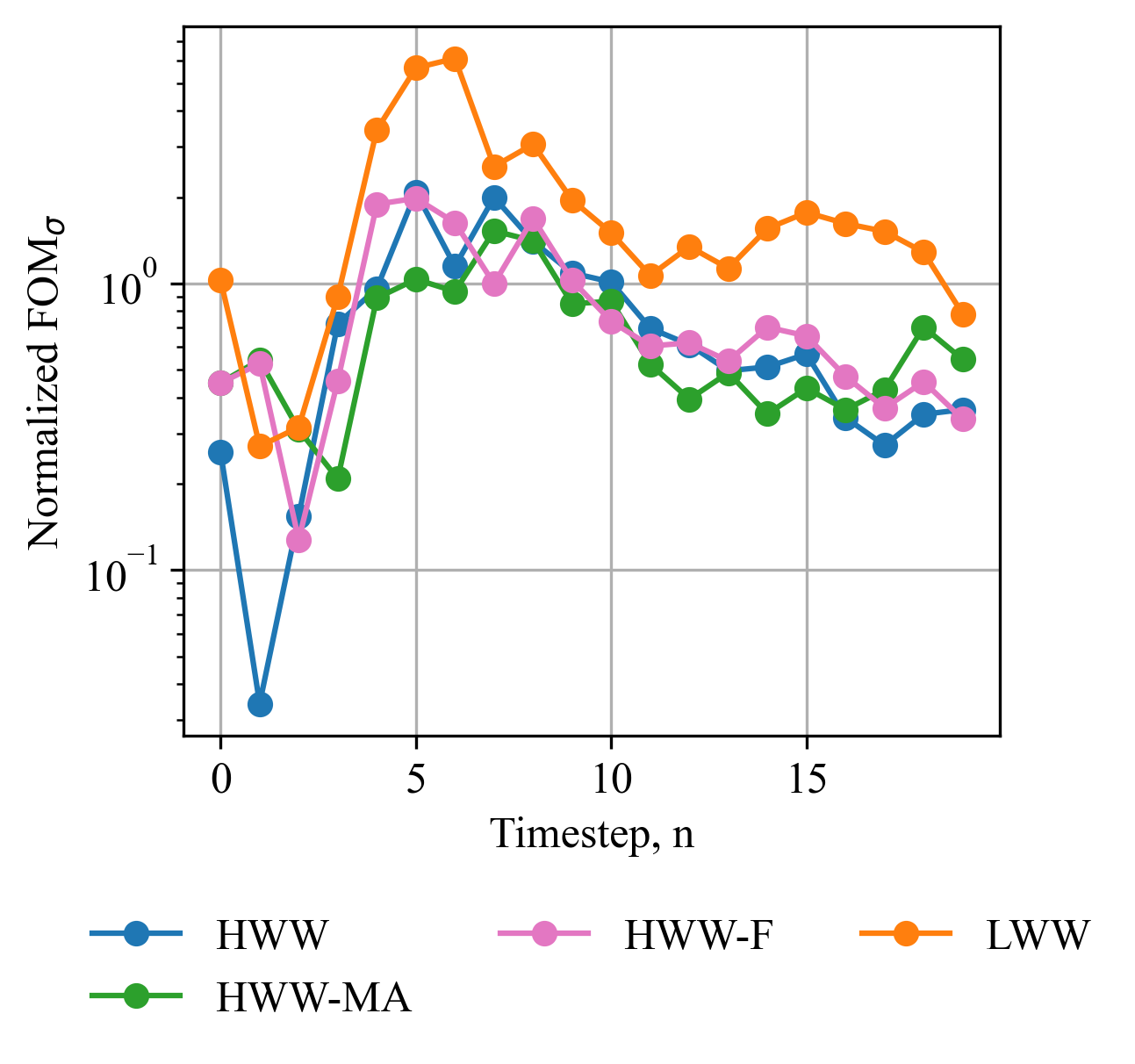} 
   \caption{ \label{fig:FOM_sdev_normalized}
   Normalized $FOM_{\sigma}$ for the HWW algorithm with and without filtering techniques and the LWW algorithm
 } 
\vspace{0.75cm}
    \centering
    \includegraphics[width=0.6\textwidth]{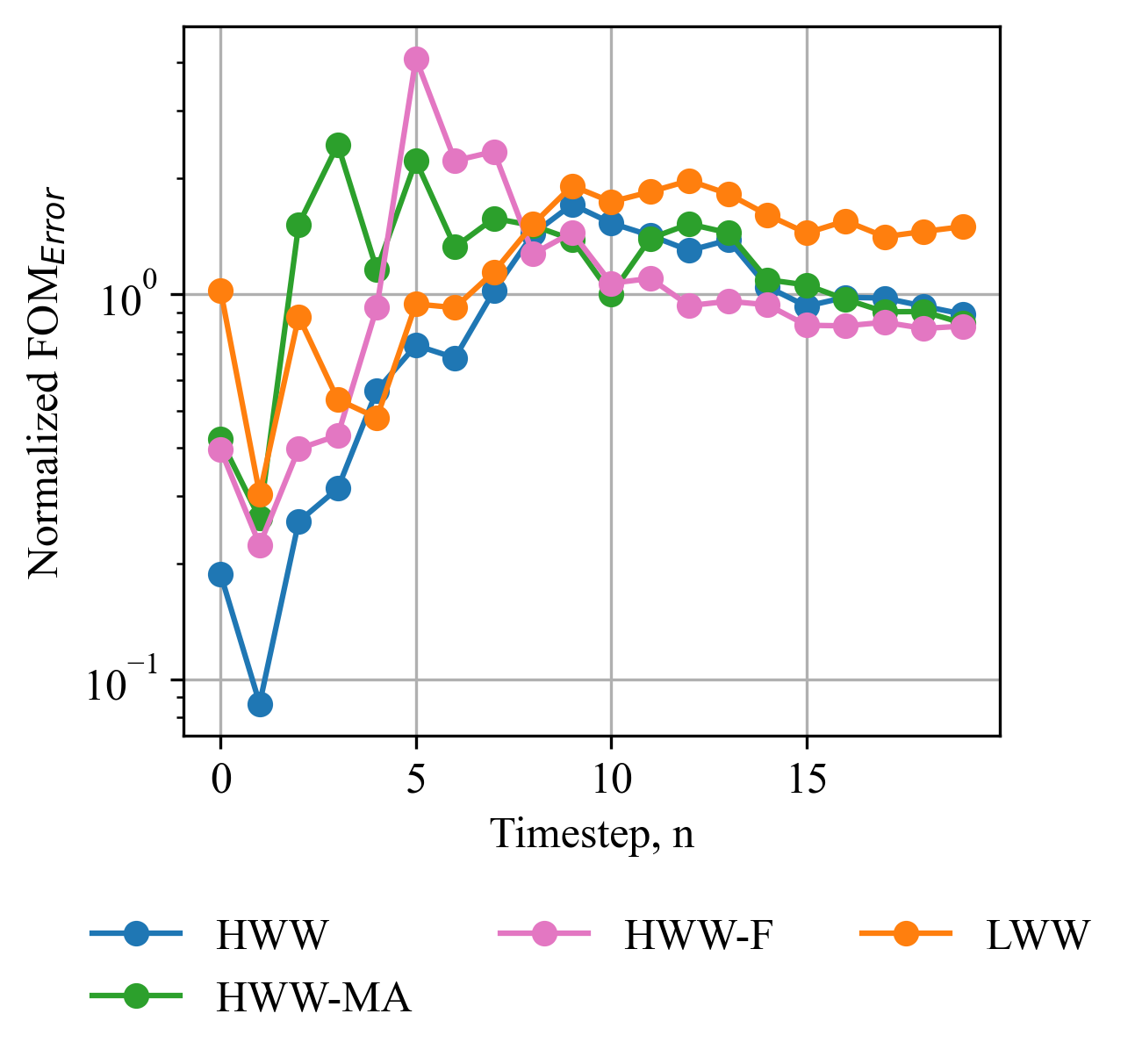} 
   \caption{ \label{fig:FOM_error_normalized}
   Normalized $FOM_{Error}$ for the HWW algorithm with and without filtering techniques, the LWW algorithm
 } 
\end{figure}

In Figures \ref{fig:FOM_sdev} and \ref{fig:FOM_error} are the FOMs for the HWW algorithms with and without filtering, LWW and analog Monte Carlo. The LWW and analog methods have shorter runtimes, since they require less tallies. This causes their $FOM_{\sigma}$ to be greater than the hybrid methods. 
In Figures \ref{fig:FOM_sdev_normalized} and \ref{fig:FOM_error_normalized} the FOMs for HWW and LWW are normalized by the analog FOM. Figure \ref{fig:FOM_error_normalized} shows that weight window methods generally improve $FOM_{Error}$ and the filtered hybrid methods improve $FOM_{Error}$ between $n=2$ and $n=15$. The decrease in $FOM_{Error}$ in the early time steps comes from splitting that results in additional particles being simulated, increasing runtime. During this period of time, the analog particle distribution is sufficient to capture the solution accurately, so weight windows decrease efficiency. The HWW algorithm with moving average filter achieves the greatest $FOM_{Error}$ increase of the hybrid methods, with an average of 1.25x improvement over analog. HWW with Fourier filtering, no filtering, and LWW achieved averages of 1.15x, 0.92x, 1.30x respectively.

\section{Conclusion}
\label{sec:discussion}
We have developed non-analog Monte Carlo methods for time-dependent transport problems with weight windows for global variance reduction. Weight windows are defined with the solution of a hybrid Monte Carlo / deterministic scheme based on the 
low-order second-moment equations approximated with the second-order discretization in time and space. One of the methods applies filtering techniques to closures and time-step initial conditions of the low-order equations to reduce the effect of stochastic noise on the auxiliary solution.

A parametric study of HWW showed that choice of $\rho$ and $\varepsilon_{min}$ has significant impact on particle splitting. If $\varepsilon_{min}$ is too small, excessive splitting will occur at the wavefront. Using too many updates to the hybrid solution, $u_{ww}^n$, did lead to some instability.  These results informed the choice of parameters used in the FOM calculations. The optimal choice of parameters is problem dependent and largely influenced by the accuracy of the auxiliary solution.

Filtering the Monte Carlo quantities of the low order problem improved the hybrid solution. Both Fourier and moving average filtering reduced noise and preserved the underlying structure of the closure. The moving average filter proved be effective and is computationally cheap. 

We have shown that constraining particles to an accurate hybrid auxiliary solution (HWW) provides a more uniform relative variance than analog Monte Carlo or LWW. The hybrid solution computed with unfiltered quantities had relative L2 error that was smaller than the Monte Carlo solution from the same histories. When the filtering techniques were applied, the hybrid solution error decreased further. Using HWW led to a more uniform spatial distribution of particles, and therefore relative stochastic error. In problems where the solution is required at both high and low flux regions, HWW enable efficient solution in the low flux regions. The relative error based figure of merit for HWW showed improvement over analog Monte Carlo. We automate the process of defining a deterministic problem and providing the solution to the Monte Carlo solver. By forming a new hybrid problem at each Monte Carlo census, accumulated discretization error is eliminated. This also allows for integration with multiphysics codes, where material properties and sources are evaluated at the census and not known beforehand.

This method is limited by the additional cost of the hybrid problem. A large portion of this cost is the tallies needed to define the low order problem. These tallies could be made cheaper or the hybrid problem could be formulated to require less additional tallies. Monotonization techniques can be applied to modify the second-order scheme for discretization of the LOSM equations \cite{duraisamy-2003,pg-dya-jcp-2020}
The developed algorithms can be extended to multidimensional problems to further realize the benefits of the self-adjoint second moment operator and filtering.

\section{Acknowledgments}
This work was supported by the Center for Exascale Monte-Carlo Neutron Transport (CEMeNT) a PSAAP-III project funded by   
National Nuclear Security Administration of Department of Energy, grant number: DE-NA003967.

\section{Authorship Contribution Statement}
Caleb Shaw:   conceptualization, methodology, formal analysis, investigation, software, data curation,  validation, visualization, writing - original draft, and writing - review \& editing.
Dmitriy Anistratov:  conceptualization, methodology, formal analysis, investigation, project administration, funding acquisition, resources, supervision, and writing - review \& editing.

\bibliographystyle{elsarticle-num}
\bibliography{sas-dya-ww-arxiv-v2}

\begin{thebibliography}{10}
\expandafter\ifx\csname url\endcsname\relax
  \def\url#1{\texttt{#1}}\fi
\expandafter\ifx\csname urlprefix\endcsname\relax\def\urlprefix{URL }\fi
\expandafter\ifx\csname href\endcsname\relax
  \def\href#1#2{#2} \def\path#1{#1}\fi

\bibitem{Book_MonteCarlo_1991_Lux_MonteCarloParticleTransportMethodsNeutronandPhotonCalculations}
I.~Lux, L.~Koblinger, Monte Carlo Particle Transport Methods: Neutron and
  Photon Calculations, CRC Press, 1991.

\bibitem{TechReport_2024_LANL_LA-UR-24-24602Rev.1_KuleszaAdamsEtAl}
J.~A. Kulesza, T.~R. Adams, J.~C. Armstrong, S.~R. Bolding, F.~B. Brown, J.~S.
  Bull, T.~P. Burke, A.~R. Clark, R.~A. Forster, III, J.~F. Giron, A.~S.
  Grieve, C.~J. Josey, R.~L. Martz, G.~W. McKinney, E.~J. Pearson, M.~E.
  Rising, C.~J. Solomon, Jr., S.~Swaminarayan, T.~J. Trahan, C.~A. Weaver,
  S.~C. Wilson, A.~J. Zukaitis, {MCNP\textsuperscript{\textregistered} Code
  Version 6.3.1 Theory \& User Manual}, Tech. Rep. LA-UR-24-24602, Rev.~1, Los
  Alamos National Laboratory, Los Alamos, NM, USA (May 2024).

\bibitem{osti_527548}
K.~A. {Van Riper}, T.~J. Urbatsch, P.~D. Soran, K.~Parsons, J.~E. Morel, G.~W.
  McKinney, S.~R. Lee, L.~A. Crotzer, R.~Alcouffe, F.~W. Brinkley, T.~E. Booth,
  J.~Anderson, {AVATAR} -- automatic variance reduction in {M}onte {C}arlo
  calculations, Tech. Rep. LA-UR-97-0919 (05 1997).

\bibitem{wagner-haghighat}
J.~C. Wagner, A.~Haghighat, Automated variance reduction of {M}onte {C}arlo
  shielding calculations using the discrete ordinates adjoint function, Nuclear
  Science and Engineering 128 (1998) 186--208.

\bibitem{wagner-peplow-mosher-2014}
J.~Wagner, D.~Peplow, S.~Mosher, Fw-cadis method for global and regional
  variance reduction of {M}onte {C}arlo radiation transport calculations,
  Nuclear Science and Engineering 176 (2014) 37--57.

\bibitem{Cooper_Larsen_2001}
M.~A. Cooper, E.~W. Larsen, Automated weight windows for global {M}onte {C}arlo
  particle transport calculations, Nuclear Science and Engineering 137~(1)
  (2001) 1–13.

\bibitem{Cooper_diss_1999}
M.~A. Cooper, An automated variance reduction method for global {M}onte {C}arlo
  neutral particle transport problems, Ph.D. thesis, University of Michigan
  (1999).

\bibitem{gol'din-cmmp-1964}
V.~Y. Gol'din, A quasi-diffusion method of solving the kinetic equation, Comp.
  Math. and Math. Phys. 4 (1964) 136--149.

\bibitem{auer-mihalas-1970}
L.~H. Auer, D.~Mihalas, On the use of variable {E}ddington factors in non-{LTE}
  stellar atmospheres computations, Monthly Notices of the Royal Astronomical
  Society 149 (1970) 65--74.

\bibitem{Wollaber_2008}
A.~B. Wollaber, Advanced {M}onte {C}arlo methods for thermal radiation
  transport, Ph.D. thesis, University of Michigan (2008).

\bibitem{mcclarren2013temperature}
R.~G. McClarren, T.~J. Urbatsch, Temperature-extrapolation method for implicit
  {M}onte {C}arlo - radiation hydrodynamics calculations, in: Proceedings of
  the International Conference on Mathematics and Computational Methods Applied
  to Nuclear Science \& Engineering (M\&C 2013), American Nuclear Society, Sun
  Valley, Idaho, USA, 2013, pp. on CD--ROM.

\bibitem{pozulp-mc2023}
M.~M. Pozulp, T.~S. Haut, P.~S. Brantley, J.~L. Vujic, An implicit {M}onte
  {C}arlo acceleration scheme, in: Proceedings of The International Conference
  on Mathematics and Computational Methods Applied to Nuclear Science and
  Engineering, Niagara Falls, Canada, August 13-17, 2023.

\bibitem{pozulp-mc2025}
M.~M. Pozulp, T.~S. Haut, P.~S. Brantley, S.~Olivier, J.~L. Vujic, A hybrid
  second moment method for thermal radiative transfer, in: Proceedings of the
  International Conference on Mathematics and Computational Methods Applied to
  Nuclear Science and Engineering, Denver, Co, April 27-30, 2025, pp.
  1518--1527.

\bibitem{smm-1976}
E.~E. Lewis, W.~F. {Miller Jr}, Comparison of {P$_1$} synthetic acceleration
  techniques, Transactions of the American Nuclear Society 23 (6 1976).

\bibitem{vnn-dya-ans-annual-2024}
V.~N. Novellino, D.~Y. Anistratov, Analysis of hybrid {MC}/deterministic
  methods for transport problems based on low-order equations discretized by
  finite volume schemes, Transactions of the American Nuclear Society 130
  (2024) 408--411.

\bibitem{vnn-dya-arxiv-2025}
V.~N. Novellino, D.~Y. Anistratov, Multi-level hybrid {M}onte {C}arlo /
  deterministic methods for particle transport problems, arXiv:2510.09545v1
  (2025).

\bibitem{Joss}
J.~P. Morgan, I.~Variansyah, S.~L. Pasmann, K.~B. Clements, B.~Cuneo, A.~Mote,
  C.~Goodman, C.~Shaw, J.~Northrop, R.~Pankaj, E.~Lame, B.~Whewell, R.~G.
  McClarren, T.~S. Palmer, L.~Chen, D.~Y. Anistratov, C.~T. Kelley, C.~J.
  Palmer, K.~E. Niemeyer, {M}onte {C}arlo / {D}ynamic {C}ode ({MC/DC}): An
  accelerated {P}ython package for fully transient neutron transport and rapid
  methods development, Journal of Open Source Software 9~(96) (2024) 6415.

\bibitem{Landman_McClarren_Madsen_Long_2014}
J.~T. Landman, R.~G. McClarren, J.~R. Madsen, A.~R. Long, Analysis of lagged
  weight windows for implicit {M}onte {C}arlo variance reduction, Transactions
  of the American Nuclear Society 111~(1) (2014) 648–650.

\bibitem{Cefus-Larsen-ttsp-1989}
G.~R. Cefus, E.~W. Larsen, Stability analysis of the quasidiffusion and second
  moment methods for iteratively solving discrete-ordinates problems, Transport
  Theory and Statistical Physics 18~(5-6) (1989) 493--511.

\bibitem{Smith}
S.~W. Smith, The Scientist and Engineer’s Guide to Digital Signal Processing,
  California Technical Publishing, 1997.

\bibitem{azurv1}
B.~Ganapol, R.~Baker, J.~A. Dahl, R.~E. Alcouffe, Homogeneous infinite media
  time-dependent analytical benchmarks, Technical Report LA-UR-01-1854, Los
  Alamos National Laboratory (January 2001).

\bibitem{osti_1889957}
J.~A. Kulesza, T.~R. Adams, J.~C. Armstrong, S.~R. Bolding, F.~B. Brown, J.~S.
  Bull, T.~P. Burke, A.~R. Clark, R.~A.~F. III, J.~F. Giron, T.~S. Grieve,
  C.~J. Josey, R.~L. Martz, G.~W. McKinney, E.~J. Pearson, M.~E. Rising,
  C.~J.~S. Jr., S.~Swaminarayan, T.~J. Trahan, S.~C. Wilson, A.~J. Zukaitis,
  {MCNP} code version 6.3.0 theory and user manual, Tech. rep., Los Alamos
  National Laboratory (LANL), Los Alamos, NM (United States) (09 2022).

\bibitem{duraisamy-2003}
K.~Duraisamy, J.~D. Baeder, J.-G. Liu, Concepts and application of
  time-limiters to high resolution schemes, J. of Scientific Computing 19
  (2003) 139--162.

\bibitem{pg-dya-jcp-2020}
P.~Ghassemi, D.~Y. Anistratov, Multilevel quasidiffusion method with
  mixed-order time discretization for multigroup thermal radiative transfer
  problems, Journal of Computational Physics 409 (2020) 109315.

\end{thebibliography}
\end{document}